\documentclass[11pt,preprint]{imsart}
\RequirePackage[OT1]{fontenc}
\RequirePackage{amsthm,amsmath}
\RequirePackage[colorlinks,citecolor=blue,urlcolor=blue]{hyperref}
\RequirePackage{hypernat}
\usepackage{graphicx}
\usepackage{enumerate}
\usepackage{latexsym}
\usepackage{amssymb}
\numberwithin{equation}{section}

\def\bbR{\mathbb{R}}

\usepackage{algpseudocode}
\usepackage[ruled, slide, lined, linesnumbered]{algorithm2e}

\def\bbN{\mathbb{N}}

\startlocaldefs \theoremstyle{plain}
\newtheorem{Theorem}{Theorem}[section]

\newtheorem{Proposition}[Theorem]{Proposition}

\newtheorem{lemma}{\indent Lemma}
\newtheorem{Remark}{Remark}[section]
\theoremstyle{definition}

\newtheorem{proposition}{\indent Proposition}

\makeatletter\def\theequation{\arabic{section}.\arabic{equation}}
\@addtoreset{equation}{section}\makeatother

\endlocaldefs
\usepackage{multirow}
\begin{document}
\begin{frontmatter}

\title{Nonparametric Bayesian Aggregation for Massive
Data}\runtitle{Nonparametric Bayesian Aggregation}
\author
{
    Zuofeng Shang\thanks{Associate Professor, Department of Mathematical Sciences, New Jersey Institute of Technology, Newark, NJ 07102. E-mail: zuofengshang@gmail.com. This work was supported by NSF grants DMS 1821157 and DMS 1764280.},~
	Botao Hao\thanks{Ph.D student, Department of Statistics, Purdue University, West Lafayette, IN 47906. E-mail: haobotao000@gmail.com.},~
	Guang Cheng\thanks{Professor, Department of Statistics, Purdue University, West Lafayette, IN 47906. E-mail: chengg@purdue.edu. This work was supported by the National Science Foundation under Grants DMS-1712907, DMS-1811812, DMS-1821183, and Office of Naval Research, (ONR N00014-18-2759).}
}

\today

\runauthor{Shang et al.}

\begin{abstract}
We develop a set of scalable Bayesian inference procedures for a general class of nonparametric regression models. Specifically, nonparametric Bayesian inferences are separately performed on each subset randomly split from a massive dataset, and then the obtained local results are aggregated into global counterparts. This aggregation step is explicit without involving any additional computation cost. By a careful partition, we show that our aggregated inference results obtain an oracle rule in the sense that they are equivalent to those obtained directly from the entire data (which are computationally prohibitive). For example, an aggregated credible ball achieves desirable credibility level and also frequentist coverage while possessing the same radius as the oracle ball.
\end{abstract}



\end{frontmatter}

\section{Introduction}

With rapid development in modern technology, massive data sets are becoming more and more common.
An important feature of massive data is their large volume which hinders applications of traditional statistical methods. For example, due to huge data amount and limited CPU memory, it is often impossible to process the entire data in a single machine. In the parallel computing environment, a common practice is to distribute massive data to multiple processors, and then aggregate local results in an efficient way. A series of frequentist methods such as \cite{KTSJ12,MHM10,ZDW15,ZCL14} have been proposed in this Divide-and-Conquer (D\&C) framework.

In Bayesian community, there are quite a few computational or methodological works developed for massive data such as scalable algorithms for Bayesian variable selection (\cite{vRD15,WPD14}) and scalable posterior sampling in parametric models (\cite{WD14,WGHD15}). 
Theoretical guarantees of D\&C methods
have been recently obtained in robust estimation (\cite{MSLD14}), approximation of posterior distributions (\cite{LSD17}),
posterior interval estimation (\cite{SLD18}), credible sets of signal in Gaussian white noise (\cite{Sv17,Sv18}).
Rather, the present paper puts focus on uncertainty quantification of the
model parameter in general nonparametric regression, primarily in theoretical aspects. 
For instance, how to aggregate individual 
posterior means into a global one that maintains frequentist optimality?
How to aggregate individual credible balls into a global one with a minimal possible radius?
And how many divisions and what kind of priors should be chosen to guarantee Bayesian and frequentist validity of the aggregated ball? 
We attempt to address these questions in a univariate nonparametric regression setup. 

Specifically, we develop a set of aggregation procedures in Bayesian nonparametric regression. 
As a first step, nonparametric Bayesian regression is separately fitted based on each subsample randomly split from a massive dataset. 
A variety of finite sample valid credible balls (credible intervals) for regression functions (their linear functionals \cite{RR12}, e.g., local values) are then constructed from each individual posterior distribution based on MCMC. In the second step, we aggregate these credible balls (credible intervals) into global counterparts analytically without involving any additional computation. For example, the center of an aggregated ball is obtained by weighted averaging Fourier coefficients of all individual (approximate) posterior modes, while the radius is given through an explicit formula on individual radii. A notable advantage of this distributed strategy is its dramatically faster computational speed, and this computational advantage becomes more obvious as data size grows.

Our aggregation procedures are proven to obtain an oracle rule in the sense that they are equivalent to those obtained directly from the entire data, i.e., called as oracle results which are computationally prohibitive in practice. For example, 
our aggregated posterior means are proven to achieve optimal estimation rate, and
our aggregated credible ball achieves desirable credibility level and also frequentist coverage while possessing asymptotically
the same radius as the oracle ball. These oracle results hold when the assigned Gaussian process priors in each subset are properly chosen and the number of subsets does not grow too fast. A fundamental theory underlying Bayesian aggregation is a {\em uniform} version of nonparametric Gaussian approximation theorem, also called as Bernstein-von Mises theorem. Developed based on our recent work \cite{SC14}, this theory states  that a sequence of individual posterior distributions converge to Gaussian processes uniformly over the number of subsets.

The rest of this paper is organized as follows. 
Section \ref{sec:for:nonparametric:model} describes our Bayesian nonparametric model
with a Gaussian process prior, based on which our main results are developed in Section~\ref{sec:main}. Specifically, a uniform nonparametric Gaussian approximation theorem is established in Section \ref{sec:UBvM}, and all the Bayesian aggregation procedures together with their theoretical guarantee are provided in Sections~\ref{sec:aggr:post:mean}--\ref{sec:asymp:post:infer}. 
Section \ref{sec:simulations} provides a simulation study to justify our methods. Section~\ref{sec:real} applies the proposed procedures to a real dataset of large size. Main proofs are provided in Appendix. Other results and additional plots are given in a supplementary document \cite{SCBigRate}.

\section{Nonparametric Bayesian Aggregation: An Illustration}\label{sec:gar}

In this section, we provide a concrete example to demonstrate the intuition of our nonparametric Bayesian aggregation procedure.
Our example is based on the special uniform design and periodic Sobolev space 
which makes our aggregation procedure explicit and easy to understand.
Section \ref{sec:toy:np:bayes:model} describes our nonparametric Bayesian model,
and Section \ref{sec:toy:agg:proc} demonstrates our algorithm and its numeric performance.
General aggregation procedures will be proposed
in Sections \ref{sec:for:nonparametric:model} and \ref{sec:main} with asymptotic properties investigated as well.

\subsection{Nonparametric Bayesian model}\label{sec:toy:np:bayes:model}
Suppose that we observe the data $Z_i=(Y_i,X_i)$, $i=1,\ldots,N$,
generated from the following Gaussian regression model with uniform design
\begin{equation}\label{toy:example:model}
Y_i|f, X_i\sim N(f(X_i),1),\,\,\,\,X_1,\ldots,X_N\overset{iid}{\sim} Unif[0,1].
\end{equation}
Randomly split $\{1,2,\ldots,N\}$ 
into $s$ subsets $I_1,I_2,\ldots,I_s$ with $|I_1|=\cdots=|I_s|=n$ (so $N=ns$). 
Denote $\textbf{D}_j=\{Z_i| i\in I_j\}$ 
the $j$-th subsample for $j=1,\ldots,s$
and $\textbf{D}=\cup_{j=1}^s\textbf{D}_j$ the entire sample. 

Suppose that $f$ belongs to an $m$-order periodic Sobolev space $S_0^m[0,1]$ 
where
$S_0^m[0,1]$ is the collection of all functions on $[0,1]$ of the form
\begin{equation}\label{per:spline:expression}
f(x)=\sqrt{2}
\sum_{k=1}^\infty f_k\cos(2\pi kx)+\sqrt{2}\sum_{k=1}^\infty g_k\sin(2\pi kx)
\end{equation}
with real coefficients $f_k,g_k$ satisfying
\begin{equation}\label{J:penalty}
\sum_{k=1}^{\infty}(f_k^2+g_k^2)(2\pi k)^{2m}<\infty.
\end{equation}
Here, $m>1/2$ is a constant describing the smoothness of the functions.
Wahba (1990) \cite{W90} introduced a Gaussian process (GP) prior on $f$ which has an
interesting smoothing spline interpretation. 
Specifically, she assumed that the coefficients $f_k,g_k$ in (\ref{per:spline:expression}) 
are independent and normally distributed as follows:
\begin{equation}\label{GP:prior}
f_k,g_k\sim N\left(0,[(2\pi k)^{2m+\beta}+n\lambda(2\pi k)^{2m}]^{-1}\right),\,\,\,\,k=1,2,\ldots,
\end{equation}
where $\beta>1$ and $\lambda\ge0$ are predecided constants. 
In particular, $\beta$ represents the ``relative smoothness'' of the prior to the parameter space
and $\lambda$ represents the amount of rescaling. Rescaling priors are also considered by \cite{Sv17,Sv18} 
for constructing credible sets of signals in Gaussian white noise.
It can be examined that if $f$ satisfies (\ref{per:spline:expression}) and
(\ref{GP:prior}), then $f$
is a Gaussian process with mean zero and isotropic covariance function
\begin{equation}\label{eqn:kernal}
K_0(x,x')=2\sum_{k=1}^\infty\frac{\cos(2\pi k(x-x'))}{(2\pi k)^{2m+\beta}+n\lambda(2\pi k)^{2m}},\,\,\,\,x,x'\in[0,1].
\end{equation}
Wahba (\cite{W90}) showed that the above GP prior (\ref{GP:prior}) 
generates a posterior distribution corresponding to 
a penalized likelihood function (with $\lambda$ the penalty parameter). 
This provides a Bayesian interpretation for smoothing splines. 
Below we provide some details to justify this argument.

Let $\Pi_\lambda$ denote the probability distribution of $f$ under (\ref{GP:prior}).
To derive the posterior distribution, we need to find the ``prior density'' of $f$.
Unlike the parametric settings where the prior densities are Radon-Nikodym (RN) derivatives 
w.r.t. Lebesgue measure, in the current infinite-dimensional setting
it is impossible to do so since there is no Lebesgue measure on $S^m_0[0,1]$
(see \cite{BTJ92}).
Instead, we need to characterize the prior density of $f$ as an RN derivative w.r.t. other
kinds of measures such as Gaussian measure.
Following Wahba (\cite{W90}), $\Pi_\lambda$ and $\Pi\equiv\Pi_0$ (corresponding to $\lambda=0$) are equivalent probability measures, and the RN derivative of $\Pi_\lambda$ w.r.t. $\Pi$ is
\begin{eqnarray}\label{expression:prior}
\frac{d\Pi_\lambda}{d\Pi}(f)&=&\prod_{k=1}^\infty\left(1+n\lambda(2\pi k)^{-\beta}\right)^{-1}
\times\exp\left(-\frac{n\lambda}{2}\sum_{k=1}^\infty(f_k^2+g_k^2)(2\pi k)^{2m}\right)\nonumber\\
&=&\prod_{k=1}^\infty\left(1+n\lambda(2\pi k)^{-\beta}\right)^{-1}
\times\exp\left(-\frac{n\lambda}{2}\int_0^1f^{(m)}(x)^2dx\right)\nonumber\\
&=&\prod_{k=1}^\infty\left(1+n\lambda(2\pi k)^{-\beta}\right)^{-1}
\times\exp\left(-\frac{n\lambda}{2}J(f)\right),
\end{eqnarray}
where $J(f)=\int_0^1f^{(m)}(x)^2dx$.
Note that $\prod_{k=1}^\infty\left(1+n\lambda(2\pi k)^{-\beta}\right)^{-1}$ converges thanks to $\beta>1$
so that (\ref{expression:prior}) is a valid expression.
(\ref{expression:prior}) provides an expression for the prior density of $f$,
which induces the following posterior distribution for $f$ given subsample $j$:
\begin{eqnarray}\label{toy:pen:lik}
dP(f|\textbf{D}_j)
&\propto& P(\textbf{D}_j|f)d\Pi_\lambda(f)\nonumber\\
&\propto&\exp\left(-\frac{1}{2}\sum_{i\in I_j}(Y_i-f(X_i))^2-\frac{n\lambda}{2}J(f)\right)d\Pi(f),\,\,
j=1,\ldots,s.
\end{eqnarray}
Recall that $I_j$ indexes the $j$-th subsample.
The right hand-side of (\ref{toy:pen:lik}) corresponds to
penalized likelihood function $\ell_j(f)=-\frac{1}{2n}\sum_{i\in I_j}(Y_i-f(X_i))^2-\frac{\lambda}{2}J(f)$
which has been well studied in smoothing spline literature (\cite{W90}).
Theoretically, we recommend to choose $\lambda\asymp N^{-\frac{2m}{2m+\beta}}$ 
which will be proven to yield optimal Bayesian inference; see Sections~\ref{sec:for:nonparametric:model} and \ref{sec:main}. 
The duality between the posterior and smoothing spline, i.e., (\ref{toy:pen:lik}),
enables us to easily choose $\lambda$ for practical use, e.g., GCV considered by \cite{W90}.  

\subsection{Nonparamtric Bayesian Aggregation}\label{sec:toy:agg:proc}

First of all, we calculate 
$\breve{f}_{j,n}=E\{f|\textbf{D}_j\}$, $j=1,\ldots,s$,
the posterior means based on individual posterior distributions \eqref{toy:pen:lik}. 
Then we construct a $(1-\alpha)$-th credible ball centering at $\breve{f}_{j,n}$ with radius $r_{j,n}(\alpha)$. That is, $r_{j,n}(\alpha)>0$
such that $P(f\in S_0^m[0,1]:\|f-\breve{f}_{j,n}\|_{L^2}\le r_{j,n}(\alpha)|\textbf{D}_j)=1-\alpha$,
where $\|\cdot\|_{L^2}$ is the usual $L^2$-norm, i.e., $\|f\|_{L^2}=\sqrt{\int_0^1f(x)^2dx}$.
In practice, $\breve{f}_{j,n}$ and $r_{j,n}(\alpha)$ can be both estimated by 
the posterior samples. For instance, 
generate $M$ independent samples $f_{j1},\ldots, f_{jM}$ from (\ref{toy:pen:lik});
estimate $\breve{f}_{j,n}$ by their average and
estimate $r_{j,n}(\alpha)$ by the $(1-\alpha)$-th percentile of
$\|f_{jl}-\breve{f}_{j,n}\|_{L^2}$ for $1\le l\le M$. We postpone the computational details of 
the sampling procedure to Section \ref{subsection:computational_details}.

We next present a concrete aggregation scheme (procedures (\ref{alg:i})--(\ref{alg:iii}) below) to construct a credible ball based on these individual results $\{\breve{f}_{j,n}, r_{j,n}(\alpha)\}_{j=1}^s$. 
Specifically,
an aggregated credible ball for $f$, denoted $R_N(\alpha)$, is constructed with its center/radius 
obtained through \textit{weighted} averaging the individual centers/radii. 
Unlike simple averaging commonly used in frequentist setting (see \cite{ZDW15}), our procedures for posterior mean aggregation
and radius aggregation are weighted averaging with weights $w_{s,N,\lambda,k}$ defined in (\ref{special:aggregation:mean}). These weights are used to calibrate the prior effect such that
the aggregation procedure can have satisfactory asymptotic property.
The details of our procedure are demonstrated as follows:
\begin{enumerate}[(1).]
	\item\label{alg:i} \textit{Posterior mean aggregation}. For $j$-th subsample and $k\ge 1$, find 
	\begin{equation}\label{eqn:coeff_intergral}
	\breve{f}_{j,n,k}=\sqrt{2}\int_0^1 \breve{f}_{j,n}(x)\cos(2\pi kx)dx, \ \breve{g}_{j,n,k}=\sqrt{2}\int_0^1\breve{f}_{j,n}(x)\sin(2\pi kx)dx,
	\end{equation}
	where $\breve{f}_{j,n}$ is the posterior mean based on subsample $j$. Then we aggregate these quantities through the following formulas:
	\begin{equation}
	\breve{f}_{N,\lambda,k}=\sum_{j=1}^s\breve{f}_{j,n,k}/s, \ \breve{g}_{N,\lambda,k}=\sum_{j=1}^s\breve{g}_{j,n,k}/s.
	\end{equation}
	In the end, we let
\begin{eqnarray}\label{special:aggregation:mean}
\breve{f}_{N,\lambda}(x)=\sum_{k=1}^\infty w_{s,N,\lambda,k}
\left\{\breve{f}_{N,\lambda,k}\sqrt{2}\cos(2\pi kx)
+\breve{g}_{N,\lambda,k}\sqrt{2}\sin(2\pi kx)\right\},
\end{eqnarray} 
where
$w_{s,N,\lambda,k}=\frac{s(2\pi k)^{2m+\beta}+N(1+\lambda(2\pi k)^{2m})}{(2\pi k)^{2m+\beta}+N(1+\lambda(2\pi k)^{2m})}$ for $k\ge1$.
	\item\label{alg:ii} \textit{Posterior radius aggregation}. Aggregate the radii $r_{j,n}(\alpha)$ through the following formula:
	\begin{equation}\label{special:aggregation:radius}
	r_N(\alpha)=\sqrt{A_{N,s}\left(\frac{1}{s}\sum_{j=1}^sr_{j,n}(\alpha)^2\right)+B_{N,s}},
	\end{equation}
	where
		\begin{equation}\label{eqn:radius_weight}
	\begin{split}
	&A_{N,s}=\sqrt{C_2/D_2}s^{-\frac{4m+2\beta-1}{2(2m+\beta)}},\\
	&B_{N,s}=\left(2C_1-
	2D_1\sqrt{C_2/D_2}s^{-\frac{1}{2(2m+\beta)}}\right)N^{-\frac{2m+\beta-1}{2m+\beta}},\\
	&C_k=\int_0^\infty(1+(2\pi x)^{2m}+(2\pi x)^{2m+\beta})^{-k}dx,\,\,\,\,k=1,2,\\
	&D_k=\int_0^\infty(1+(2\pi x)^{2m})^{-k}dx,\,\,\,\,k=1,2.
	\end{split}
	\end{equation} 
	\item\label{alg:iii} \textit{Aggregated credible ball}:
	 \begin{equation}\label{special:aggregate:set}
	 R_N(\alpha)=\{f\in S_0^m[0,1]:\|f-\breve{f}_{N,\lambda}\|_{L^2}\le r_N(\alpha)\}.
	 \end{equation}
\end{enumerate}

Algorithms based on weighted averaging have been proposed in numerous computational aspects. For instance, 
\cite{huang2005sampling, Neiswanger2014, Scott2016} 
proposed computational procedures for efficiently aggregating local MCMC samples in which the aggregation steps involve proper weight averaging.
Such algorithms are particularly useful to produce MCMC samples from the oracle posterior 
which can be used for various inferential purposes, e.g., estimation and testing.
The present paper focuses on inferences, e.g., construction of credible balls, in a special class of nonparametric regression models, and has more extensive theoretical guarantees. 

In practice, one can approximate the integral \eqref{eqn:coeff_intergral} through discretization; 
see Section \ref{subsection:computational_details}. Theorem \ref{cp:cr:strong} will show that $R_N(\alpha)$ given in (\ref{special:aggregate:set}) 
asymptotically covers $1-\alpha$ mass of the posterior based on the full data set and includes the true function with probability approaching one. More theoretical study on $R_N(\alpha)$ such as its center and radius can be found in Sections~\ref{sec:aggr:post:mean} and \ref{sec:cr:strong}. 
Note that these sections present an aggregation procedure in a more general context, which covers (\ref{special:aggregate:set}) as a special case.

\begin{figure}[htp]
	\begin{center}
		\includegraphics[scale=0.42]{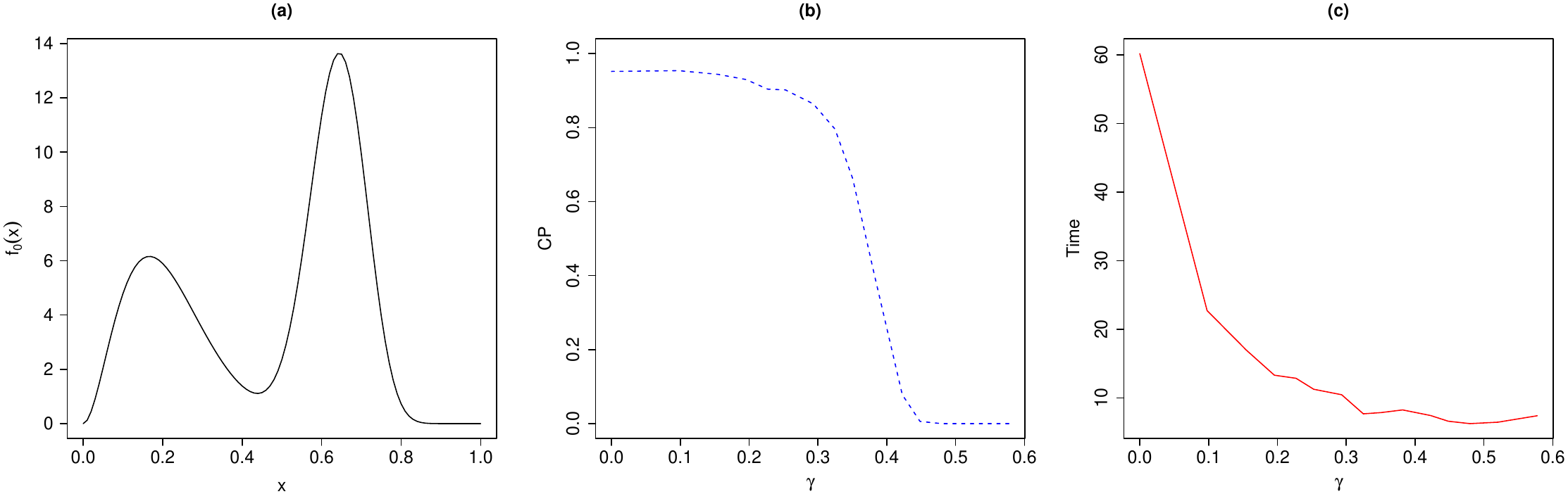}
		\caption{\textit{\footnotesize Examination of our aggregation procedures (\ref{alg:i})--(\ref{alg:iii}). Results are based on
		$N=1200$ observations generated from (\ref{toy:example:model}) and a GP prior (\ref{GP:prior})
		with $m=\beta=2$ and $\lambda=N^{-2/3}$.
	 (a) True regression function $f_0(x)=2.4\beta_{30,17}(x)+1.6\beta_{3,11}(x)$,
        where $\beta_{a,b}$ is the probability density function for $Beta(a,b)$.
        (b) Coverage probability (CP) of $R_N(0.95)$ vs. $\gamma$.
        (c) Computing time (in seconds) of $R_N(0.95)$ vs. $\gamma$.}}
		\label{fig:toy:example}
	\end{center}
\end{figure}

A toy simulation study was carried out to examine the proposed procedures (\ref{alg:i})--(\ref{alg:iii}). 
Specifically, we examine the computing time and coverage probability (CP) of $R_N(\alpha)$ for
various choices of $s$. The CP is defined as the relative frequency of
the sets that cover the truth. We choose $m=\beta=2$ in our GP prior \eqref{GP:prior}. Results are summarized in Figure \ref{fig:toy:example}. 
Plot (a) displays the true function $f_0$ under which data were generated.
Plot (b) displays how the CP varies as $\gamma:=\log(s)/\log(N)$. 
Plot (c) displays that the computing time decreases when $\gamma$ increases. There seems to be a transition for CP vs. $\gamma$, i.e.,
CP is uniformly close to one when $0\le\gamma<0.3$ and approaches zero when $\gamma>0.4$. 
In conclusion, $R_N(\alpha)$ possesses
both satisfactory frequentist coverage and computational efficiency
when $\gamma\approx 0.2$. Other choices of $\gamma$ either lower CP
or slow down the computing. Thus, under a proper choice of $s$, our aggregation procedure can
maintain good statistical properties and reduce computing burden at the same time. 
Careful readers may have noticed that the CP approaches one rather than the credibility level $(1-\alpha)$. This issue can be addressed by a modified aggregated set proposed in Section \ref{sec:cr:weak}. More comprehensive simulation results
are provided in Section \ref{sec:simulations} to examine various aggregation procedures such as the pointwise
credible intervals.

\section{A Nonparametric Bayesian Framework Based on General Design and Space}\label{sec:for:nonparametric:model}

In this section, we introduce a more general Bayesian nonparametric framework based on
general design and function space
under which the aggregation results will be obtained.
Suppose that the data $\{Y_i,X_i\}_{i=1}^N$ follow 
a nonparametric regression model:
\begin{equation}\label{basic:model}
Y_i|f,X_i\overset{ind.}{\sim} N(f(X_i),\sigma^2),\,\,\,\,X_1,\ldots,X_N\overset{iid}{\sim}\pi(x),
\end{equation}
where $\pi(\cdot)$ is a probability density on $\mathbb{I}=(0,1)$,
and $f$ belongs to an $m$-order Sobolev space $S^m(\mathbb{I})$:
\begin{eqnarray}\label{def:sob}
S^m(\mathbb{I})=\{f\in L^2(\mathbb{I})| f^{(0)},f^{(1)},
\ldots,f^{(m-1)} \textrm{are abs. cont. and}\;f^{(m)}\in L^2(\mathbb{I})\}.
\end{eqnarray}
In particular, $S^m_0[0,1]$ is a proper subset of $S^m(\mathbb{I})$.
Throughout, we let $m>1/2$ such that $S^m(\mathbb{I})$ is a reproducing kernel Hilbert space (RKHS).
For technical convenience, assume $\sigma^2=1$ 
and $0<\inf_{x\in\mathbb{I}}\pi(x)\le\sup_{x\in\mathbb{I}}\pi(x)<\infty$. 
When $\sigma^2$ is unknown, our approach can still be applied with $\sigma^2$
replaced by its consistent estimate.

For any $f,g\in S^m(\mathbb{I})$, define
$V(f,g)=E\{f(X)g(X)\}$ and $J(f,g)=\int_0^1 f^{(m)}(x)g^{(m)}(x)dx$.
Following \cite{SC13}, there exists a sequence of eigenfunctions 
$\varphi_1,\varphi_2,\ldots\in S^m(\mathbb{I})$ and a sequence of eigenvalues
$0=\rho_1=\rho_2=\cdots=\rho_m<\rho_{m+1}\le\rho_{m+2}\le\cdots$ such that $\rho_\nu\asymp \nu^{2m}$ and
\begin{equation}
V(\varphi_\nu,\varphi_\mu)=\delta_{\nu\mu},\,\,\,\,J(\varphi_\nu,\varphi_\mu)=\rho_\nu\delta_{\nu\mu},\,\,\,\,
\nu,\mu\ge1,
\end{equation}
where $\delta_{\nu\mu}$ is the Kronecker's delta.
 
We next place a prior distribution $\Pi_\lambda$ on $f$,
where $\Pi_\lambda$ is a probability measure on $S^m(\mathbb{I})$
and $\lambda\ge0$ is a hyperparameter. 
Similar to Section \ref{sec:gar},
we will characterize $\Pi_\lambda$ through its Radon-Nikodym (RN) derivative 
w.r.t. $\Pi$, with $\Pi$ a pre-given probability measure $\Pi$ on $S^m(\mathbb{I})$.
Specifically, assume that the RN derivative of $\Pi_\lambda$ w.r.t. $\Pi$ satisfies
\begin{equation}\label{RN:derivative}
\frac{d\Pi_\lambda}{d\Pi}(f)\propto\exp\left(-\frac{n\lambda}{2}J(f)\right),
\end{equation}
where $J(f)$ is defined in (\ref{J:penalty}).
Interestingly, it is possible to explicitly construct $\Pi_\lambda$ and $\Pi$ such that
(\ref{RN:derivative}) holds. To see this, let 
\begin{equation}\label{prior1}
G_\lambda(\cdot)=\sum_{\nu=m+1}^\infty w_\nu\varphi_\nu(\cdot),
\end{equation}
where $w_\nu$'s are independent of the observations satisfying
$w_\nu\sim N(0,1/(\rho_\nu^{1+\beta/(2m)}+n\lambda\rho_\nu)), \nu>m$.
Let $G(\cdot)=G_{\lambda=0}(\cdot)$. 
Suppose $\Pi_\lambda$ and $\Pi$ are probability measures 
induced by $G_\lambda$ and $G$,
i.e., $\Pi_\lambda(S)=P(G_\lambda\in S)$ and $\Pi(S)=P(G\in S)$
for any measurable $S\subseteq S^m(\mathbb{I})$.
It follows by H\'{a}jek's lemma (see \cite{SC14}) that (\ref{RN:derivative}) holds.
In (\ref{prior1}), $\lambda\ge0$ and $\beta>1$ are both hyper-parameters characterizing the smoothness of the prior.
It is easy to check that the sample path of $G_\lambda$ belongs to $S^m(\mathbb{I})$ for any $\beta>1$ almost surely. 
As demonstrated in a simulation study, 
the GCV-selected $\lambda$ is sufficient to provide satisfactory results.

\section{Main Results}\label{sec:main}

In this section, we present a series of main results that are built upon a uniform Gaussian approximation theorem (Section~\ref{sec:UBvM}). Three classes of aggregation procedures are then proposed: aggregated credible balls in both strong and weak topology, and aggregated credible intervals for linear functionals. These results can be classified into two types: {\em finite sample} construction (Sections \ref{sec:cr:strong}, \ref{sec:cr:weak} and \ref{sec:ci:functional}) and {\em asymptotic} construction (Section \ref{sec:asymp:post:infer}). The former construction is often time-consuming since its radius (interval length) is obtained through $s$ posterior sampling, while the latter employs a large-sample limit of the radius given by an explicit formula. The computational gain will be illustrated by the simulations in Section~\ref{sec:simulations}.
Similar to Section \ref{sec:gar},
let $I_1,I_2,\ldots,I_s$ be a random partition of $\{1,2,\ldots,N\}$ such that $\cup_{j=1}^s I_j=\{1,2,\ldots,N\}$ with $|I_j|=n$
for $j=1,\ldots,s$ and $N=ns$.

\subsection{A Uniform Gaussian Approximation Theorem}\label{sec:UBvM}

A fundamental theory underlying Bayesian aggregation is developed in this section. It is a {\em uniform} version of Gaussian approximation theorem that characterizes the limit shapes of a sequence of individual posterior distributions. This uniform validity holds if the number of posterior distributions does not grow too fast. Also, Bayesian aggregation procedures possess frequentist validity if $\lambda$ is chosen properly.

Similar to (\ref{toy:pen:lik}), we note that each sub-posterior distribution can be written as	
$$dP(f|\textbf{D}_j)\propto \exp(n\ell_{jn}(f))d\Pi(f),$$ where $\ell_{jn}(f)=n^{-1}\sum_{i\in I_j}(Y_i-f(X_i))^2-(\lambda/2)J(f)$. Define
\begin{equation}\label{smoothing:spline:likhood}
\widehat{f}_{j,n}=\arg\max_{f\in S^m(\mathbb{I})}\ell_{jn}(f), \,\,j=1,\ldots,s.
\end{equation}
Suppose that $\widehat{f}_{j,n}$ admits the following Fourier expansion:
\begin{eqnarray}\label{dfn:fjh}
\widehat{f}_{j,n}(\cdot)=\sum_{\nu=1}^\infty \widehat{f}_\nu^{(j)}\varphi_\nu(\cdot),\,\,1\le j\le s.
\end{eqnarray}

Define $h=\lambda^{1/(2m)}$ with $h^\ast:= N^{-\frac{1}{2m+\beta}}$. We remark that $h^\ast$ is an optimal choice for
our aggregation procedure as will be shown later.
\begin{Theorem}\label{uniform:bvm:thm}(Uniform Gaussian Approximation)
Suppose that $f_0$ admits a Fourier expansion $f_0(\cdot)=\sum_{\nu=1}^\infty f_\nu^0\varphi_\nu(\cdot)$
which further satisfies
	\begin{eqnarray*}
		\mbox{Condition (\textbf{S})}:\;\;\;\;\;\;\;\;\;
		\sum_{\nu=1}^\infty |f_\nu^0|^2\rho_\nu^{1+\frac{\beta-1}{2m}}<\infty
	\end{eqnarray*}
If the following holds 
\begin{eqnarray}\label{ratcon}
m>1+\frac{\sqrt{3}}{2}\approx 1.866, 1<\beta<2m+\frac{1}{2m}-1, s=o(N^{\frac{\beta-1}{2m+\beta}})\;\;\mbox{and}\;\;h\asymp h^\ast,
\end{eqnarray}
then we have as $N\rightarrow\infty$,
	\begin{equation}\label{uniform:bvm:expression}
	\sup_{S\in\mathcal{S}}\max_{1\le j\le s}|P(S|\textbf{D}_j)-P_{0j}(S)|=
	O_{P_{f_0}}\left(\sqrt{s}N^{-\frac{4m^2+2m\beta-10m+1}{4m(2m+\beta)}}(\log{N})^{\frac{5}{2}}\right),
	\end{equation}
	where $\mathcal S$ is the Borel $\sigma$-algebra on $S^m(\mathbb{I})$ with respect to $\Pi$, and $P_{0j}$'s are GPs defined by
	\begin{equation}\label{scalable:gaussian:measure}
	P_{0j}(S)=\frac{\int_S\exp\left(-\frac{n}{2}\|f-\widehat{f}_{j,n}\|^2\right)d\Pi(f)}
	{\int_{S^m(\mathbb{I})}\exp\left(-\frac{n}{2}\|f-\widehat{f}_{j,n}\|^2\right)d\Pi(f)},
	\,\,\,\,S\in\mathcal{S}.
	\end{equation}
\end{Theorem}
Proof of Theorem \ref{uniform:bvm:thm} is rooted in \cite{SC14} who essentially considered $s=1$.
Substantial efforts have been made here to quantify a range of partition size $s$ such that local posteriors can be uniformly approximated
by GPs. The explicit structure of the GPs provides a guideline for our aggregation procedures
which will be introduced in subsequent sections. It should be emphasized that our aggregation of GPs is
weighted-averaging which is different from product-based ones such as \cite{cao2014generalized}.

Condition (\textbf{S}) amounts to requiring known regularity of the truth $f_0\in S^{m+\frac{\beta-1}{2}}(\mathbb{I})$. 
This can be seen from the inequality $\sum_{\nu=1}^\infty |f_\nu^0|^2\nu^{2m+\beta-1}<\infty$ since $\rho_\nu\asymp\nu^{2m}$. 
This condition essentially means that $f_0$ has derivatives up to order $m+\frac{\beta-1}{2}$
(when this order is integer-valued).
Combined with (\ref{ratcon}) this means that the regularity of $f_0$ belongs to $(m,2m+\frac{1}{4m}-1)$,
i.e., the truth function is jointly confined by both functional space and the prior.
The $\|\cdot\|$-norm used in  (\ref{scalable:gaussian:measure}) is defined as follows. For any $g,\widetilde{g}\in S^m(\mathbb{I})$,
define
\begin{eqnarray}
\langle g,\widetilde{g}\rangle=V(g,\widetilde{g})+\lambda J(g,\widetilde{g})\label{dfn:inn}
\end{eqnarray}
and its squared norm $\|g\|^2=\langle g, g\rangle$. Clearly, $\langle\cdot,\cdot\rangle$ is a valid inner product on
$S^m(\mathbb{I})$.

\begin{Remark}
We remark that (\ref{ratcon}) can be replaced by a more general rate condition:
\[
		nh^{2m+1}\ge1,\,\,a_n=O(\widetilde{r}_n),\,\,
		b_n\le1,\,\,r_n^2b_n\le\widetilde{r}_n^2,\,\, n\widetilde{r}_n^2b_n=o(1),
\]
where $r_n=(nh)^{-1/2}+h^m, \widetilde{r}_n=(nh/\log{2s})^{-1/2}+h^{m+\frac{\beta-1}{2}}, 
a_n=n^{-1/2}h^{-\frac{6m-1}{4m}}r_n\log{N}, 
b_n=n^{-1/2}h^{-\frac{6m-1}{4m}}(\log{N})^{3/2}$.
Here, we provide a technical explanation for the terms $r_n, \widetilde{r}_n, a_n, b_n$.
Specifically, $r_n$ can be viewed as the rate of convergence of local ordinary penalized MLE (\ref{smoothing:spline:likhood}),
$\widetilde{r}_n$ can be viewed as the posterior contraction rate of the local Bayesian mode,
$a_n, b_n$ are error bounds of the higher-order remainders in the Taylor expansions of the individual penalized likelihood functions.
Uniform Gaussian approximation for general $h$ (not necessarily $h\asymp h^*$)
can be established under such condition.
\end{Remark}

Theorem 3.5 in \cite{SC14} shows that $P_{0j}$ (conditional on $\textbf{D}_j$) is induced by a Gaussian process, denoted as $W^j$, in the sense that $P_{0j}(S)=P(W^j\in S|\textbf{D}_j)$ for any $S\in\mathcal{S}$. Define
\begin{eqnarray}\label{gamnu}
\tau_\nu^2=\rho_\nu^{1+\frac{\beta}{2m}},\,\,\nu\ge1.
\end{eqnarray}
Then we have
\begin{eqnarray*}
	W^j(\cdot)&=&\sum_{\nu=1}^\infty (a_{n,\nu}\widehat{f}_\nu^{(j)}+b_{n,\nu}\tau_\nu v_\nu)\varphi_\nu(\cdot),
	\,\,j=1,2,\ldots,s,
\end{eqnarray*}
where $a_{n,\nu}=n(1+\lambda\rho_\nu)(\tau_\nu^2+n(1+\lambda\rho_\nu))^{-1}$, $b_{n,\nu}=(\tau_\nu^2+n(1+\lambda\rho_\nu))^{-1/2}$ and $v_\nu\sim N(0, \tau_\nu^{-2})$. For convenience, define the mean functions of $W^j$ as
\begin{eqnarray}\label{wjmean}
\widetilde{f}_{j,n}(\cdot):=\sum_{\nu=1}^\infty a_{n,\nu}\widehat{f}_\nu^{(j)}\varphi_\nu(\cdot),\,\,j=1,\ldots,s,
\end{eqnarray}
such that we can re-express $W^j$ as
\[
W^j=\widetilde{f}_{j,n}+W_n,\,\,j=1,\ldots,s,
\]
where $W_n(\cdot):=\sum_{\nu=1}^\infty b_{n,\nu}\tau_\nu v_\nu\varphi_\nu(\cdot)$ is a zero-mean GP. Note that the posterior mode $\widetilde{f}_{j,n}$ is very close to $\widehat f_{j,n}$ since $\|\widetilde{f}_{j,n}-\widehat{f}_{j,n}\|=o_{P_{f_0}}(1)$ uniformly for $1\le j\le s$; see the proof of Theorem \ref{cp:cr:strong}. The above characterization of $W^j$ is useful for the subsequent Bayesian aggregation procedures.

\subsection{Aggregated posterior means}\label{sec:aggr:post:mean}

In this section, we propose a method to aggregate the posterior means $\breve{f}_{j,n}:=E\{f|\textbf{D}_j\}$, for $j=1,\ldots,s$. The aggregated mean function, denoted as $\breve{f}_{N,\lambda}(\cdot)$, can be viewed as a nonparametric Bayesian estimate of $f$, and will be used to construct aggregated credible balls/intervals to be introduced later.

Our aggregation procedure is 
\begin{equation}\label{aggreg:center}
\breve{f}_{N,\lambda}(\cdot)=\sum_{\nu=1}^\infty\frac{a_{N,\nu}}{a_{n,\nu}}
V\left(\frac{1}{s}\sum_{j=1}^s\breve{f}_{j,n},\varphi_\nu\right)\varphi_\nu(\cdot).
\end{equation}
Note that when the model is Gaussian and $f\in S_0^m(0,1)$,
(\ref{aggreg:center}) becomes (\ref{special:aggregation:mean}).
Next we will show that the aggregation procedure
(\ref{aggreg:center}) yields minimax optimality in the following theorem.
\begin{Theorem}\label{uniform:post:mean}
Under conditions of Theorem \ref{uniform:bvm:thm}, the following result holds:
\begin{equation}
\max_{1\le j\le s}\|\breve{f}_{j,n}-\widetilde{f}_{j,n}\|
=O_{P_{f_0}}\left(\widetilde{r}_n\sqrt{s}N^{-\frac{4m^2+2m\beta-10m+1}{4m(2m+\beta)}}(\log{N})^{\frac{5}{2}}\right),
\label{uniform:post:mean:eqn}
\end{equation}
If, in addition, $3/2<\beta<2m+1/(2m)-3/2$ and
$s$ satisfies
\begin{equation}\label{s:for:aggregation}
s=o\left(N^{\frac{4m^2+2m\beta-11m+1}{8m(2m+\beta)}}(\log{N})^{-\frac{3}{2}}\right),
\end{equation}
then it holds that
\begin{equation}
\|\breve{f}_{N,\lambda}-f_0\|_2=O_{P_{f_0}}\left(N^{-\frac{2m+\beta-1}{2(2m+\beta)}}\right),
\label{big:rate:eqn}
\end{equation}
where $\|f\|_2=\sqrt{V(f)}$ denotes the $V$-norm.
\end{Theorem} 
According to \cite{VZ08:2}, the rate in (\ref{big:rate:eqn}) is minimax optimal
given Condition (\textbf{S}).

\subsection{Aggregated credible region in strong topology}\label{sec:cr:strong}

In this section, we construct an aggregated credible region based on $s$ individual credible regions (w.r.t. a weighted $\ell^2$-norm). Specifically, $s$ radii are combined in an explicit manner. This aggregated region possesses nominal posterior mass asymptotically, and is further proven to cover the true function with probability tending to one. This nice frequentist property is achieved as long as $s$ is not diverging fast and the assigned GP prior in each subset is chosen by setting $h\asymp h^\ast$, i.e., $\lambda\asymp N^{-2m/(2m+\beta)}$. The conservative frequentist coverage can be improved to the nominal level if we use a weaker norm in defining credible region; see Section~\ref{sec:cr:weak}.

Based on each subset $\textbf{D}_j$, the individual credible ball is constructed as follows: 
$$R_{j,n}(\alpha)=\{f\in S^m(\mathbb{I}): \|f-\breve{f}_{j,n}\|_2\le r_{j,n}(\alpha)\}.$$ The credible ball centers around the posterior mean $\breve{f}_{j,n}$, while its radius $r_{j,n}(\alpha)$ is directly sampled from MCMC such that $P(R_{j,n}(\alpha)|\textbf{D}_j)=1-\alpha$ for any $\alpha\in(0,1)$. We will construct an ``aggregated" region
centering at $\breve{f}_{N,\lambda}$ with radius explicitly constructed as follows:
\begin{eqnarray}\label{bar:r:N}
&&r_N(\alpha)=
\sqrt{
	\frac{1}{N}\left[\zeta_{1,N}+\sqrt{\frac{\zeta_{2,N}}{\zeta_{2,n}}}
	\left(\frac{n}{s}\sum_{j=1}^s r_{j,n}^2(\alpha)-\zeta_{1,n}\right)\right]},\label{bar:r:N:2}
\end{eqnarray}
where
\[
\zeta_{k,n}=\sum_{\nu=1}^\infty\left(\frac{n}{\tau_\nu^2+n(1+\lambda\rho_\nu)}\right)^k\;\mbox{for}
\;k=1,2.
\]
The final aggregated credible region is obtained as
\begin{equation}\label{strong:topo:cr}
R_N(\alpha):=\{f\in S^m(\mathbb{I}): \|f-\breve{f}_{N,\lambda}\|_2\le r_N(\alpha)\}.
\end{equation}

Our theorem below confirms that $R_N(\alpha)$ indeed possesses (asymptotic) posterior mass $(1-\alpha)$, and more importantly, proves that it covers the true function $f_0$ with probability tending to one.
\begin{Theorem}\label{cp:cr:strong}
	Suppose that 
	$f_0$ satisfies Condition (\textbf{S}),
	$m>1+\frac{\sqrt{3}}{2}$,
	$3/2<\beta<2m+1/(2m)-3/2$, $s=o(N^{\frac{\beta-1}{2m+\beta}})$, (\ref{s:for:aggregation})
	and $h\asymp h^\ast$.
	Then for any $\alpha\in (0,1)$, $P(R_N(\alpha)|\textbf{D})=1-\alpha+o_{P_{f_0}}(1)$ and
	$\lim_{n\rightarrow\infty}P_{f_0}(f_0\in R_N(\alpha))=1$.
\end{Theorem}

From the proof of Theorem~\ref{cp:cr:strong}, we point out that when $s=1$, the posterior mass of the aggregated credible region is exactly $1-\alpha$, consistent with \cite{SC14}. This remark also applies to other aggregated procedures to be presented later.

\begin{Remark}\label{rem}
	When $h\asymp h^\ast$, the radius of the aggregated ball $r_N(\alpha)\asymp N^{-\frac{2m+\beta-1}{2(2m+\beta)}}$ according to the discussions in Section \ref{sec:asymp:post:infer}. This is the optimal rate at which a posterior ball contracts based on the entire sample; see \cite{VZ08:2}. 
\end{Remark}

\subsection{Aggregated credible region in weak topology}\label{sec:cr:weak}

In this section, we invoke a weaker norm (than that used in Section~\ref{sec:cr:strong}) to construct an aggregated credible region. Under this new norm (inspired by \cite{CN13,CN14}), it is proven that the frequentist coverage \textit{exactly} matches with the asymptotic credibility level. The requirement on $s$ and $h$ in this section remains the same as Section \ref{sec:cr:strong}.

We define a weaker norm than $\|\cdot\|_2$, denoted $\|\cdot\|_\omega$. For any $f\in S^m(\mathbb{I})$ with $f=\sum_\nu f_\nu\varphi_\nu$, define $\|f\|_\omega^2=\sum_{\nu=1}^\infty \omega_\nu f_\nu^2$, where $\omega_\nu=(\nu (\log{2\nu}))^{-\tau}$ for some constant $\tau>1$. Since $\omega_\nu<1$ for all $\nu\ge1$, we
have $\|f\|_\omega\le \|f\|_2$. Under the new $\|\cdot\|_\omega$-norm, each individual $(1-\alpha)$ credible region is constructed as
$$R^\omega_{j,n}(\alpha)=\{f\in S^m(\mathbb{I}): \|f-\breve{f}_{j,n}\|_\omega\le r_{\omega,j,n}(\alpha)\},$$ where $r_{\omega,j,n}(\alpha)$ is directly obtained from posterior sampling such that $P(R^\omega_{j,n}(\alpha)|\textbf{D}_j)=1-\alpha$.

Under $\|\cdot\|_\omega$-norm, the aggregated credible region is constructed as:
\begin{equation}\label{weak:topo:cr}
R^\omega_N(\alpha):=\{f\in S^m(\mathbb{I}): \|f-\breve{f}_{N,\lambda}\|_\omega\le r_{\omega,N}(\alpha)\},
\end{equation}
where the radius is given as
\begin{equation}\label{bar:r:omega:N}
r_{\omega,N}(\alpha)=\sqrt{\frac{1}{s^2}\sum_{j=1}^s r_{\omega,j,n}^2(\alpha)}.
\end{equation}
Interestingly, Section \ref{sec:asymp:post:infer} illustrates that the aggregated radius $r_{\omega,N}(\alpha)$ contracts at root-$N$ rate.

Our theorem below shows that the frequentist covergage of $R^\omega_N(\alpha)$
exactly matches with the asymptotic posterior mass, both of which achieve the nominal level $(1-\alpha)$.
\begin{Theorem}\label{cp:cr:weak}
	Suppose that
	$f_0$ satisfies Condition (\textbf{S}),
	$m>1+\sqrt{3}/2$,
	$2\le\beta<\frac{(2m-1)^2}{2m}$, 
	$s=o(N^{\frac{\beta-1}{2m+\beta}})$, 
	$s=o(N^{\frac{4m^2+2m\beta-12m+1}{8m(2m+\beta)}}(\log{N})^{-\frac{3}{2}})$,
	and $h\asymp h^\ast$.
	Then for any $\alpha\in (0,1)$, $P(R^\omega_N(\alpha)|\textbf{D})=1-\alpha+o_{P_{f_0}}(1)$ and
	$\lim_{n\rightarrow\infty}P_{f_0}(f_0\in R^\omega_N(\alpha))=1-\alpha$.
\end{Theorem}


\subsection{Aggregated credible interval for linear functional}\label{sec:ci:functional}

In this section, we construct aggregated credible intervals for a class of linear functionals of $f$, denoted as $F(f)$. Examples include the evaluation functional, i.e., $F(f)=f(x)$, and integral functional, i.e., $F(f)=\int_0^1 f(x)dx$. Specifically, the interval is centered at $F(\breve f_{N,\lambda})$ with an length aggregated through $s$ lengths obtained from posterior sampling. Posterior and frequentist coverage properties of this aggregated interval depends on the functional form $F(\cdot)$. Again, our theory holds when $s$ is mildly diverging and $h\asymp h^\ast$.

Let $F:S^m(\mathbb{I})\mapsto\bbR$ be a linear $\Pi$-measurable functional satisfying the following
Condition (\textbf{F}):
$\sup_{\nu\ge1}|F(\varphi_\nu)|<\infty$, and
there exist constants $\kappa>0$ and $r\in[0,1]$ such that for any $f\in S^m(\mathbb{I})$,
\begin{equation}\label{boundedness:F}
|F(f)|\le\kappa h^{-r/2}\|f\|.
\end{equation}
It follows by \cite{SC14} that the evaluation functional satisfies
Condition (\textbf{F}) with $r=1$ and the integral functional satisfies
Condition (\textbf{F}) with $r=0$.

Based on each $\textbf{D}_j$, we obtain from posterior samples the following $(1-\alpha)$ credible interval:
\[
CI_{j,n}^F(\alpha):=\{f\in S^m(\mathbb{I}): |F(f)-F(\breve{f}_{j,n})|\le r_{F,j,n}(\alpha)\},
\]
where $r_{F,j,n}(\alpha)$ is a radius such that $P(CI_{j,n}^F(\alpha)|\textbf{D}_j)=1-\alpha$. The aggregated credible interval is constructed as
\begin{equation}\label{aggregated:CI}
CI_N^F(\alpha):=\{f\in S^m(\mathbb{I}): |F(f)-F(\breve{f}_{N,\lambda})|\le \bar{r}_{F,N}(\alpha)\}
\end{equation}
where
\begin{equation}\label{bar:r:F:N}
r_{F,N}(\alpha)=\frac{\theta_{1,N}}{\theta_{1,n}}\sqrt{\frac{1}{s}\sum_{j=1}^s r_{F,j,n}(\alpha)^2}\;\;\;\mbox{and}
\;\;\;\theta_{k,n}^2=\sum_{\nu=1}^\infty
\frac{F(\varphi_\nu)^2}{(\tau_\nu^2+n(1+\lambda\rho_\nu))^k}\;\mbox{for}\;k=1,2.
\end{equation}
The shrinking rate of $\bar{r}_{F,N}(\alpha)$ depends on the functional form $F$; see  Section \ref{sec:asymp:post:infer}.

Our theorem below investigates the asymptotic properties of $CI_F^N(\alpha)$ in terms of both posterior and frequentist coverage.
\begin{Theorem}\label{cp:ci:functional}
	Suppose that
	$f_0=\sum_{\nu=1}^\infty f_\nu^0\varphi_\nu$
	satisfies Condition (\textbf{S}$'$): $\sum_{\nu=1}^\infty |f_\nu^0|^2\nu^{2m+\beta}<\infty$,
	$E_{f_0}\{\epsilon^4|X\}\le M_4$ a.s. for some constant $M_4>0$,
	$N^k\theta_{k,N}^2\gtrsim h^{-r}$ for $k=1,2$,
	$m>1+\frac{\sqrt{3}}{2}$,
	$2\le\beta<\frac{(2m-1)^2}{2m}$, $s=o(N^{\frac{\beta-1}{2m+\beta}})$,
	$s=o(N^{\frac{4m^2+2m\beta-12m+1}{8m(2m+\beta)}}(\log{N})^{-\frac{3}{2}})$, (\ref{s:for:aggregation})
	and $h\asymp h^\ast$.
	Then for any $\alpha\in (0,1)$,
	$P(CI_N^F(\alpha)|\textbf{D})=1-\alpha+o_{P_{f_0}}(1)$,
	and $\lim\inf_{N\rightarrow\infty}P_{f_0}(f_0\in CI_N^F(\alpha))\ge 1-\alpha$ given that Condition (\textbf{F}) holds.
	Moreover, if $0<\sum_{\nu=1}^\infty F(\varphi_\nu)^2<\infty$,
	then $\lim_{N\rightarrow\infty}P_{f_0}(f_0\in CI_N^F(\alpha))=1-\alpha$.
\end{Theorem}
\noindent Note that Condition (\textbf{S}$'$) is slightly stronger than Condition (\textbf{S}) required in Theorem~\ref{uniform:bvm:thm}. 
Indeed, this condition essentially means that $f_0$ has derivatives up to order $m+\frac{\beta}{2}$
(when this order is integer-valued).
Hence, Theorem \ref{cp:ci:functional} requires a more smooth true function $f_0$.

It was shown in \cite{SC14} that the integral functional $F_x(f):=\int_0^x f(z)dz$ for any $x\in[0,1]$ satisfies (\ref{boundedness:F}) with $r=0$ and $0<\sum_{\nu=1}^\infty F_x(\varphi_\nu)^2<\infty$. Therefore, the $(1-\alpha)$-th credible interval
of $F_x(f)$ achieves exactly $(1-\alpha)$ frequentist coverage, while that for the evaluation functional is more conservative. These theoretical findings will be empirically verified in Section \ref{sec:simulations} .

\subsection{Asymptotic aggregated inference}\label{sec:asymp:post:infer}

In practice, the centers $\breve{f}_{N,\lambda}$, $F(\breve{f}_{N,\lambda})$
and the radii $r_{j,n}(\alpha)$, $r_{\omega,j,n}(\alpha)$, $r_{F,j,n}(\alpha)$
in Sections \ref{sec:cr:strong} -- \ref{sec:ci:functional} are directly obtained from posterior samples.
Sometimes posterior sampling is time consuming and inefficient, particularly as $s\to\infty$. This computational consideration motivates us to propose an \textit{asymptotic} approach in which one replaces the above centers/radii by their large sample limits. Our new asymptotic inference procedures dramatically improve the computing speed,
as displayed in simulations; see Section \ref{sec:simulations}.

Define
\begin{equation}\label{bar:r:N:1}
\widetilde{f}_{N,\lambda}(\cdot)=\sum_{\nu=1}^\infty \frac{a_{N,\nu}}{a_{n,\nu}}
V\left(\frac{1}{s}\sum_{j=1}^s\widetilde{f}_{j,n},\varphi_\nu\right)\varphi_\nu(\cdot).
\end{equation}
Clearly, $\widetilde{f}_{N,\lambda}$ is a counterpart of $\breve{f}_{N,\lambda}$ (\ref{aggreg:center})
with $\breve{f}_{j,n}$ therein replaced by $\widetilde{f}_{j,n}$.
By a careful examination of the proofs of Theorems \ref{cp:cr:strong} -- \ref{cp:ci:functional}, it can be shown that the following limits hold:
\begin{eqnarray}\label{some:limits}
\|\breve{f}_{N,\lambda}-\widetilde{f}_{N,\lambda}\|&=&o_{P_{f_0}}(N^{-1/2}h^{-1/4}),\nonumber\\
\max_{1\le j\le s}
\bigg|\frac{nr_{j,n}^2(\alpha)-\zeta_{1,n}}{\sqrt{2\zeta_{2,n}}}-z_\alpha\bigg|&=&o_{P_{f_0}}(1),\nonumber\\
\max_{1\le j\le s}|\sqrt{n}r_{\omega,j,n}(\alpha)-\sqrt{c_\alpha}|&=&o_{P_{f_0}}(1),\nonumber\\
\max_{1\le j\le s}|r_{F,j,n}(\alpha)/\theta_{1,n}-z_{\alpha/2}|&=&o_{P_{f_0}}(1),
\end{eqnarray}
where $z_\alpha=\Phi^{-1}(1-\alpha)$ with $\Phi(\cdot)$ being the c.d.f. of standard normal random variable,
and $c_\alpha>0$ satisfies $P(\sum_{\nu=1}^\infty d_\nu\eta_\nu^2\le c_\alpha)=1-\alpha$
with $\eta_\nu$ being independent standard normal random variables.

It yields from (\ref{some:limits}) that
the following approximation relationships hold
uniformly for $1\le j\le s$:
\begin{eqnarray*}
	r_{j,n}(\alpha)\approx\sqrt{\frac{\zeta_{1,n}+\sqrt{2\zeta_{2,n}}z_\alpha}{n}},\;\;\;
	r_{\omega,j,n}(\alpha)\approx\sqrt{\frac{c_\alpha}{n}}\;\;\;\mbox{and}\;\;\;
	r_{F,j,n}(\alpha)\approx\theta_{1,n}z_{\alpha/2},
\end{eqnarray*}
which further implies (by the aggregation formulae (\ref{bar:r:N}), (\ref{bar:r:omega:N}) and (\ref{bar:r:F:N}))
\begin{equation}\label{apprad}
\begin{aligned}
&&r_N(\alpha)\approx r_N^\dag(\alpha):=\sqrt{\frac{\zeta_{1,N}+\sqrt{2\zeta_{2,N}}z_\alpha}{N}},\\
&&r_{\omega,N}(\alpha)\approx r_{\omega,N}^\dag(\alpha):=\sqrt{\frac{c_\alpha}{N}},\\
&&r_{F,N}(\alpha)\approx r_{F,N}^\dag(\alpha):=\theta_{1,N}z_{\alpha/2}.
\end{aligned}
\end{equation}
Thus, we have the following {\em asymptotic} counterparts of $R_N(\alpha)$, $R^\omega_N(\alpha)$  
and $CI_N^F(\alpha)$:
\begin{eqnarray}
&&R_N^\dag(\alpha):=\{f\in S^m(\mathbb{I}):
\|f-\widetilde{f}_{N,\lambda}\|_2\le r_N^\dag(\alpha)\},\label{asymp:strong:topo:cr}\\
&&R^{\dag\omega}_N(\alpha):=\{f\in S^m(\mathbb{I}):
\|f-\widetilde{f}_{N,\lambda}\|_\omega\le r_{\omega,N}^\dag(\alpha)\},\label{asymp:weak:topo:cr}\\
&& CI_N^{\dag F}(\alpha):=\{f\in S^m(\mathbb{I}): |F(f)-F(\widetilde{f}_{N,\lambda})|\le r_{F,N}^\dag(\alpha)\}.\label{asymp:aggregated:CI}
\end{eqnarray}

Our theorem below shows that the posterior coverage and frequentist coverage of the above computationally efficient alternatives remain the same as those for $R_N(\alpha)$, $R^\omega_N(\alpha)$ and $CI_N^F(\alpha)$ under the same set of conditions.
\begin{Theorem}\label{cp:cr:strong:asymp}
	Suppose that all assumptions in Theorems \ref{cp:cr:strong} -- \ref{cp:ci:functional} hold. Then for any $\alpha\in (0,1)$, $R_N^\dag(\alpha)$, $R^{\dag\omega}_N(\alpha)$ and $CI_N^{\dag F}(\alpha)$ possess exactly the same posterior and frequentist properties as  $R_N(\alpha)$, $R^\omega_N(\alpha)$ and $CI_N^F(\alpha)$, respectively.
\end{Theorem}

As a byproduct, (\ref{apprad}) implies the contraction rate of each aggregated credible ball/interval in Sections~\ref{sec:cr:strong} -- \ref{sec:asymp:post:infer}. It is easy to see that $r_{\omega,N}(\alpha)\asymp N^{-1/2}$. As for $r_{F,N}(\alpha)$, it depends on the functional form $F$. For example, when $F$ is an evaluation functional, it holds that $\theta_{1,N}^2\asymp (Nh)^{-1}$,
leading to $N^{-\frac{2m+\beta-1}{2(2m+\beta)}}$ when $h\asymp h^*$;
when $F$ is an integral functional, we have $r_{F,N}(\alpha)\asymp N^{-1/2}$ since $\theta_{1,N}^2\asymp N^{-1}$. As for $r_N(\alpha)$, it can be shown by a simple fact $\zeta_{1,N}, \zeta_{2,N}\asymp h^{-1}$ that $r_N(\alpha)\asymp (Nh)^{-1/2}\asymp N^{-\frac{2m+\beta-1}{2(2m+\beta)}}$ when $h\asymp h^*$. This contraction rate turns out to be optimal based on the entire sample; see \cite{VZ08:2}. However, if we choose $h$ in the scale of subsample size $n$, e.g., $h\asymp n^{-\frac{1}{2m+\beta}}$, similar arguments show that $r_N(\alpha)\asymp N^{-\frac{2m+\beta-1}{2(2m+\beta)}}s^{-\frac{1}{2(2m+\beta)}}$. Hence, such a region contracts faster than the optimal rate, which results in unsatisfactory frequentist coverage.

Table \ref{table:comparison:aggregation} summarizes six aggregated credible regions/intervals from Sections \ref{sec:cr:strong} -- \ref{sec:ci:functional} in terms of their centers and radii.

\begin{table}[htp]
	\begin{center}
		\caption{Summary of Aggregated $(1-\alpha)$ Credible Regions/Intervals}
		\begin{tabular}{ccccc}\\ \hline
			Type&Name& Notation &Center&Radius\\ \hline\\
			
			\multirow{3}{30mm}{Finite-sample}&strong CR for $f$ & $R_N(\alpha)$ & $\breve{f}_{N,\lambda}$ & $r_N(\alpha)$ \\
			&weak CR for $f$   & $R^\omega_N(\alpha)$& $\breve{f}_{N,\lambda}$ &$r_{\omega,N}(\alpha)$\\
			&CI for $F(f)$     & $CI_N^F(\alpha)$ & $F(\breve{f}_{N,\lambda})$ & $r_{F,N}(\alpha)$\\ \\
			
			\multirow{3}{30mm}{Asymptotic}& strong CR for $f$ & $R_N^\dag(\alpha)$ &
			$\widetilde{f}_{N,\lambda}$ & $r_N^\dag(\alpha)$ \\
			&weak CR for $f$   & $R^{\dag\omega}_N(\alpha)$& $\widetilde{f}_{N,\lambda}$ &$r_{\omega,N}^\dag(\alpha)$\\
			&CI for $F(f)$     & $CI_N^{\dag F}(\alpha)$ & $F(\widetilde{f}_{N,\lambda})$ & $r_{F,N}^\dag(\alpha)$\\
			
			\hline
		\end{tabular}
		\label{table:comparison:aggregation}
	\end{center}
\end{table}%

\section{Simulation Study}\label{sec:simulations}

In this section, statistical properties of the proposed aggregated procedures are examined using a simulation study.
We generated samples from the following model
\begin{equation}\label{normal:regression:model}
Y_{ij}=f_0(X_{ij})+\epsilon_{ij},\,\,i=1,2,\ldots,n, j=1,2,\ldots,s,
\end{equation}
where $X_{ij}\overset{iid}{\sim}Unif[0,1]$,
$\epsilon_{ij}\overset{iid}{\sim}N(0,1)$,
and $\epsilon_{ij}$ are independent of $X_{ij}$. 
The true regression function was
chosen to be $f_0(x)=2.4\beta_{30,17}(x)+1.6\beta_{3,11}(x)$,
where $\beta_{a,b}$ is the probability density function for $Beta(a,b)$.

Consider GP prior $f\sim \sum_{\nu=1}^n w_\nu\varphi_\nu$, where $w_\nu$ are defined in (\ref{prior1}).
The proposed Bayesian procedures were examined.
Specifically, we computed the frequentist
coverage proportions (CP) of the credible regions (\ref{strong:topo:cr}),
(\ref{weak:topo:cr}), (\ref{asymp:strong:topo:cr}), (\ref{asymp:weak:topo:cr}),
and credible intervals (\ref{aggregated:CI}), (\ref{asymp:aggregated:CI}).
In particular, (\ref{strong:topo:cr}),
(\ref{weak:topo:cr}) and (\ref{aggregated:CI}) were constructed based
on posterior samples, as described in
Sections \ref{sec:aggr:post:mean}--\ref{sec:ci:functional}; whereas
(\ref{asymp:strong:topo:cr}), (\ref{asymp:weak:topo:cr})
and (\ref{asymp:aggregated:CI}) were constructed based on 
asymptotic theory developed in Section \ref{sec:asymp:post:infer}.
To ease presentation,
we call (\ref{strong:topo:cr}) and
(\ref{weak:topo:cr}) as finite-sample credible regions (FCR), and call (\ref{asymp:strong:topo:cr}) and (\ref{asymp:weak:topo:cr}) as asymptotic credible regions (ACR).

The calculation of CP was based on $500$ independent experiments. Specifically, the CP is the proportion of the credible regions/intervals containing $f_0$/$F(f_0)$ (for a linear functional $F$).
Two types of $F$ were considered: (1) the evaluation functional $F_x(f)=f(x)$ for any $x\in[0,1]$,
and (2) the integral functional $F_x(f)=\int_0^x f(z)dz$ for any $x\in[0,1]$. In both cases, we consider $F_x$ with $x$ being 15 evenly spaced points in [0.05,0.95]. To make the study more complete,
a set of credibility levels were examined, i.e., $1-\alpha=0.1,0.3,0.5,0.7,0.9,0.95$.
In each experiment, $N=1200$ independent samples were generated from the model (\ref{normal:regression:model}). For ACR and FCR, we chose the number of divisions
$s=1,2,3,4,5,6,8,10,12,15,20,24,30,40,60$.
Define $\gamma=\log{s}/\log{N}$.
Note that $s=1$ (equivalently, $\gamma=0$) means ``no division." 

\begin{figure}

\centering

\includegraphics[width=1\textwidth,height=1\textwidth]{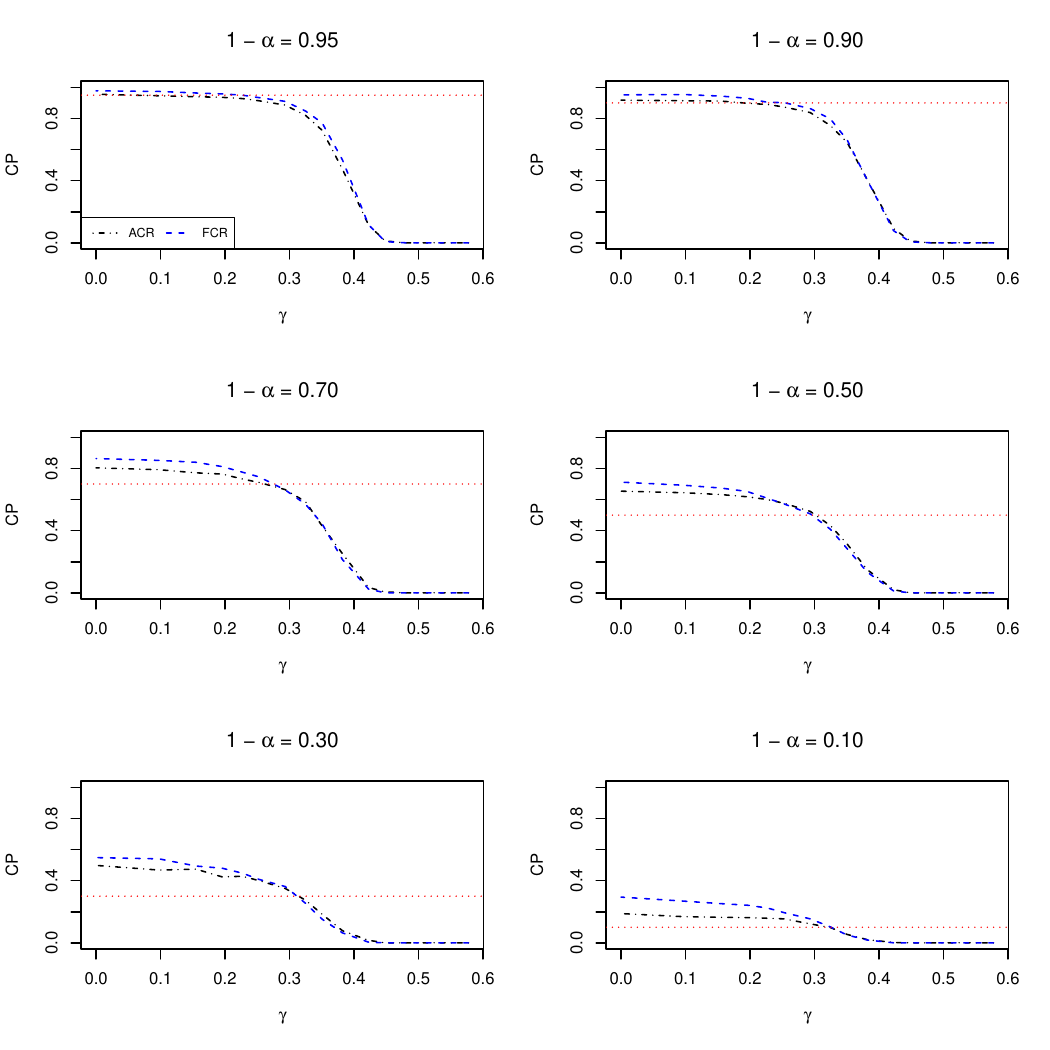}
\caption{\textit{\footnotesize CP of ACR and FCR based on strong topology.
		Dotted red lines indicate credibility levels.}}
\label{fig:strong:cr}
\end{figure}

\begin{figure}
\centering
\includegraphics[width=1\textwidth,height=1\textwidth]{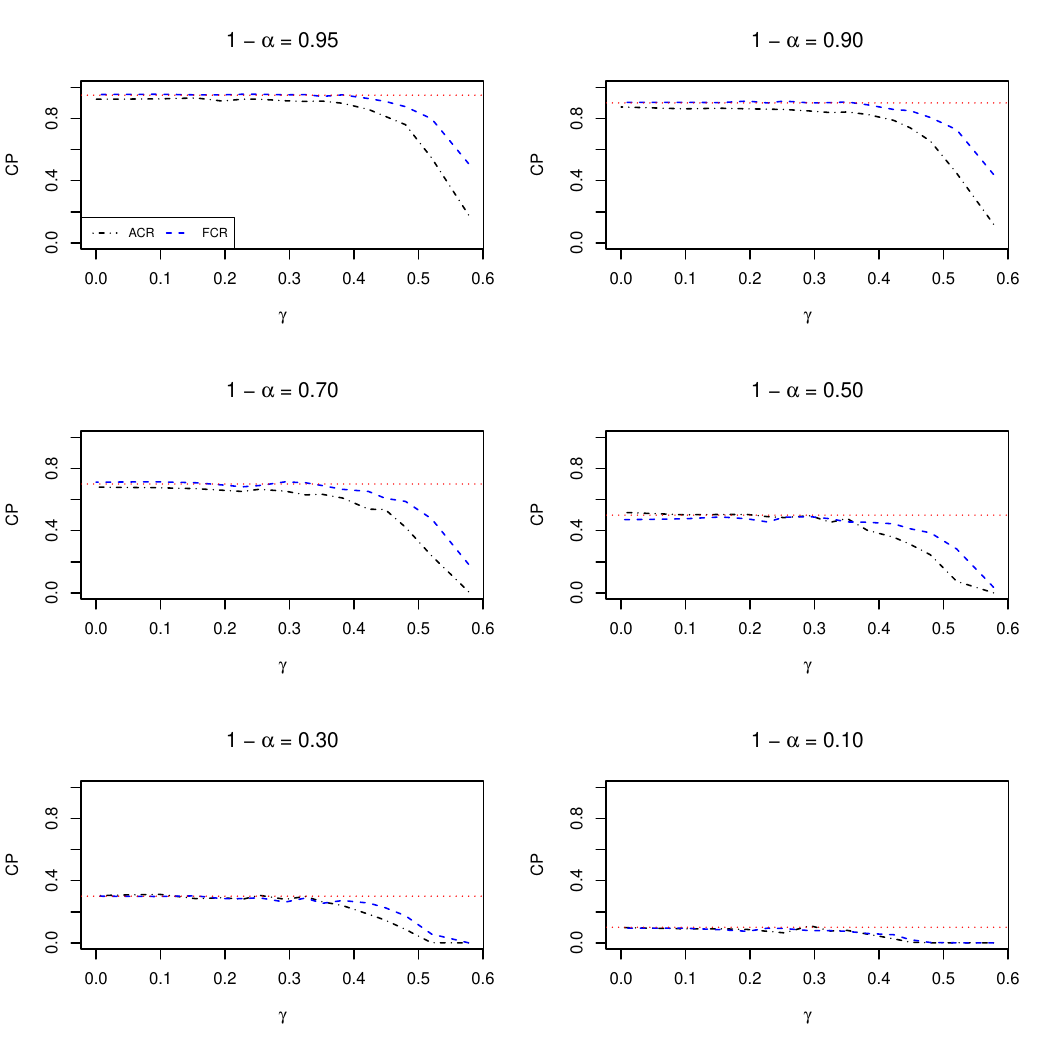}
\caption{\textit{\footnotesize CP of ACR and FCR based on weak topology. Dotted red lines indicate credibility levels.}}
\label{fig:weak:cr}
\end{figure}

Figure \ref{fig:strong:cr}
demonstrates the results
for FCR and ACR based on strong topology,
i.e., (\ref{strong:topo:cr}) and (\ref{asymp:strong:topo:cr}).
The red dotted line indicates the $(1-\alpha)$ credibility level.
It can be seen that the CP of both FCR and ACR is above the credibility levels when $\gamma$ is small, while it suddenly drops to zero as $\gamma$ is beyond some threshold, say $0.3$. This observation supports our theory that $s$ should not grow too fast, and that the credible regions based on strong topology tends to be more ``conservative."
Figure \ref{fig:weak:cr} demonstrates the results
for FCR and ACR based on weak topology,
i.e., (\ref{weak:topo:cr}) and (\ref{asymp:weak:topo:cr}).
We observe that the CP of both ACR and FCR approaches the desired credibility levels when $\gamma\le 0.3$, but quickly drops to zero when $\gamma$ becomes large. This observation also supports our theory that the use of weak topology leads to a more satisfactory frequentist coverage.

\begin{figure}
\centering
\includegraphics[width=1\textwidth,height=1\textwidth]{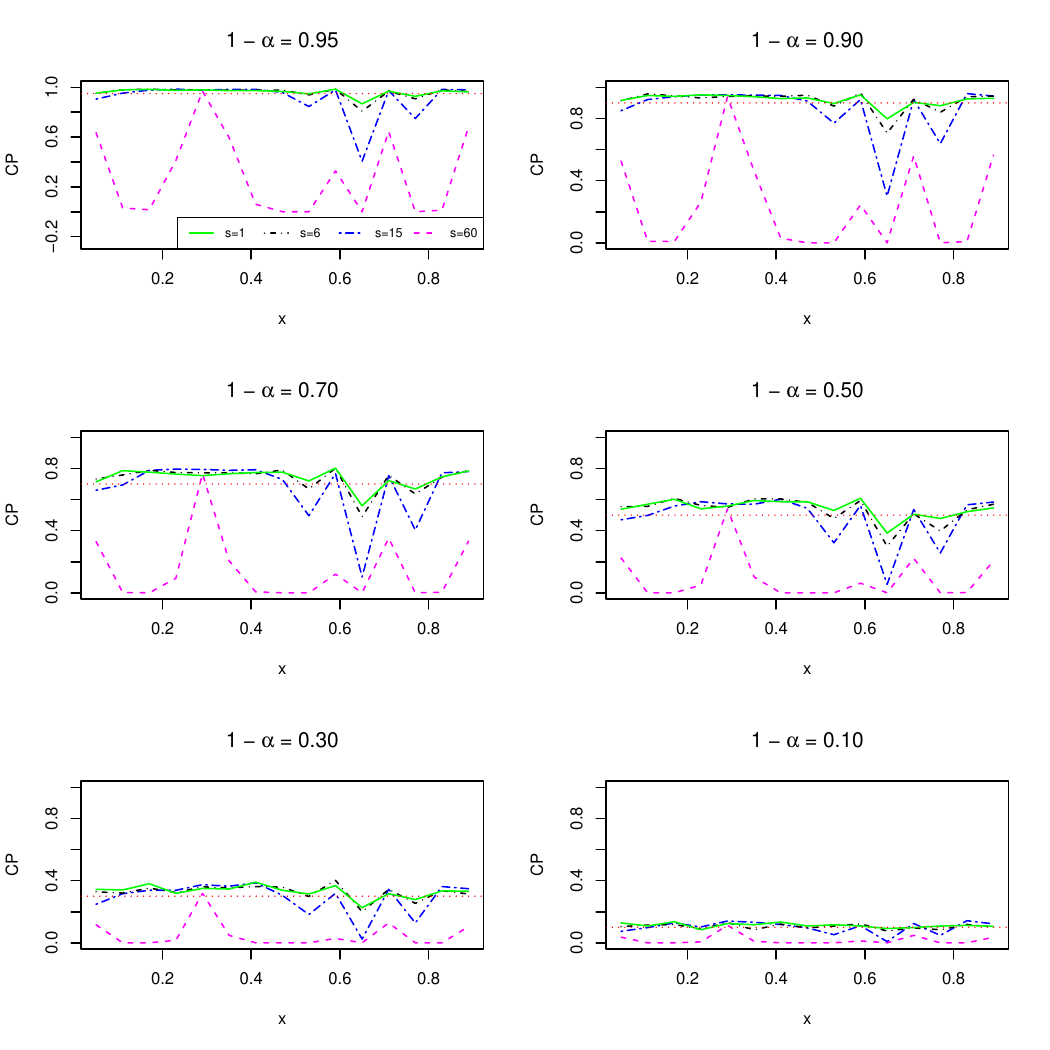}
\caption{\textit{\footnotesize CP of $F_x(f)=f(x)$ against $x$ based on posterior samples of $f$. Dotted red lines indicate credibility levels.}}
		\label{fig:finite:pointwise:CI}
\end{figure}

\begin{figure}
\centering
\includegraphics[width=1\textwidth,height=\textwidth]{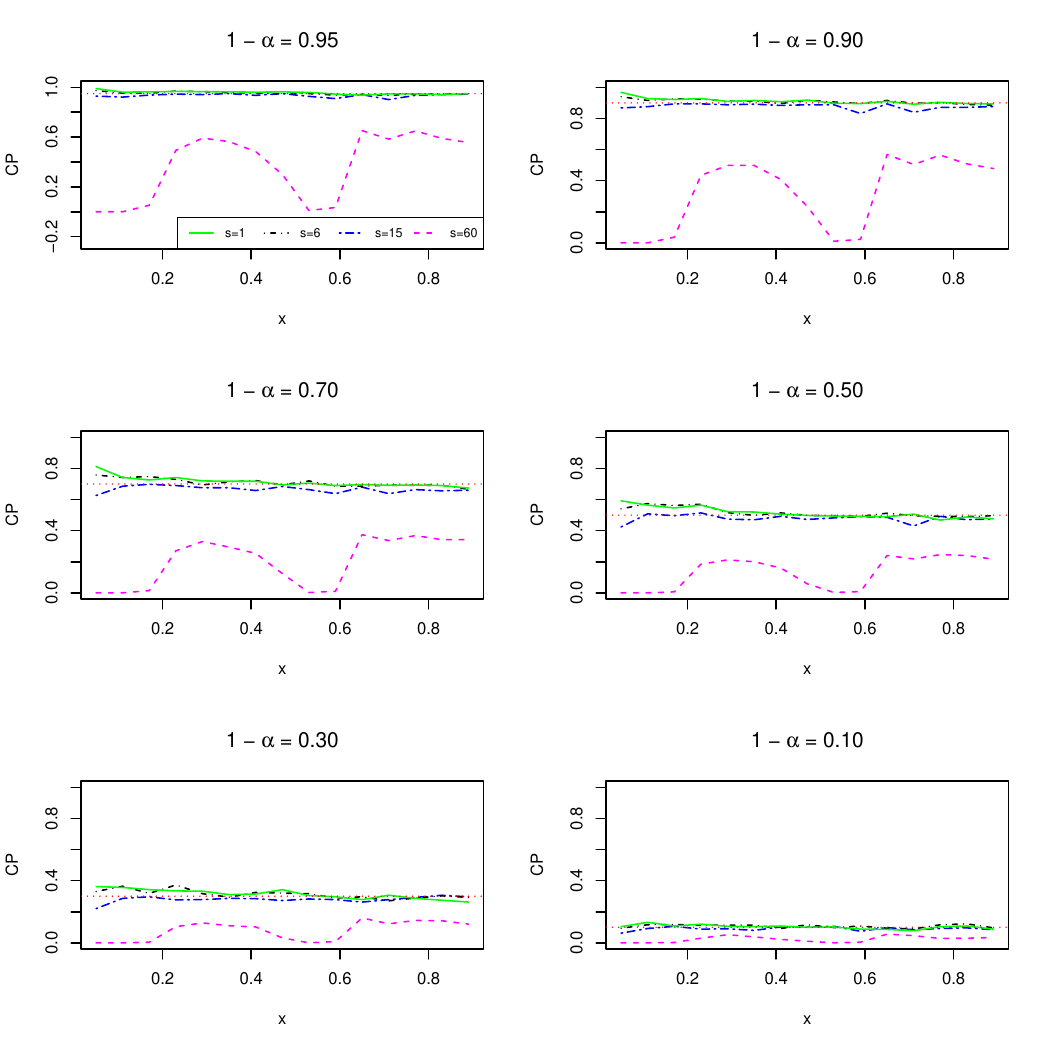}
\caption{\textit{\footnotesize CP of $F_x(f)=\int_0^xf(z)dz$ against $x$ based on posterior samples of $f$. Dotted red lines indicate credibility levels.}}
\label{fig:finite:integral:CI}

\end{figure}

For credible intervals of linear functionals, we chose the number of divisions
$s=1,6,15,60$.
Figures \ref{fig:finite:pointwise:CI} and \ref{fig:finite:integral:CI}
display the results for evaluation functional and integral functional, respectively,
based on posterior samples.
It can be seen that when $s=60$, the CP of the credible intervals for the evaluation functional
drops to zero at most of the $x$ points, indicating the failure in covering the true
values of the function. However, when $s=1,6,15$, the CP is above
the credibility levels except for the points where
the true function $f_0$ has peaks; see (a) of Figure \ref{fig:toy:example}.
The observation that the CP stays above $(1-\alpha)$ coincides with our theory
that the credible interval of the evaluation functional is
conservative. On the other hand, 
it can be seen that when $s=60$, the CP of the credible intervals for the integral functional
becomes far below the credibility levels
at most $x$. However, when $s=1,6,15$, the CP is close to the credibility levels at all $x$.
This finding coincides with our theory that
the the credible interval of the integral functional
achieves exactly $(1-\alpha)$ frequentist coverage. The above results 
also support our claim that $s$ cannot grow too fast for guaranteeing frequency validity.
Credible intervals based on asymptotic theory, i.e., (\ref{asymp:aggregated:CI}),
were summarized in Figures \ref{fig:asymp:pointwise:CI} and \ref{fig:asymp:integral:CI} of the supplement document
\cite{SCBigRate}. Interpretations of these results are similar to those based on finite posterior samples.


\begin{figure}[htp]
	\begin{center}
		\includegraphics[scale=0.4]{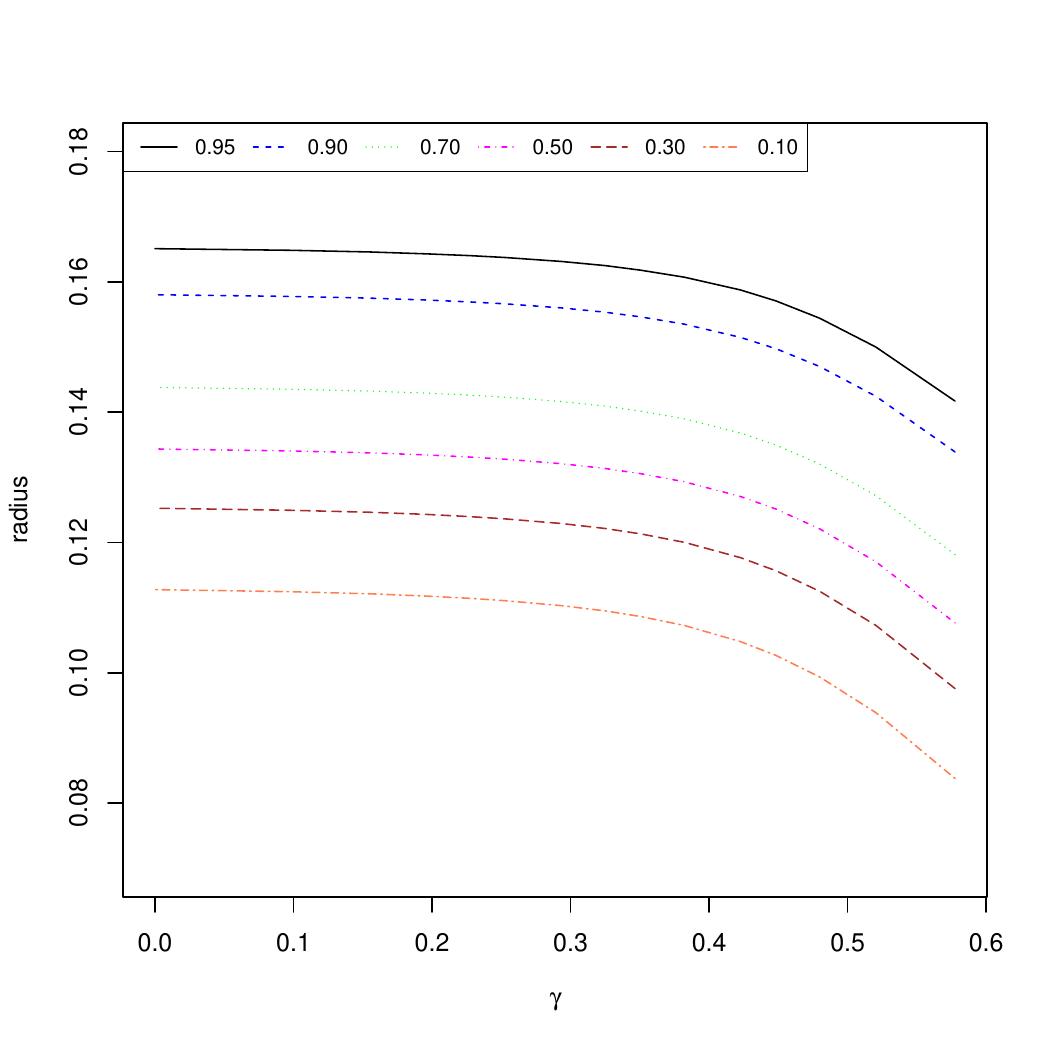}
		\caption{\textit{\footnotesize Radius of credible region (\ref{strong:topo:cr}) against
				$\gamma$. Legend indicates the credibility levels $1-\alpha$.}}
		\label{fig:strong_compare_radius}
	\end{center}
\end{figure}

\begin{figure}[htp]
	\begin{center}
		\includegraphics[scale=0.4]{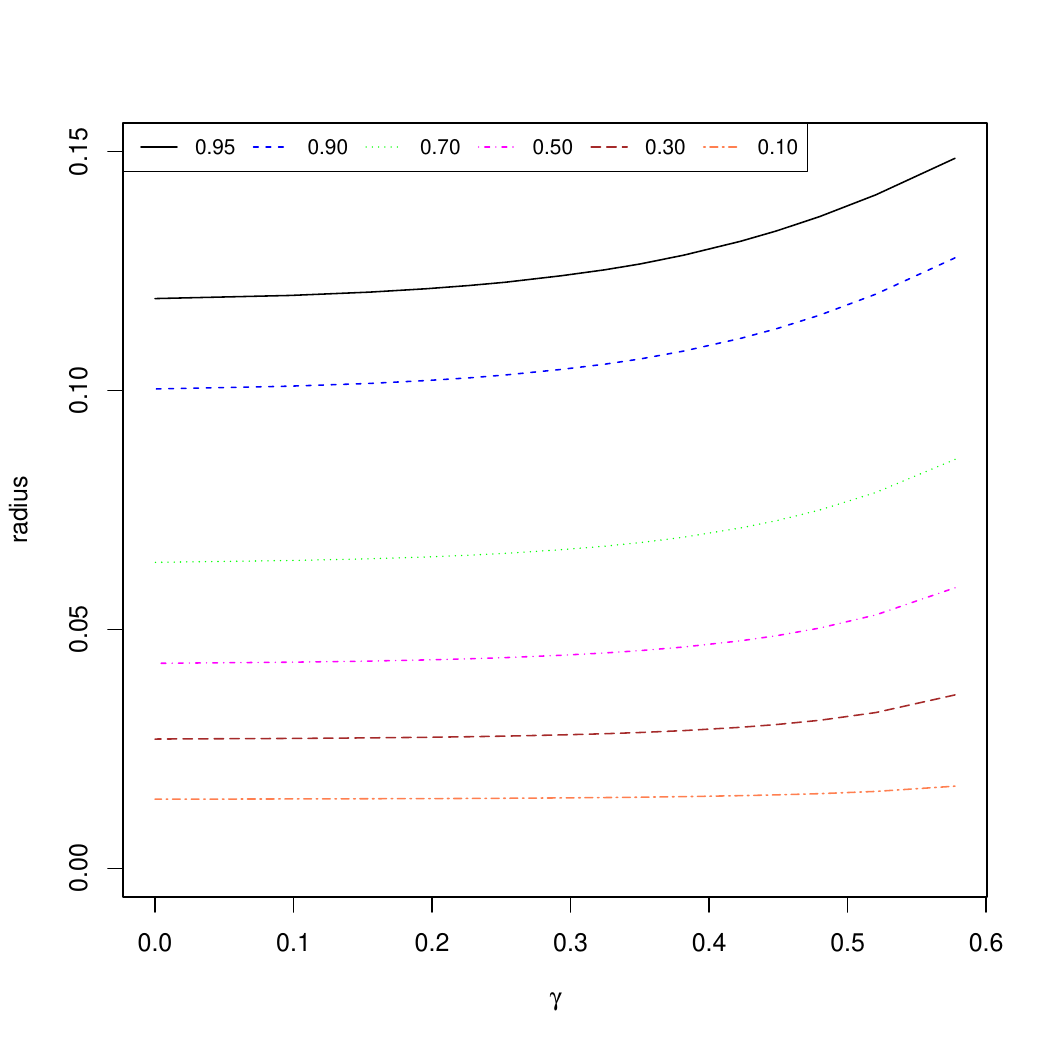}
		\caption{\textit{\footnotesize Radius of credible region (\ref{weak:topo:cr}) against
				$\gamma$. Legend indicates the credibility levels $1-\alpha$.}}
		\label{fig:weak_compare_radius}
	\end{center}
\end{figure}

\begin{figure}[htp]
	\begin{center}
		\includegraphics[scale=0.65]{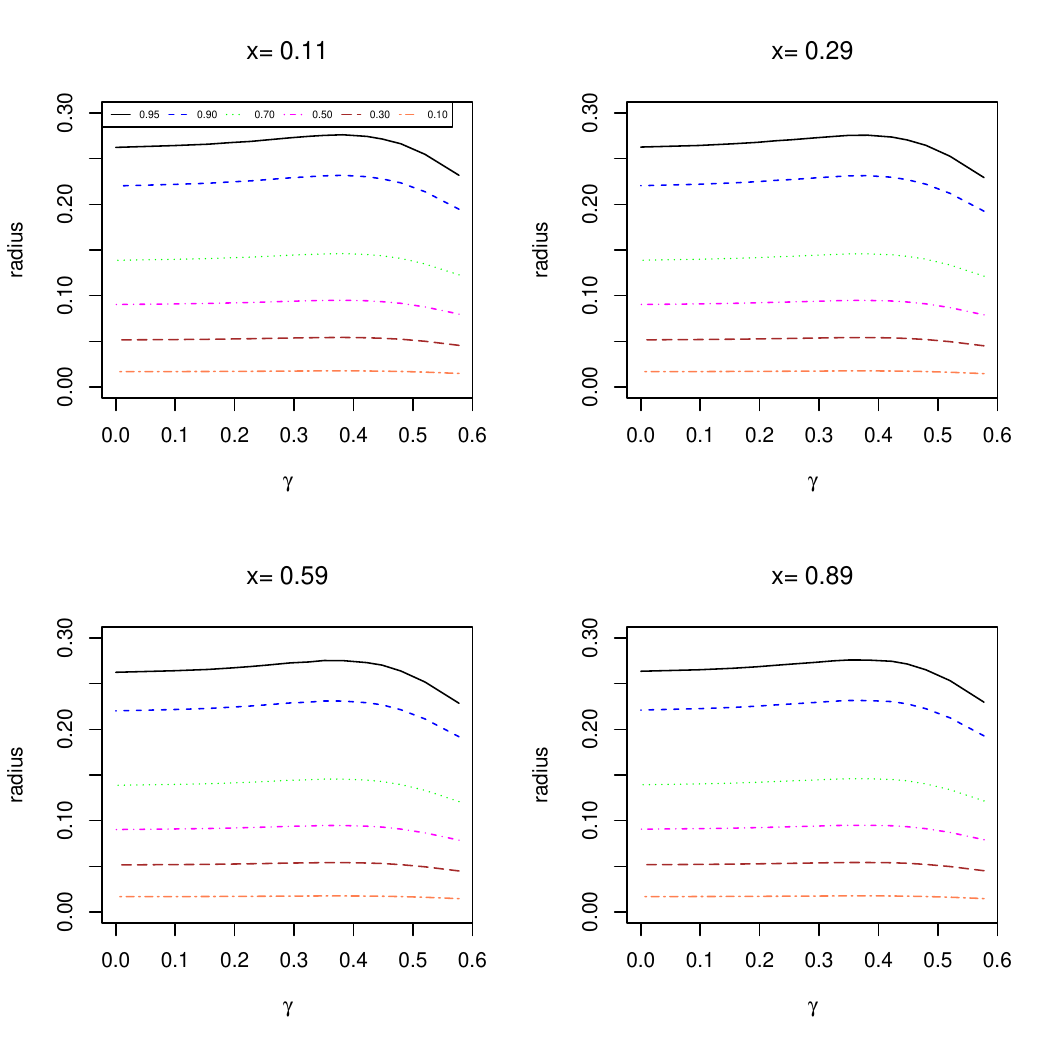}
		\caption{\textit{\footnotesize Radius of credible interval (\ref{aggregated:CI})
				for pointwise functional $F_x(f)=f(x)$ against
				$\gamma$. Legend indicates the credibility levels $1-\alpha$.
				Four values of $x$ are considered.}}
		\label{fig:pointwise_compare_radius}
	\end{center}
\end{figure}

\begin{figure}[htp]
	\begin{center}
		\includegraphics[scale=0.65]{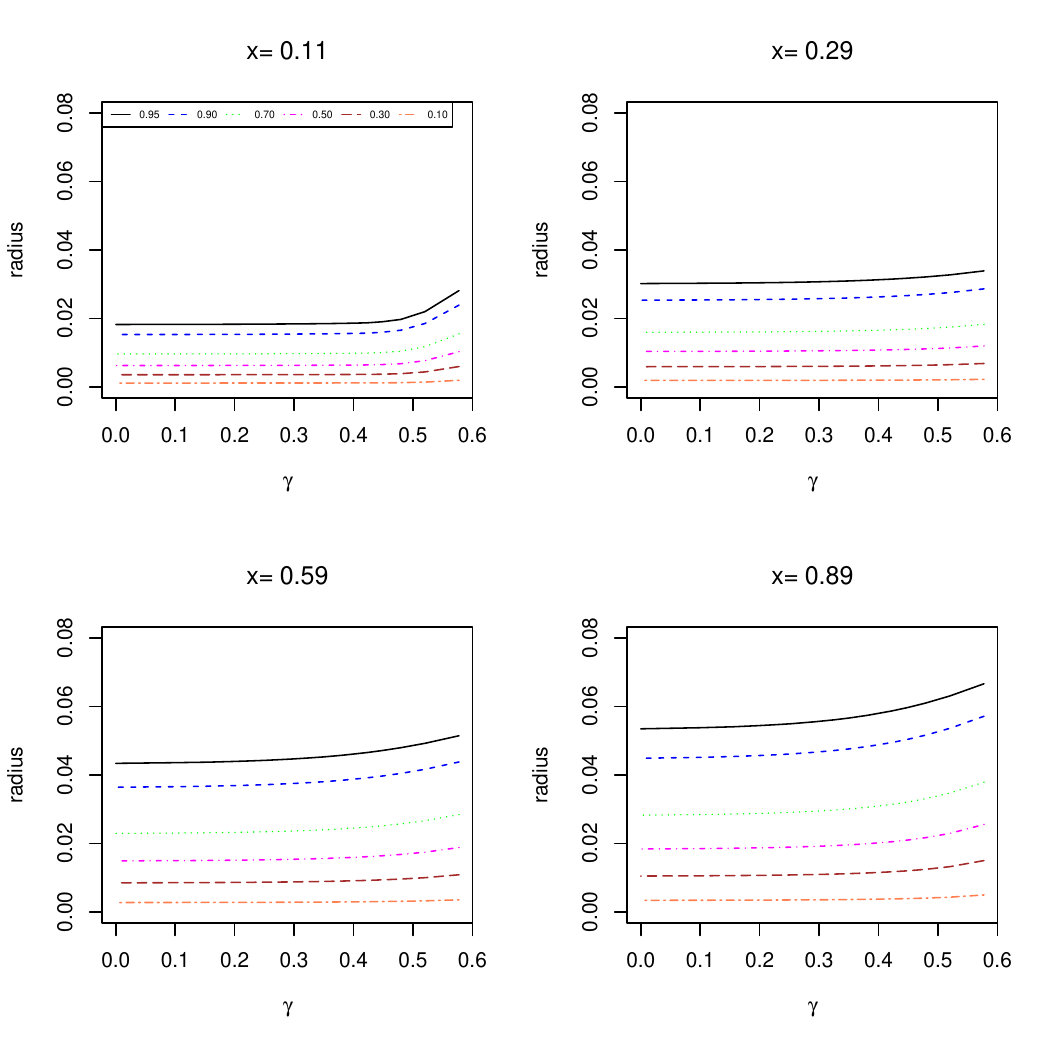}
		\caption{\textit{\footnotesize Radius of credible interval (\ref{aggregated:CI})
				for integral functional $F_x(f)=\int_0^x f(z)dz$ against
				$\gamma$. Legend indicates the credibility levels $1-\alpha$.
				Four values of $x$ are considered.}}
		\label{fig:integral_compare_radius}
	\end{center}
\end{figure}

The supplement document \cite{SCBigRate} also includes
Figures \ref{fig:strong_compare_radius} -- \ref{fig:integral_compare_radius}
which demonstrate how the radii/lengths of the aggregated credible
regions/intervals change along with $\gamma$, the size of the subsample. It can be observed that when $\gamma\le 0.3$, indicating that the full sample is divided into at most twelve subsamples, the radii of the aggregated
regions/intervals are almost identical to the radii of
the regions/intervals directly constructed from the full sample, i.e., $\gamma=0$.
This means that our aggregated procedures, based on a suitable amount of divisions, indeed mimic the oracle procedures.
However, when $\gamma$ increases to $0.6$, the distinctions between the
the aggregated and oracle procedures quickly become obvious.

We also repeated the above study for
$N=1800$ and $2400$. The plots corresponding to these studies are given in supplement document;
see Section \ref{sec:suppl:additional:simulation} of \cite{SCBigRate}.
The interpretations of these additional results are similar as above.

\begin{figure}[h]
\centering
\includegraphics[scale=0.65]{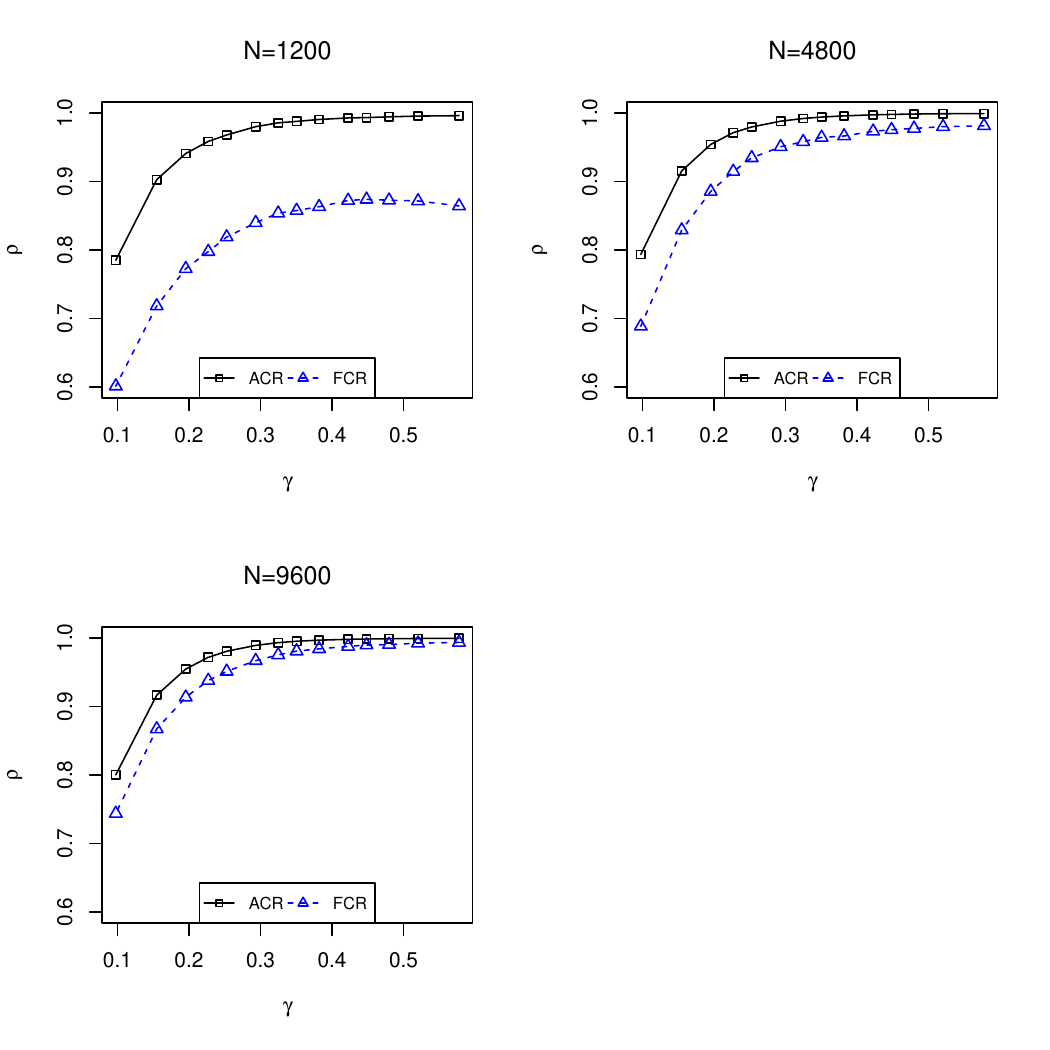}
\caption{\textit{\footnotesize $\rho$ versus $\gamma$ based on FCR and ACR for single experiment.}}
\label{fig:running:time}

\end{figure}

To the end of this section, computing efficiency is investigated. Figure \ref{fig:running:time} displays the results based on a single experiment for various choices of $N$. Specifically, we look at the value of the quantity $\rho=1-(T/T_0)$ versus a collection of $\gamma$'s for FCR and ACR, where $T_0$ ($T$) is the computing time without using D\&C (based on D\&C). We observe that
$T$ is substantially smaller than $T_0$, and this computation efficiency (as reflected by the value of $\rho$) becomes more obvious as $\gamma$ grows for each fixed $N$. This can also be seen as $N$ grows for each fixed $\gamma$. However, this reduction in computing time does not affect the performances of the aggregated credible regions  when $0\le \gamma\le 0.3$, as demonstrated in Figures \ref{fig:strong:cr}, \ref{fig:weak:cr}, \ref{fig:strong_compare_radius}--\ref{fig:integral_compare_radius}.

\section{Real Data Analysis}\label{sec:real}
In this section, we apply our methods to Million Song Data (MSD) and Flight Delay Data (FDD).

\subsection{Million Song Data}
As a real application, we apply our aggregation procedure to analyze MSD. 
The MSD is a perfect example of large dataset, a freely-available
collection of audio features and metadata for a million contemporary
popular music tracks. Each observation is a song
track released between the year 1922 and 2011. The response variable $Y_i$ is the year when the
song was released and the covariate $X_i$ is the timbre average of the song. The main purpose is to explore a
relationship, denoted as $f$, between song features and years in a nonparametric regression model, i.e.,
\textrm{$\textrm{year}=f(\textrm{timbre})$+error}.
The above model is useful to predict production year based on song timbre.
Due to enormous sample size,
processing the entire data is infeasible. 
In frequentist setting,
a distributed kernel ridge regression method 
was proposed by \cite{ZDW15,ZWJ15} for estimation
purposes (without quantifying uncertainty). 

\begin{figure}[h]
	\centering
		\begin{minipage}[t]{0.45\linewidth}
		\centering
		\includegraphics[scale=0.3]{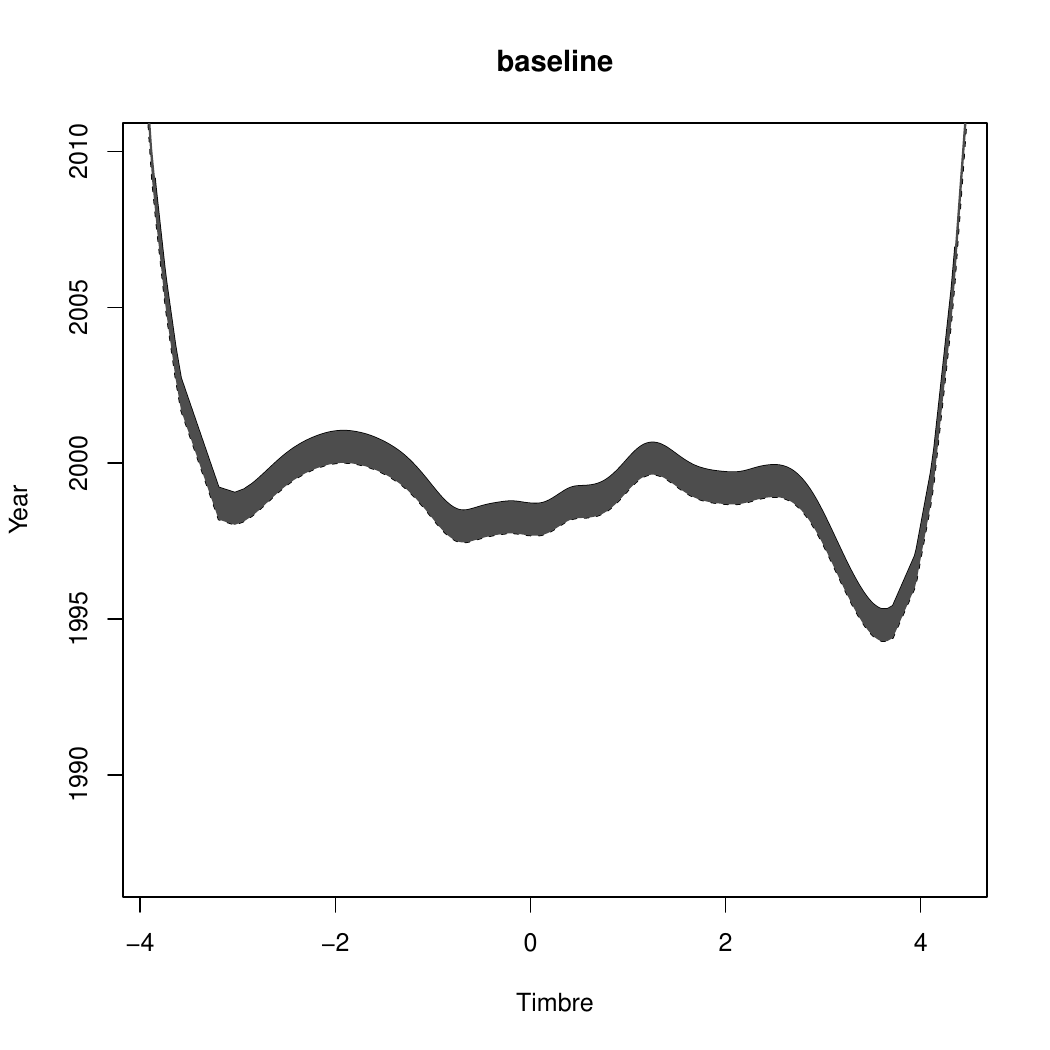}
		\end{minipage}
	\centering
	    \begin{minipage}[t]{0.45\linewidth}
		\centering
		\includegraphics[scale=0.3]{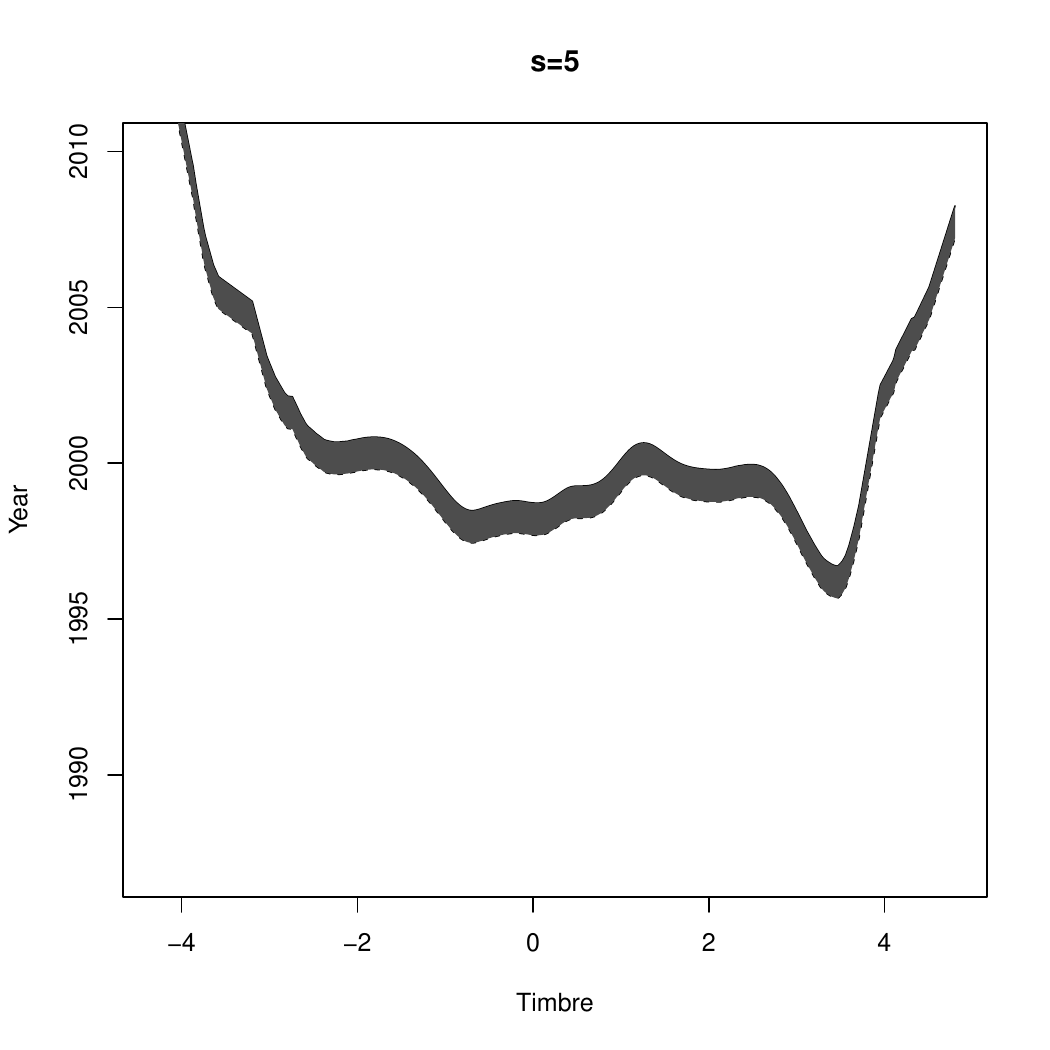}
		\end{minipage}	 
		\begin{minipage}[t]{0.45\linewidth}
		\centering
		\includegraphics[scale=0.3]{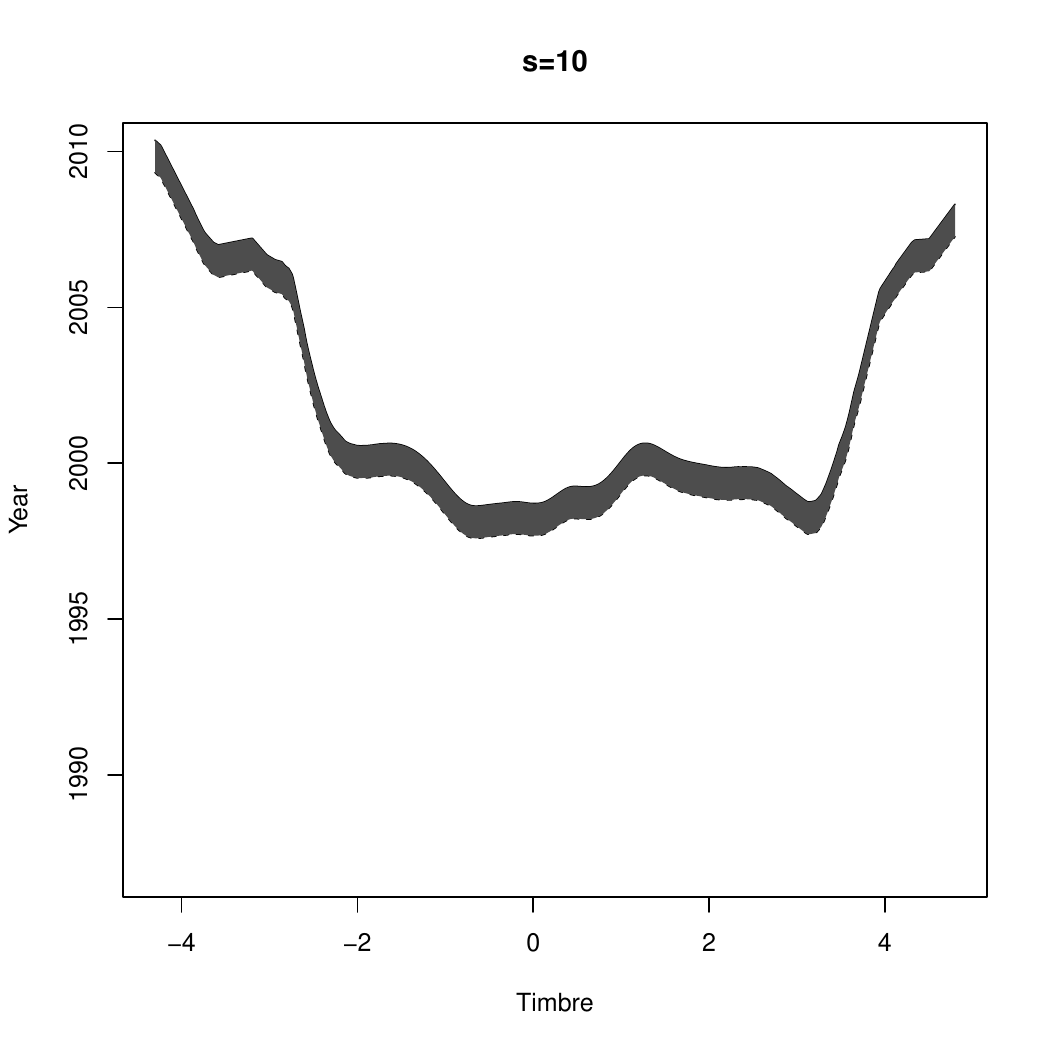}
		\end{minipage}
		\begin{minipage}[t]{0.45\linewidth}
		\centering
		\includegraphics[scale=0.3]{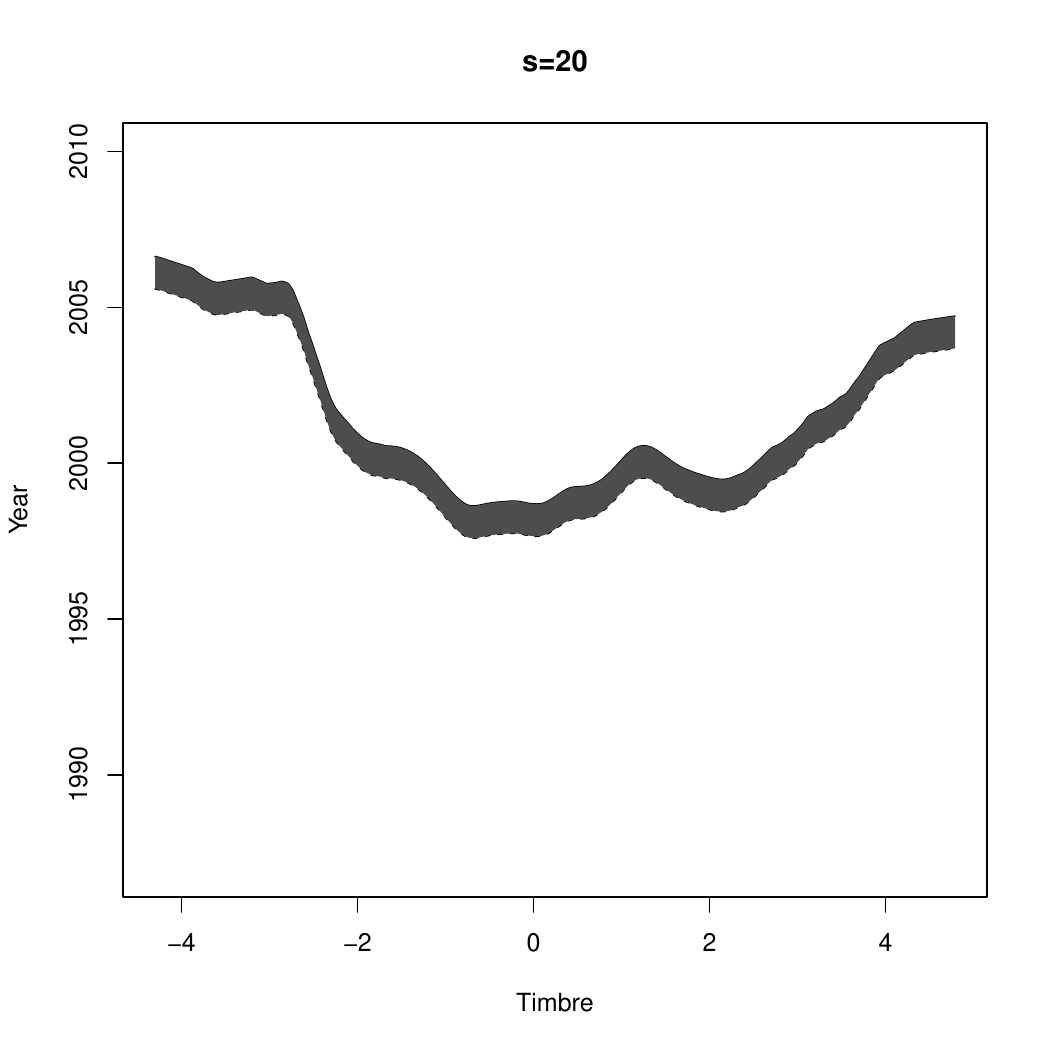}
		\end{minipage}
		\caption{\textit{\footnotesize 95\% Credible
		sets (grey areas) for $f$ based on a subset of 10,000 samples in Million Song Data. The first plot refers to the baseline method
		where the whole samples were used. The rest three plots refer to the aggregation procedure which was applied to 5, 10, 20 random splits.}}
\label{fig:MSDCR}
	\end{figure}

In the Bayesian setup, we applied our aggregation procedure to construct 95\% credible sets for $f$ based on a subset of $N=10,000$ songs released from the year 1996 to 2010.
We randomly split the observations to $s=5, 10, 20$ subsets. 
We also compared our results with the baseline method in which all ten thousand observations were used. 
Credible sets are displayed as gray areas in Figure \ref{fig:MSDCR}. 
We find that the shapes of all credible sets are overall the same when the timbre ranges from -4 to 4, e.g., all 
display a W-shape, although the results are a bit sensitive near the endpoints.
Therefore, the overall pattern of the sets appears to be
insensitive to the above selections of $s$. 
		
\subsection{Flight Delay Data}

\begin{figure}[h]
	\centering
		\begin{minipage}[t]{0.45\linewidth}
		\centering
		\includegraphics[scale=0.3]{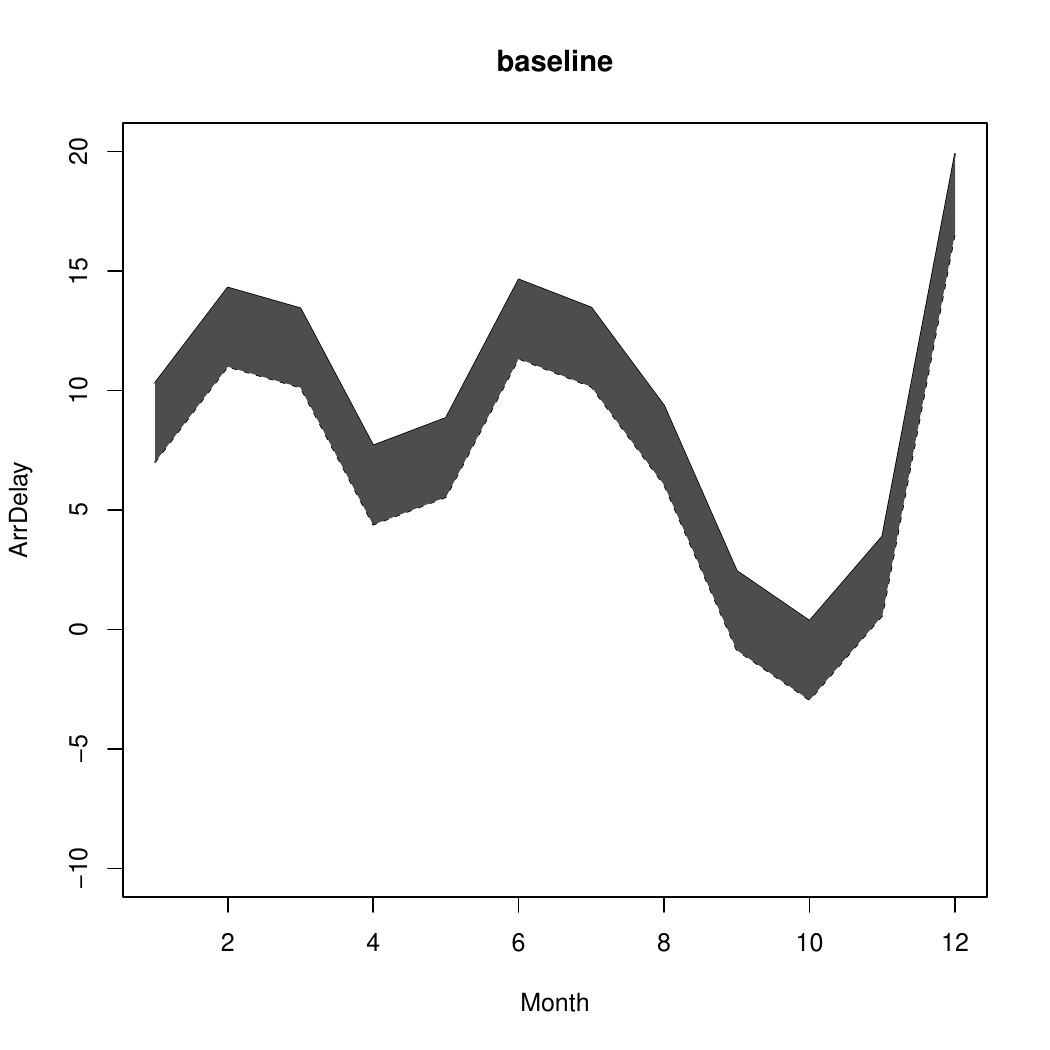}
		\end{minipage}
	\centering
	    \begin{minipage}[t]{0.45\linewidth}
		\centering
		\includegraphics[scale=0.3]{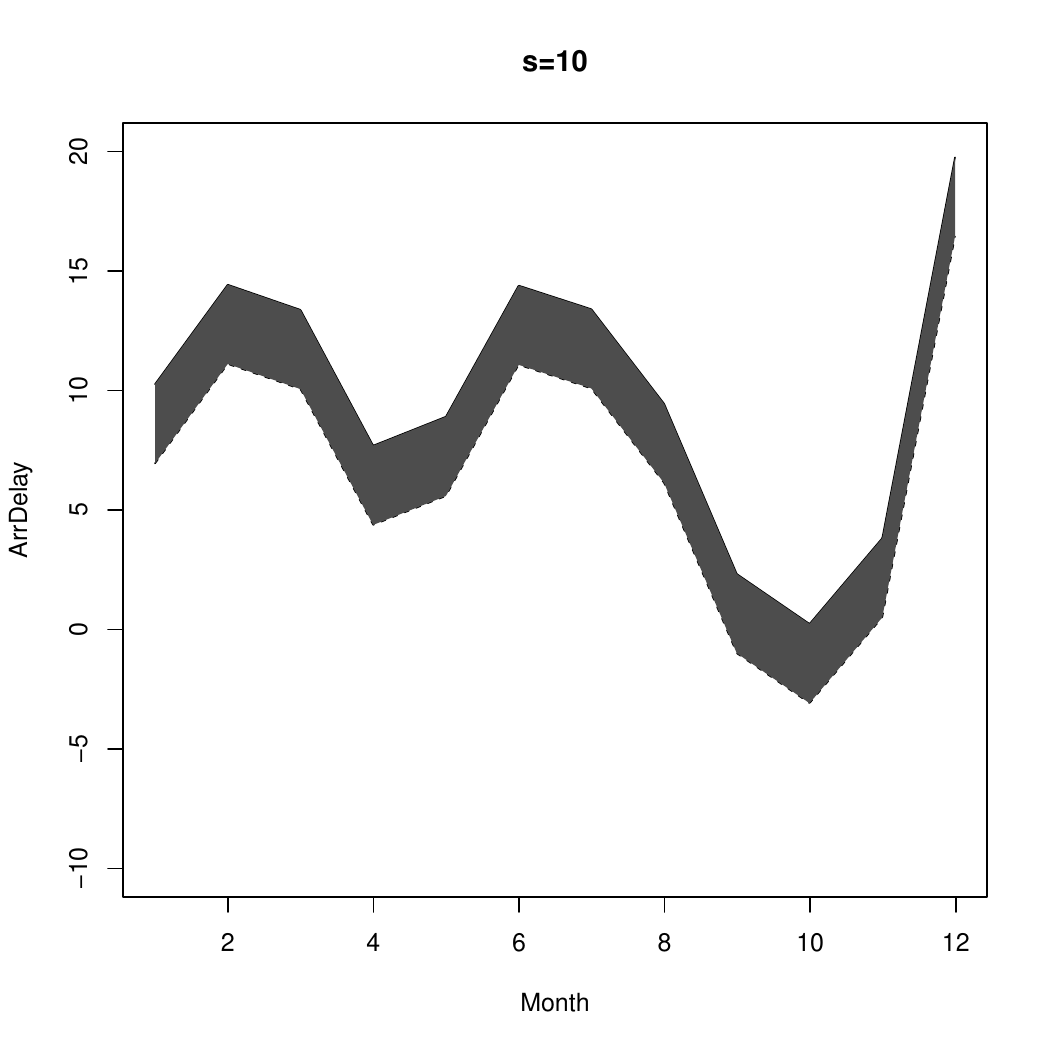}
		\end{minipage}	 
		\begin{minipage}[t]{0.45\linewidth}
		\centering
		\includegraphics[scale=0.3]{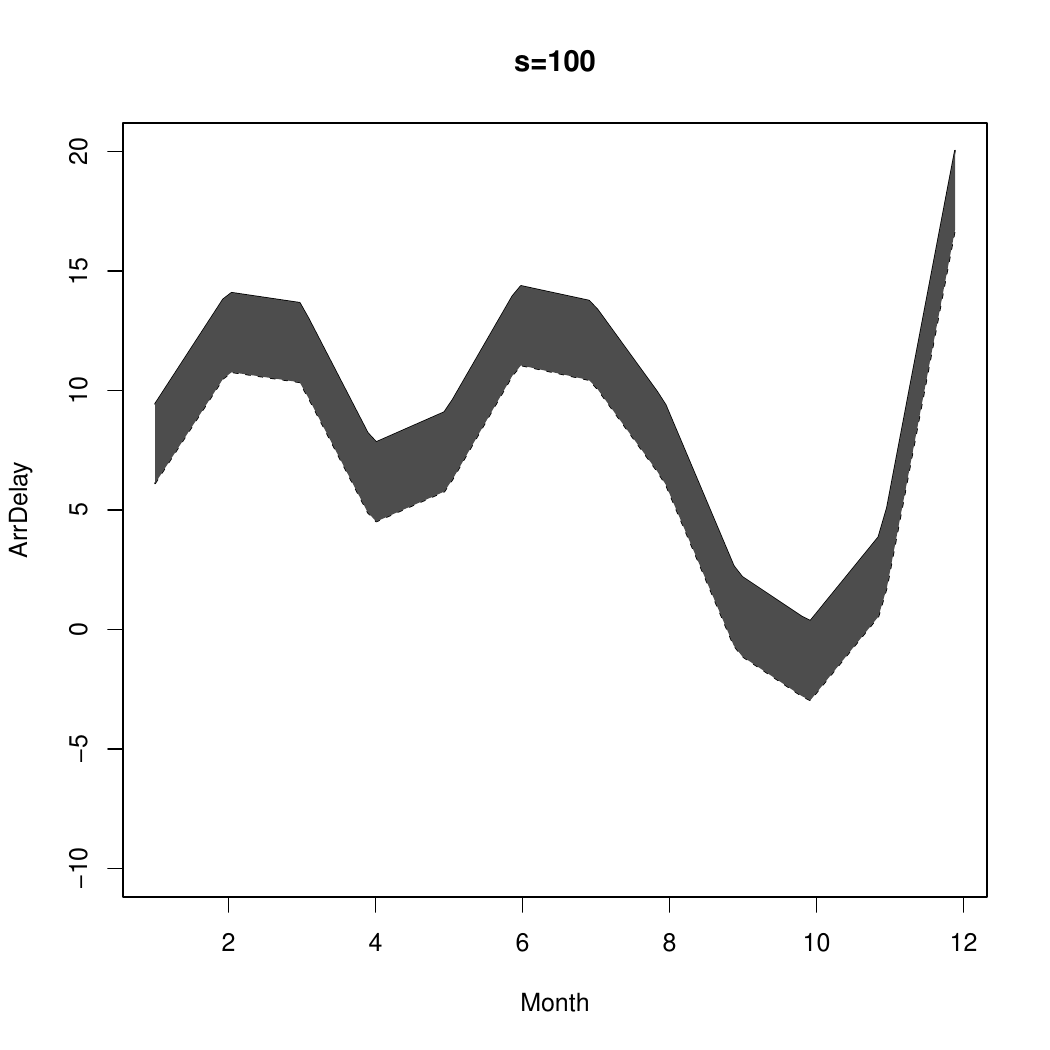}
		\end{minipage}
		\begin{minipage}[t]{0.45\linewidth}
		\centering
		\includegraphics[scale=0.3]{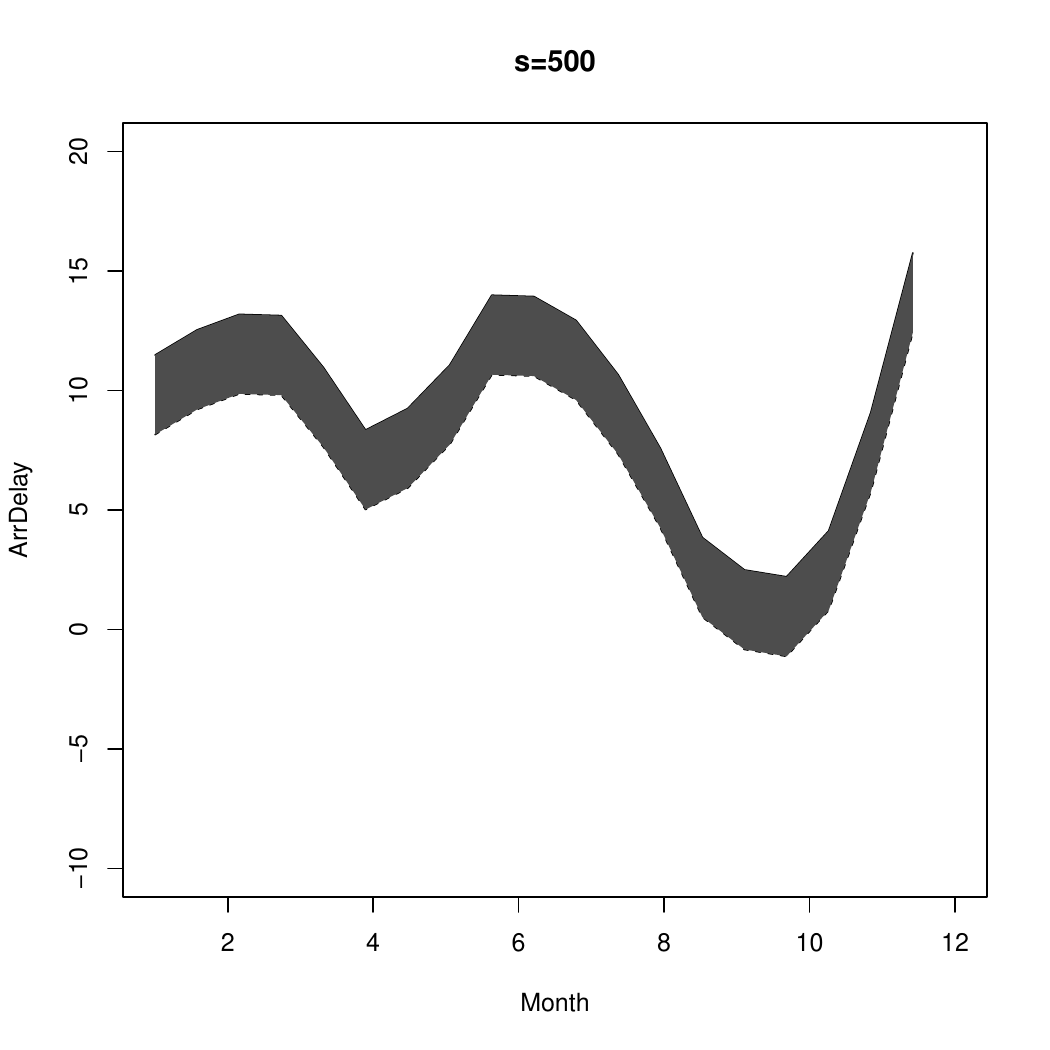}
		\end{minipage}
		\caption{\textit{\footnotesize 95\% Credible
		sets (grey areas) for $f$ based on a subset of 10,000 samples in Flight Delay Data. The first plot refers to the baseline 
		method where the whole samples were used. The rest three plots refer to the aggregation procedure which was applied to 10, 100, 500 random splits.}}
\label{fig:FDDCR}
	\end{figure}

We applied our aggregation procedure to one more real data set, the FDD. 
The data consists of flight arrival and departure information for all commercial flights within the United States, from October 1987 to April 2008. The main purpose is to find the key factors that have an impact on the flight delay. We considered the
relationship (denoted $f$) between month and the length of the flight delay, i.e.,
\textrm{$\textrm{length of flight delay}=f(\textrm{month})$+error}. Negative length of delay implies that the flight arrived earlier. We applied the same Bayesian aggregation procedure as described in MSD to a randomly
selected subset of $N=10,000$ flight information in the year 2007. We randomly split the observations to  $s=10, 100, 500$ subsamples, based on which the aggregated credible sets for $f$ were constructed. We also compared the results with the baseline where all the ten thousand samples were used. Credible sets are displayed as gray areas in Figure \ref{fig:FDDCR}. Again, the shapes of the four credible sets appear to be almost the same for all $s$.  

\subsection{Computation Efficiency}
\begin{figure}
	\centering
\includegraphics[scale=0.5]{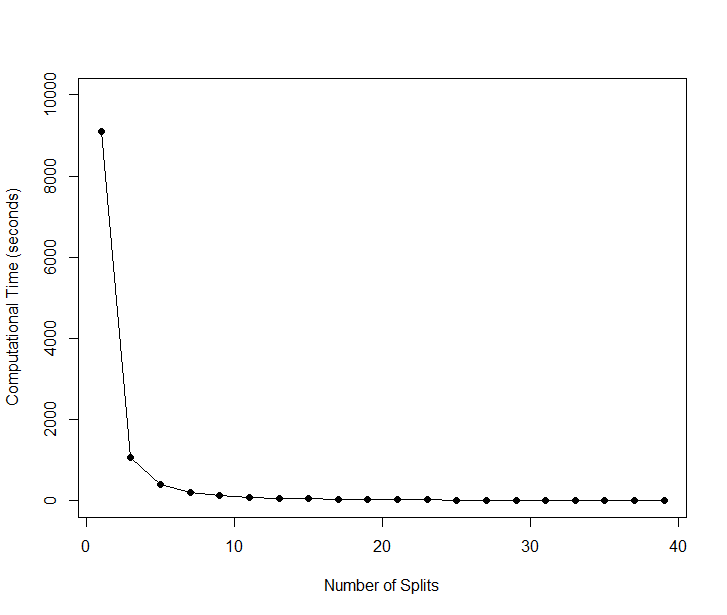}
	\caption{\textit{\footnotesize Computational time of aggregation procedures for MSD.}}
\label{fig:time_MSD}
	\end{figure}

\begin{figure}
	\centering
\includegraphics[scale=0.5]{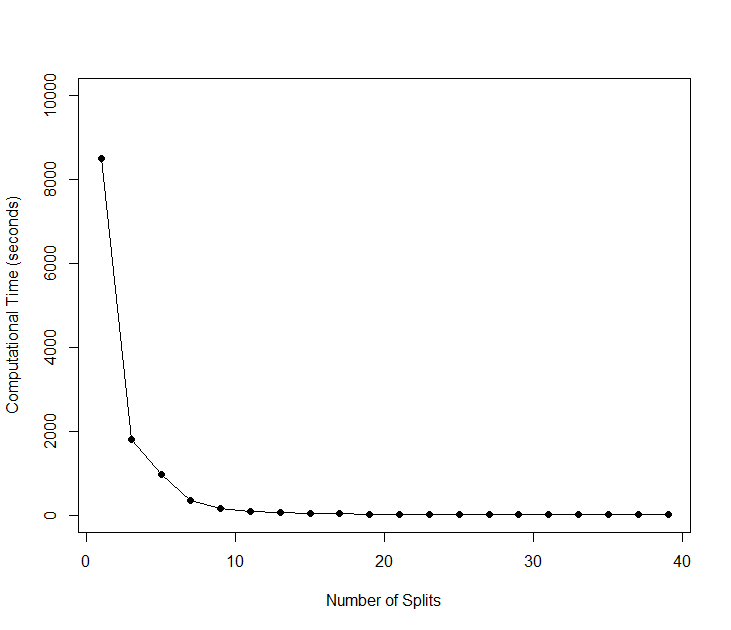}
	\caption{\textit{\footnotesize Computational time of aggregation procedures for FDD.}}
\label{fig:time_FDD}
	\end{figure}

We compare the overall execute computation time of both MSD and FDD on different numbers of splits, e.g. computational time per machine $\times$ number of machines in Figure \ref{fig:time_MSD}-\ref{fig:time_FDD}. It can be seen that the computing time dramatically decreases as the number of splits increases, which reflects the scalability of our proposed algorithm. 

\section{Conclusions}

This paper proposes algorithms for aggregating individual posterior results such as modes, balls, intervals, into their global counterparts.
The algorithms are easy-to-implement which are particularly useful in big data scenarios.
We also experimented the proposed algorithms through simulated and real data sets.
A notable contribution of this article is to provide rigorously justified theoretical guarantees.
The major tool for proving our theoretical results is a uniform Gaussian approximation theorem
which shows that the individual posterior distributions converge uniformly to Gaussian processes
provided that the number of subsets is not too large.

\section*{Acknowledgments}

Shang's research is sponsored by NSF DMS-1764280 and NSF DMS-1821157.
Cheng's research is sponsored by NSF (CAREER Award DMS-1151692, DMS-1418042) and Office of Naval Research (ONR N00014-15-1-2331).
The authors thank the Associate Editor Eric P. Xing and two anonymous referees for their helpful comments and suggestions
that significantly improve the quality of this paper.

\section{APPENDIX}

\setcounter{subsection}{0}
\renewcommand{\thesubsection}{A.\arabic{subsection}}
\setcounter{equation}{0}
\renewcommand{\theequation}{A.\arabic{equation}}
\setcounter{lemma}{0}
\renewcommand{\thelemma}{A.\arabic{lemma}}
\setcounter{proposition}{0}
\renewcommand{\theproposition}{A.\arabic{proposition}}
This appendix section contains the proofs of the main results. Section \ref{app:proofs:in:UBvM}
contains proof of Theorem \ref{uniform:bvm:thm} and relevant preliminary results.
Section \ref{app:proofs:in:aggr:post:mean} includes the proof of Theorem \ref{uniform:post:mean}.
Sections \ref{app:proofs:in:cr:strong} and \ref{sec:cr:weak} includes the proofs
of Theorems \ref{cp:cr:strong} and \ref{cp:cr:weak}, i.e., coverage properties of the credible
sets based on strong and weak topology respectively.

All proofs crucially depend on an eigensystem designed for simultaneous diagonalization of the two bilinear functionals $U,V$ induced from likelihood and prior, respectively. In fact, $(\varphi_\nu, \rho_\nu)$ is a solution of the following ordinary differential system (whose existence and uniqueness is guaranteed by \cite{Birk1908}):
\begin{eqnarray}\label{eigen:problem}
&&(-1)^m \varphi_\nu^{(2m)}(\cdot)=\rho_\nu \pi(\cdot)
\varphi_\nu(\cdot),\nonumber\\
&&\varphi_\nu^{(j)}(0)=\varphi_\nu^{(j)}(1)=0,\,\,\,\,
j=m,m+1,\ldots,2m-1,
\end{eqnarray}
Properties of this eigen-system are summarized in Proposition~\ref{A2}, whose proof can be found in \cite[Proposition 2.2]{SC13}.
\begin{proposition}\label{A2}
It holds that $\sup_{\nu\in\bbN}\|\varphi_\nu\|_\infty<\infty$,
and that the sequence $\rho_\nu$ is nondecreasing with
$\rho_1=\cdots=\rho_m=0$, and $\rho_\nu>0$ for $\mu>m$.
Moreover, $\rho_\nu\asymp \nu^{2m}$ and
\begin{equation}\label{regular}
V(\varphi_\mu,\varphi_\nu)=\delta_{\mu\nu},\,\,\,\,
J(\varphi_\mu,\varphi_\nu)=\rho_\mu\delta_{\mu\nu},\,\,\mu,\nu\in\bbN,
\end{equation}
where $\delta_{\mu\nu}$ is the Kronecker's delta. In particular, any
$f\in S^m(\mathbb{I})$ admits a Fourier expansion
$f=\sum_{\nu} V(f,\varphi_\nu)\varphi_\nu$ with convergence held in the $\|\cdot\|$-norm.
\end{proposition}

\subsection{Proofs in Section \ref{sec:UBvM}}\label{app:proofs:in:UBvM}
The proof of Theorem \ref{uniform:bvm:thm} requires the following technical result which  derives a local contraction rate $\widetilde r_n$ uniformly over $s$:
$\widetilde{r}_n=(nh/\log{2s})^{-1/2}+h^{m+\frac{\beta-1}{2}}$.
The proof can be found in (\cite{SCBigRate}).
\begin{Proposition}\label{refined:contraction:rate}
If $f_0$ satisfies Condition (\textbf{S}) and the following Rate Condition (\textbf{R}) holds:
\[
nh^{2m+1}\ge1,\,\,a_n=O(\widetilde{r}_n),\,\,
b_n\le1,\,\,r_n^2b_n\le\widetilde{r}_n^2.
\]
Let $a\ge0$ be a fixed constant. Then for any $\varepsilon\in (0,1)$, there exist positive
constants $M',N'$ s.t. for any $n\ge N'$,
\begin{equation}\label{refined:rate:eqn}
P_{f_0}\left(\max_{1\le j\le s}\{E\{\|f-f_0\|^aI(\|f-f_0\|\ge M'\widetilde{r}_n)|\textbf{D}_j\}\ge
M's^2\exp(-n\widetilde{r}_n^2/\log(2s))\right)\le\varepsilon
\end{equation}
\end{Proposition}
We remark that Proposition \ref{refined:contraction:rate} 
significantly generalizes the classical results in \cite{GGV00,VZ08}. 

\begin{proof}[Proof of Theorem \ref{uniform:bvm:thm}]
Let $M_1,M_2$ be large positive constants.
For any fixed constant $a\ge0$, consider three events:
\begin{eqnarray*}
\mathcal{E}_n'&=&\{\max_{1\le j\le s}\|\widehat{f}_{j,n}-f_0\|\le M_1\widetilde{r}_n\}\\
\mathcal{E}_n''&=&\{\max_{1\le j\le s}
E\{\|f-f_0\|^aI(\|f-f_0\|\ge M_2\widetilde{r}_n)|\textbf{D}_j\}\le M_2 s^2\exp(-n\widetilde{r}_n^2/\log(2s))\}\\
\mathcal{E}_n'''&=&\{\max_{1\le j\le s}E_{0j}\{\|f-f_0\|^aI(\|f-f_0\|\ge M_2\widetilde{r}_n)\}\le
M_2\exp(-n\widetilde{r}_n^2)\}
\end{eqnarray*}
where $E_{0j}$ means expectation taken under $P_{0j}$.
It follows from \cite{SCBigRate} and Proposition \ref{refined:contraction:rate} 
that we can choose $M_1>M_2$ (both large enough) 
s.t. $P_{f_0}(\mathcal{E}_n'\cap\mathcal{E}_n'')\ge 1-\varepsilon_1/2$
where $\varepsilon_1>0$ is an arbitrary constant.
Meanwhile, by (\cite{SCBigRate}) we have, on $\mathcal{E}_n'$, for any $1\le j\le s$,
\begin{eqnarray}\label{bvm:eqn1}
&&E_{0j}\{\|f-f_0\|^aI(\|f-f_0\|\ge M_2\widetilde{r}_n)\}\nonumber\\
&=&\frac{\int_{\|f-f_0\|\ge M_2\widetilde{r}_n}\|f-f_0\|^a
\exp\left(-\frac{n}{2}\|f-\widehat{f}_{j,n}\|^2\right)d\Pi(f)}
{\int_{S^m(\mathbb{I})}
\exp\left(-\frac{n}{2}\|f-\widehat{f}_{j,n}\|^2\right)d\Pi(f)}\nonumber\\
&\le&\frac{\int_{\|f-f_0\|\ge M_2\widetilde{r}_n}\|f-f_0\|^a
\exp\left(-\frac{n}{2}\|f-\widehat{f}_{j,n}\|^2\right)d\Pi(f)}
{\int_{\|f-f_0\|\le \widetilde{r}_n}
\exp\left(-\frac{n}{2}\|f-\widehat{f}_{j,n}\|^2\right)d\Pi(f)}\nonumber\\
&\le&\exp\left(-\left((M_2-M_1)^2/2-(M_1+1)^2/2-c_3/4\right)n\widetilde{r}_n^2\right)C(a,\Pi),
\end{eqnarray}
where $c_3>0$ is a universal constant and $C(a,\Pi)=\int_{S^m(\mathbb{I})}\|f-f_0\|^ad\Pi(f)$.
We can choose $M_2>C(a,\Pi)$ 
so that the quantity (\ref{bvm:eqn1})
is less than $M_2\exp(-n\widetilde{r}_n^2)$.
So $\mathcal{E}_n'$ implies $\mathcal{E}_n'''$,
so that $P_{f_0}(\mathcal{E}_n''')\ge P_{f_0}(\mathcal{E}_n'\cap\mathcal{E}_n'')\ge 1-\varepsilon_1/2$.
Define
$\mathcal{E}_n=\mathcal{E}_n'\cap\mathcal{E}_n''\cap\mathcal{E}_n'''$,
then it can be seen that $P_{f_0}(\mathcal{E}_n)\ge 1-\varepsilon_1$.

Let $T_{j}$ be defined as
\begin{eqnarray}\label{Tj1:Tj2}
T_{j2}(f)&=&-\frac{1}{2n}\sum_{i\in I_j}
[(\Delta f)(X_i)^2-E_X\{(\Delta f)(X)^2\}].
\end{eqnarray}
Following Lemma \ref{pre:lemma:5}, for any $1\le j\le s$,
\begin{equation}\label{bvm:eqn2}
\ell_{jn}(f)-\ell_{jn}(\widehat{f}_{j,n})+\frac{1}{2}\|f-\widehat{f}_{j,n}\|^2=T_{j}(f).
\end{equation}
It follows from the proof of Proposition \ref{refined:contraction:rate} that on $\mathcal{E}_n$,
for any $f\in S^m(\mathbb{I})$ satisfying $\|f-f_0\|\le M_2\widetilde{r}_n$ and $1\le j\le s$,
\begin{equation}\label{bvm:eqn3}
|T_{j}(f)|\le D\times\widetilde{r}_n^2b_n,
\end{equation}
where $D=D(M_1,M_2)$
is a positive constant
depending only on $M_1,M_2$.
Recall that our assumption says that $\varepsilon_2
\equiv nD\widetilde{r}_n^2b_n=o(1)$.

For $1\le j\le s$, define
\begin{eqnarray*}
&& J_{nj1}=\int_{S^m(\mathbb{I})}\exp\left(n(\ell_{jn}(f)-\ell_{jn}(\widehat{f}_{j,n}))\right)d\Pi(f),\\
&& J_{nj2}=\int_{S^m(\mathbb{I})}\exp\left(-\frac{n}{2}\|f-\widehat{f}_{j,n}\|^2\right)d\Pi(f),
\end{eqnarray*}
\begin{eqnarray*}
&& \bar{J}_{nj1}=\int_{\|f-f_0\|\le M_2\widetilde{r}_n}
\exp\left(n(\ell_{jn}(f)-\ell_{jn}(\widehat{f}_{j,n}))\right)d\Pi(f),\\
&& \bar{J}_{nj2}=\int_{\|f-f_0\|\le M_2\widetilde{r}_n}
\exp\left(-\frac{n}{2}\|f-\widehat{f}_{j,n}\|^2\right)d\Pi(f).
\end{eqnarray*}
For simplicity, let $\varepsilon_3=M_2s^2\exp(-n\widetilde{r}_n^2/\log(2s))$.
On $\mathcal{E}_n$ (with $a=0$) and for any $1\le j\le s$,
\[
0\le\frac{J_{nj1}-\bar{J}_{nj1}}{J_{nj1}}\le M_2s^2\exp(-n\widetilde{r}_n^2/\log(2s))
=\varepsilon_3,\,\,\,\,
0\le\frac{J_{nj2}-\bar{J}_{nj2}}{J_{nj2}}\le\exp(-n\widetilde{r}_n^2)\le\varepsilon_3.
\]
By some algebra, it can be shown that the above inequalities lead to
\begin{equation}\label{bvm:eqn4}
(1-\varepsilon_3)\cdot\frac{\bar{J}_{nj2}}{\bar{J}_{nj1}}\le
\frac{J_{nj2}}{J_{nj1}}\le\frac{1}{1-\varepsilon_3}\cdot\frac{\bar{J}_{nj2}}{\bar{J}_{nj1}}.
\end{equation}

Meanwhile, on $\mathcal{E}_n$ and for any $1\le j\le s$,
using (\ref{bvm:eqn3}) and the elementary inequality $|\exp(x)-1|\le 2|x|$ for $|x|\le\log{2}$, we get that
\begin{eqnarray*}
|\bar{J}_{nj2}-\bar{J}_{nj1}|&\le&\int_{\|f-f_0\|\le M_2\widetilde{r}_n}\exp\left(-\frac{n}{2}\|f-\widehat{f}_{j,n}\|^2\right)\times |\exp(nT_{j}(f))-1| d\Pi(f)\\
&\le&2\varepsilon_2\bar{J}_{nj2},
\end{eqnarray*}
leading to that
\begin{equation}\label{bvm:eqn5}
\frac{1}{1+2\varepsilon_2}\le\frac{\bar{J}_{nj2}}{\bar{J}_{nj1}}\le\frac{1}{1-2\varepsilon_2}.
\end{equation}
Combining (\ref{bvm:eqn4}) and (\ref{bvm:eqn5}), on $\mathcal{E}_n$ and for any $1\le j\le s$,
$\frac{1-\varepsilon_3}{1+2\varepsilon_2}\le
\frac{J_{nj2}}{J_{nj1}}\le\frac{1}{(1-2\varepsilon_2)(1-\varepsilon_3)}$.
When $n$ is large, $\varepsilon_3\le\varepsilon_2$ and both quantities are small, the above
inequalities lead to
\begin{equation}\label{bvm:eqn6}
-4\varepsilon_2\le\frac{1-\varepsilon_3}{1+2\varepsilon_2}-1\le
\frac{J_{nj2}}{J_{nj1}}-1\le \frac{1}{(1-2\varepsilon_2)(1-\varepsilon_3)}-1\le 4\varepsilon_2
\end{equation}

For simplicity, denote $R_{nj}(f)=nT_{j}(f)$.
For any $S\in\mathcal{S}$, let $S'=S\cap\{f\in S^m(\mathbb{I}): \|f-f_0\|\le M_2\widetilde{r}_n\}$.
Then on $\mathcal{E}_n$, we get that
$\max_{1\le j\le s}|P(S|\textbf{D}_j)-P_{0j}(S)|\le\max_{1\le j\le s}|P(S'|\textbf{D}_j)-P_{0j}(S')|+2\varepsilon_3$.
Moreover, it follows from (\ref{bvm:eqn6}) that on $\mathcal{E}_n$ and for any $1\le j\le s$,
\begin{eqnarray*}
&&|P(S'|\textbf{D}_j)-P_{0j}(S')|\\
&=&\bigg|\int_{S'}\left(\frac{\exp(n(\ell_{jn}(f)-\ell_{jn}(\widehat{f}_{j,n})))}{J_{nj1}}
-\frac{\exp\left(-\frac{n}{2}\|f-\widehat{f}_{j,n}\|^2\right)}{J_{nj2}}\right)d\Pi(f)\bigg|\\
&\le&\int_{S'}\exp\left(-\frac{n}{2}\|f-\widehat{f}_{j,n}\|^2\right)\times
\bigg|\frac{\exp(R_{nj}(f))}{J_{nj1}}-\frac{1}{J_{nj2}}\bigg|d\Pi(f)\\
&\le&\int_{S'}\exp\left(-\frac{n}{2}\|f-\widehat{f}_{j,n}\|^2\right)
\times\frac{|\exp(R_{nj}(f))-1|}{J_{nj2}}d\Pi(f)\\
&&+\int_{S'}\exp\left(-\frac{n}{2}\|f-\widehat{f}_{j,n}\|^2\right)
\times\exp(R_{nj}(f))\times\bigg|\frac{1}{J_{nj1}}-\frac{1}{J_{nj2}}\bigg|d\Pi(f)\\
&\le&2\varepsilon_2\frac{\int_{S'}\exp\left(-\frac{n}{2}\|f-\widehat{f}_{j,n}\|^2\right)d\Pi(f)}{J_{nj2}}\\
&&+\exp(\varepsilon_2)\times\bigg|\frac{1}{J_{nj1}}-\frac{1}{J_{nj2}}\bigg|
\times\int_{S'}\exp\left(-\frac{n}{2}\|f-\widehat{f}_{j,n}\|^2\right)d\Pi(f)\\
&\le&2\varepsilon_2+\exp(\varepsilon_2)\times\bigg|\frac{J_{nj2}}{J_{nj1}}-1\bigg|\le
2\varepsilon_2+4\varepsilon_2\exp(\varepsilon_2)\le14\varepsilon_2.
\end{eqnarray*}
Note that the right hand side is free of $S$.
Then we get that on $\mathcal{E}_n$,
$\sup_{S\in\mathcal{S}}
\max_{1\le j\le s}|P(S|\textbf{D}_j)-P_{0j}(S)|\le 14\varepsilon_2+2\varepsilon_3\le 16\varepsilon_2$.
This implies that for sufficiently large $n$,
\begin{eqnarray*}
&&P_{f_0}\left(\sup_{S\in\mathcal{S}}
\max_{1\le j\le s}|P(S|\textbf{D}_j)-P_{0j}(S)|>16\varepsilon_2\right)\\
&\le& P_{f_0}(\mathcal{E}_n^c)+P_{f_0}\left(\mathcal{E}_n,
\sup_{S\in\mathcal{S}}
\max_{1\le j\le s}|P(S|\textbf{D}_j)-P_{0j}(S)|>16\varepsilon_2\right)
=P_{f_0}(\mathcal{E}_n^c)\le\varepsilon_1.
\end{eqnarray*}
The desirable result follows by the simple fact 
$\varepsilon_2\lesssim \sqrt{s}N^{-\frac{4m^2+2m\beta-10m+1}{4m(2m+\beta)}}(\log{N})^{\frac{5}{2}}$
when $h\asymp h^*$.
\end{proof}

\subsection{Proofs in Section \ref{sec:aggr:post:mean}}\label{app:proofs:in:aggr:post:mean}

\begin{proof}[Proof of Theorem \ref{uniform:post:mean}]
We first show (\ref{uniform:post:mean:eqn}).
Let $A_n=\{f\in S^m(\mathbb{I}): \|f-f_0\|\ge M\widetilde{r}_n\}$
and $B_j=\{f\in S^m(\mathbb{I}): dP(f|\textbf{D}_j)\ge dP_{0j}(f)\}$
for $1\le j\le s$.
By Proposition \ref{refined:contraction:rate}, Theorem \ref{uniform:bvm:thm} and (\ref{bvm:eqn1}) with $a=1$ therein, we can choose
$M>0$ sufficiently large such that
\begin{eqnarray*}
&&\max_{1\le j\le s}\|E(f|\textbf{D}_j)-E_{0j}(f)\|\\
&=&\max_{1\le j\le s}\|\int (f-f_0)dP(f|\textbf{D}_j)-\int (f-f_0)dP_{0j}(f)\|\\
&\le&\max_{1\le j\le s}\|\int_{A_n}(f-f_0)dP(f|\textbf{D}_j)\|
+\max_{1\le j\le s}\|\int_{A_n}(f-f_0)dP_{0j}(f)\|\\
&&+\max_{1\le j\le s}\|\int_{A_n^c}(f-f_0)(dP(f|\textbf{D}_j)-dP_{0j}(f))\|\\
&\le&\max_{1\le j\le s}E\{\|f-f_0\|I(f\in A_n)|\textbf{D}_j\}
+\max_{1\le j\le s}E_{0j}\{\|f-f_0\|I(f\in A_n)\}\\
&&+M\widetilde{r}_n\max_{1\le j\le s}\int_{A_n^c}|dP(f|\textbf{D}_j)-dP_{0j}(f)|\\
&=&O_{P_{f_0}}\left(s^2\exp(-n\widetilde{r}_n^2/\log(2s))+\exp(-n\widetilde{r}_n^2)+\widetilde{r}_n\sqrt{s}N^{-\frac{4m^2+2m\beta-10m+1}{4m(2m+\beta)}}(\log{N})^{\frac{5}{2}}\right)\\
&=&O_{P_{f_0}}\left(\widetilde{r}_n\sqrt{s}N^{-\frac{4m^2+2m\beta-10m+1}{4m(2m+\beta)}}(\log{N})^{\frac{5}{2}}\right)
\equiv O_{P_{f_0}}(L_N),
\end{eqnarray*}
where the second last equality uses Theorem \ref{uniform:bvm:thm} and the fact that, uniformly for $j$,
\begin{eqnarray*}
&&\int_{A_n^c}|dP(f|\textbf{D}_j)-dP_{0j}(f)|\\
&=&|P(A_n^c\cap B_j|\textbf{D}_j)-P_{0j}(A_n^c\cap B_j)|+|P(A_n^c\cap B_j^c|\textbf{D}_j)-P_{0j}(A_n^c\cap B_j^c)|.
\end{eqnarray*}
Then (\ref{uniform:post:mean:eqn}) follows from the trivial fact that 
$E_{0j}\{f\}=E(W^j|\textbf{D}_j)=\widetilde{f}_{j,n}$.

Next we show (\ref{big:rate:eqn}).
By direct examinations we can verify the following Rate Conditions ($\textbf{R}$):
\[
n\widetilde{r}_n^2b_n=o(1),
N\widetilde{r}_{N}^2b_{N}=o(1), Nh^{1/2}a_N^2=o(1), Nh^{1/2}a_n^2=o(1).
\]
Define $Rem_{j,n}=\widehat{f}_{j,n}-f_0-S_{j,n}(f_0)$ for $j=1,2,\ldots,s$.
It follows by Lemma \ref{pre:lemma:2} of \cite{SCBigRate} that
$\max_{1\le j\le s}\|Rem_{j,n}\|=O_{P_{f_0}}(a_n)$.

It is easy to see that $a_{N,\nu}/a_{n,\nu}\le s$ for all $\nu\ge1$.
Then it holds from (\ref{bar:r:N:1}) that
\begin{eqnarray}\label{breve:vs:tilde}
\|\breve{f}_{N,\lambda}-\widetilde{f}_{N,\lambda}\|^2&=&\sum_{\nu\ge1}\left(\frac{a_{N,\nu}}{a_{n,\nu}}\right)^2
V\left(\frac{1}{s}\sum_{j=1}^s(\breve{f}_{j,n}-\widetilde{f}_{j,n}),\varphi_\nu\right)^2(1+\lambda\rho_\nu)\nonumber\\
&\le&s^2\|\frac{1}{s}\sum_{j=1}^s(\breve{f}_{j,n}-\widetilde{f}_{j,b})\|^2
=O_{P_{f_0}}\left(s^2L_N^2\right)=o_{P_{f_0}}(N^{-1}h^{-1/2}).
\end{eqnarray}
The last equality owes to 
the condition $s^4\log(2s)=o\left(N^{\frac{4m^2+2m\beta-11m+1}{2m(2m+\beta)}}(\log{N})^{-5}\right)$
and $\beta>3/2$.

By direct examinations, we have
\begin{eqnarray}\label{an:important:eqn:cr:strong}
\widetilde{f}_{N,\lambda}-f_0&=&\sum_{\nu=1}^\infty\left(a_{N,\nu}
\left(\frac{1}{s}\sum_{j=1}^sV(\widehat{f}_{j,n},\varphi_\nu)\right)-f_\nu^0\right)\varphi_\nu\nonumber\\
&=&\sum_{\nu=1}^\infty
\left(a_{N,\nu}\left(\frac{1}{s}\sum_{j=1}^sV(Rem_{j,n}+f_0+S_{j,n}(f_0),\varphi_\nu)\right)
-f_\nu^0\right)\varphi_\nu\nonumber\\
&=&\sum_{\nu=1}^\infty a_{N,\nu} V(\frac{1}{s}\sum_{j=1}^s Rem_{j,n},\varphi_\nu)\varphi_\nu+\sum_{\nu=1}^\infty
(a_{N,\nu}-1)f_\nu^0\varphi_\nu\nonumber\\
&&+\sum_{\nu=1}^\infty a_{N,\nu} V(\frac{1}{N}\sum_{i=1}^N\epsilon_i K_{X_i},\varphi_\nu)\varphi_\nu
-\sum_{\nu=1}^\infty a_{N,\nu} V(\mathcal{P}_\lambda f_0,\varphi_\nu)\varphi_\nu.
\end{eqnarray}
Denote the four terms in the above equation by $T_1,T_2,T_3,T_4$.

Since $a_{N,\nu}\le1$, it is easy to see that
\begin{eqnarray}\label{cp:cr:strong:eqn3.5}
\|T_1\|_2^2&=&\sum_{\nu=1}^\infty a_{N,\nu}^2 |V(\frac{1}{s}\sum_{j=1}^s Rem_{j,n},\varphi_\nu)|^2\nonumber\\
&\le&\sum_{\nu=1}^\infty |V(\frac{1}{s}\sum_{j=1}^s Rem_{j,n},\varphi_\nu)|^2=\|\frac{1}{s}\sum_{j=1}^s Rem_{j,n}\|_2^2\le (\max_{1\le j\le s}\|Rem_{j,n}\|)^2=O_{P_{f_0}}(a_n^2).\nonumber\\
\end{eqnarray}

Using $h\asymp N^{-1/(2m+\beta)}$ and a direct algebra we get that
\begin{eqnarray*}
\|T_2\|_2^2&=&\sum_{\nu=1}^\infty(a_{N,\nu}-1)^2|f_\nu^0|^2
\asymp\sum_{\nu=1}^\infty\left(\frac{\nu^{2m+\beta}}{\nu^{2m+\beta}+N(1+\lambda\nu^{2m})}\right)^2|f_\nu^0|^2=
o(N^{-\frac{2m+\beta-1}{2m+\beta}})=o(N^{-1}h^{-1}).
\end{eqnarray*}
Meanwhile, it follows by Proposition \cite{SCBigRate} that
\begin{eqnarray*}
\|T_4\|_2^2&=&\sum_{\nu=1}^\infty a_{N,\nu}^2 |f_\nu^0|^2\left(\frac{\lambda\rho_\nu}{1+\lambda\rho_\nu}\right)^2
\le\sum_{\nu=1}^\infty|f_\nu^0|^2\left(\frac{\lambda\rho_\nu}{1+\lambda\rho_\nu}\right)^2\\
&\lesssim&\sum_{\nu=1}^\infty|f_\nu^0|^2(h\nu)^{2m+\beta-1}\frac{(h\nu)^{2m-\beta+1}}{(1+(h\nu)^{2m})^2}
=o(N^{-\frac{2m+\beta-1}{2m+\beta}})=o(N^{-1}h^{-1}).
\end{eqnarray*}

Define $R(x,x')=\sum_{\nu=1}^\infty a_{N,\nu}\frac{\varphi_\nu(x)\varphi_\nu(x')}{1+\lambda\rho_\nu}$
for any $x,x'\in\mathbb{I}$. Also define $R_x(\cdot)=R(x,\cdot)$. It is easy to see that
$R_x\in S^m(\mathbb{I})$ for any $x\in\mathbb{I}$.
Then it can be shown that $T_3=\frac{1}{N}\sum_{i=1}^N\epsilon_i R_{X_i}$, leading to
\[
\|T_3\|_2^2=V(T_3,T_3)=\frac{1}{N^2}\sum_{i=1}^N \epsilon_i^2 V(R_{X_i},R_{X_i})
+\frac{2}{N^2}\sum_{i<k}\epsilon_i\epsilon_k V(R_{X_i},R_{X_k}).
\]
Since $E_{f_0}\{\epsilon^2V(R_X,R_X)\}=O(h^{-1})$, we have
$E_{f_0}\{\|T_3\|_2^2\}=O(N^{-1}h^{-1})$.
Therefore, $\|\widetilde{f}_{N,\lambda}-f_0\|_2^2=
O_{P_{f_0}}(N^{-1}h^{-1})=O_{P_{f_0}}\left(N^{-\frac{2m+\beta-1}{2m+\beta}}\right)$.
This together with (\ref{breve:vs:tilde}) leads to (\ref{big:rate:eqn}). 
\end{proof}

\subsection{Proofs in Section \ref{sec:cr:strong}}\label{app:proofs:in:cr:strong}

Before proving Theorem \ref{cp:cr:strong}, we give some preliminary notation and results. Define an ``oracle''
penalized likelihood
$\ell_{N,\lambda}(f)=-\frac{1}{2N}
\sum_{i=1}^N(Y_i-f(X_i))^2-\frac{\lambda}{2}J(f)$.
Applying Theorem~\ref{uniform:bvm:thm} to $s=1$, we have
\begin{equation}\label{full:bvm:expression}
\sup_{S\in\mathcal{S}}|P(S|\textbf{D})-P_{0}(S)|=o_{P_{f_0}}(1),
\end{equation}
where
$P_{0}(S)=\frac{\int_S\exp\left(-\frac{N}{2}\|f-\widehat{f}^{or}_{N,\lambda}\|^2\right)d\Pi(f)}
{\int_{S^m(\mathbb{I})}\exp\left(-\frac{N}{2}\|f-\widehat{f}^{or}_{N,\lambda}\|^2\right)d\Pi(f)}$
and $\widehat{f}^{or}_{N,\lambda}=\arg\max_{f\in S^m(\mathbb{I})}\ell_{N,\lambda}(f)$
is the ``oracle" smoothing spline estimator based on full data. Consider a generalized Fourier expansion of
$\widehat{f}^{or}_{N,\lambda}$: 
$\widehat{f}^{or}_{N,\lambda}(\cdot)=\sum_{\nu=1}^\infty V(\widehat{f}^{or}_{N,\lambda},\varphi_\nu)\varphi_\nu(\cdot)$.
By Theorem 5.2 in \cite{SC14}, we have $P_0(S)=P(W^{or}\in S|\textbf{D})$ for any
$S\in\mathcal{S}$, where $W^{or}(\cdot)=\sum_{\nu=1}^\infty (a_{N,\nu}V(\widehat{f}^{or}_{N,\lambda},
\varphi_\nu)+b_{N,\nu}\tau_\nu v_\nu)\varphi_\nu(\cdot)$. 
Here, $a_{n,\nu}$ $b_{n,\nu}$ are analogous to ones in the definition of $W^j(\cdot)$ in Section \ref{sec:UBvM},
and $v_\nu\sim N(0, \tau_\nu^{-2})$ and $\tau_\nu^2$ are given in (\ref{gamnu}). Define the mean functions of $W^{or}$ 
as
$\widetilde{f}^{or}_{N,\lambda}(\cdot):=\sum_{\nu=1}^\infty a_{N,\nu}V(\widehat{f}^{or}_{N,\lambda},\varphi_\nu)\varphi_\nu(\cdot)$.
So we can re-express $W^{or}$ as
$W^{or}=\widetilde{f}^{or}_{N,\lambda}+W_N$,
where $W_N(\cdot):=\sum_{\nu=1}^\infty b_{N,\nu}\tau_\nu v_\nu\varphi_\nu(\cdot)$ is a zero-mean GP.

The following result describes the distribution of $W_n$ and $W_N$.
\begin{lemma}\label{a:limit:distribution}
As $N\rightarrow\infty$,
$\frac{n\|W_n\|_2^2-\zeta_{1,n}}{\sqrt{2\zeta_{2,n}}}\overset{d}{\longrightarrow}N(0,1),\,\,\,\,
\textrm{and}\,\,\,\,\frac{N\|W_N\|_2^2-\zeta_{1,N}}{\sqrt{2\zeta_{2,N}}}\overset{d}{\longrightarrow}N(0,1)$.
\end{lemma}

\begin{proof}[Proof of Theorem \ref{cp:cr:strong}]
We can show that Rate Conditions ($\textbf{R}$) hold by direct calculations.

It is sufficient to investigate the $P_{f_0}$-probability of the event
$\{\|f_0-\breve{f}_{N,\lambda}\|_2\le r_N(\alpha)\}$.
To achieve this goal, we first prove the following fact:
\begin{equation}\label{cp:cr:strong:eqn1}
\max_{1\le j\le s}|z_{j,n}(\alpha)-z_\alpha|=o_{P_{f_0}}(1),
\end{equation}
where $z_\alpha=\Phi^{-1}(1-\alpha)$ and $\Phi$ is the c.d.f. of $N(0,1)$,
and $z_{j,n}(\alpha)=(nr_{j,n}(\alpha)^2-\zeta_{1,n})/\sqrt{2\zeta_{2,n}}$.
The proof of the theorem follows by (\ref{cp:cr:strong:eqn1}) and a careful analysis
of $f_0-\breve{f}_{N,\lambda}$.

We first show (\ref{cp:cr:strong:eqn1}).
It follows by Theorem \ref{uniform:bvm:thm} that for any $j=1,2,\ldots,s$,
\begin{eqnarray*}
|P(R_{j,n}(\alpha)|\textbf{D}_j)-P_{0j}(R_{j,n}(\alpha))|
&\le&\max_{1\le k\le s}|P(R_{j,n}(\alpha)|\textbf{D}_k)-P_{0k}(R_{j,n}(\alpha))|\\
&\le&\sup_{S\in\mathcal{S}}\max_{1\le k\le s}|P(S|\textbf{D}_k)-P_{0k}(S)|
=o_{P_{f_0}}(1).
\end{eqnarray*}
Together with $P(R_{j,n}(\alpha)|\textbf{D}_j)=1-\alpha$, we have
$\max_{1\le j\le s}|P_{0j}(R_{j,n}(\alpha))-(1-\alpha)|=o_{P_{f_0}}(1)$.
Let $\Delta_j=\breve{f}_{j,n}-\widetilde{f}_{j,n}$ for $1\le j\le s$.
It is clear that
\begin{eqnarray}\label{cp:cr:strong:eqn:1}
P_{0j}(R_{j,n}(\alpha))
&=&P(W^j\in R_{j,n}(\alpha)|\textbf{D}_j)=
P(\|W_n+\Delta_j\|_2\le r_{j,n}(\alpha)|\textbf{D}_j)\nonumber\\
&=&P(\|W_n\|_2^2+2\langle W_n,\Delta_j\rangle_2+\|\Delta_j\|_2^2\le r_{j,n}(\alpha)^2|\textbf{D}_j),
\end{eqnarray}
and, for any $\varepsilon\in(0,1)$,
\begin{eqnarray}\label{cp:cr:strong:eqn:2}
&&P(|\langle W_n,\Delta_j\rangle_2|^2\ge \|\Delta_j\|_2^2/(n\varepsilon)|\textbf{D}_j)
\le n\varepsilon E\{|\langle W_n,\Delta_j\rangle_2|^2|\textbf{D}_j\}/\|\Delta_j\|_2^2\nonumber\\
&=&\frac{n\varepsilon}{\|\Delta_j\|_2^2}\sum_{\nu\ge1}b_{n,\nu}^2|V(\Delta_j,\varphi_\nu)|^2\le 
\frac{n\varepsilon}{\|\Delta_j\|_2^2}\times\frac{\|\Delta_j\|_2^2}{n}=\varepsilon,
\end{eqnarray}
and by Theorem \ref{uniform:post:mean},
$\max_{1\le j\le s}\|\Delta_j\|_2^2=O_{P_{f_0}}(L_N^2)$,
where $L_N=\widetilde{r}_n\sqrt{s}N^{-\frac{4m^2+2m\beta-10m+1}{4m(2m+\beta)}}(\log{N})^{\frac{5}{2}}$.
By (\ref{s:for:aggregation}), $\zeta_{k,n}\asymp n^{1/(2m+\beta)}$ (Lemma \ref{a:limit:distribution}),
and direct examinations it holds that
\begin{equation}\label{cp:cr:strong:eqn:3}
\max_{1\le j\le s}\frac{n\|\Delta_j\|_2^2}{\sqrt{\zeta_{2,n}}}=o_{P_{f_0}}(1).
\end{equation}
Combining (\ref{cp:cr:strong:eqn:1}) and (\ref{cp:cr:strong:eqn:2}) we get that
\begin{eqnarray*}
P_{0j}(R_{j,n}(\alpha))&\ge&\Phi_n\left(z_{j,n}(\alpha)-\frac{n\|\Delta_j\|_2^2}{\sqrt{\zeta_{2,n}}}
-\frac{2n\|\Delta_j\|_2}{\sqrt{n\varepsilon\zeta_{2,n}}}\right)-\varepsilon,\\
P_{0j}(R_{j,n}(\alpha))&\le&\Phi_n\left(z_{j,n}(\alpha)-\frac{n\|\Delta_j\|_2^2}{\sqrt{\zeta_{2,n}}}
+\frac{2n\|\Delta_j\|_2}{\sqrt{n\varepsilon\zeta_{2,n}}}\right)+\varepsilon,
\end{eqnarray*}
where $\Phi_n$ is the c.d.f. of $U_n$.
It follows by Lemma \ref{a:limit:distribution} and 
Polya's theorem (\cite{CT88}) that $\Phi_n$ uniformly converges to $\Phi(\cdot)$,
the c.d.f. of standard normal variable. 
Therefore, when $n$ becomes large enough, 
\[
\bigg|\Phi_n\left(z_{j,n}(\alpha)-\frac{n\|\Delta_j\|_2^2}{\sqrt{\zeta_{2,n}}}
-\frac{2n\|\Delta_j\|_2}{\sqrt{n\varepsilon\zeta_{2,n}}}\right)
-\Phi\left(z_{j,n}(\alpha)-\frac{n\|\Delta_j\|_2^2}{\sqrt{\zeta_{2,n}}}
-\frac{2n\|\Delta_j\|_2}{\sqrt{n\varepsilon\zeta_{2,n}}}\right)\bigg|\le\varepsilon,
\]
\[
\bigg|\Phi_n\left(z_{j,n}(\alpha)-\frac{n\|\Delta_j\|_2^2}{\sqrt{\zeta_{2,n}}}
+\frac{2n\|\Delta_j\|_2}{\sqrt{n\varepsilon\zeta_{2,n}}}\right)
-\Phi\left(z_{j,n}(\alpha)-\frac{n\|\Delta_j\|_2^2}{\sqrt{\zeta_{2,n}}}
+\frac{2n\|\Delta_j\|_2}{\sqrt{n\varepsilon\zeta_{2,n}}}\right)\bigg|\le\varepsilon,
\]
where implies that
\[
\Phi\left(z_{j,n}(\alpha)-\frac{n\|\Delta_j\|_2^2}{\sqrt{\zeta_{2,n}}}-
\frac{2n\|\Delta_j\|_2}{\sqrt{n\varepsilon\zeta_{2,n}}}\right)\le P_{0j}(R_{j,n}(\alpha))+2\varepsilon
=\Phi(z_\alpha)+2\varepsilon+o_{P_{f_0}}(1),
\]
\[
\Phi\left(z_{j,n}(\alpha)-\frac{n\|\Delta_j\|_2^2}{\sqrt{\zeta_{2,n}}}+
\frac{2n\|\Delta_j\|_2}{\sqrt{n\varepsilon\zeta_{2,n}}}\right)\ge P_{0j}(R_{j,n}(\alpha))-2\varepsilon
=\Phi(z_\alpha)-2\varepsilon+o_{P_{f_0}}(1).
\]
Since (\ref{cp:cr:strong:eqn:3}) implies that
$\frac{n\|\Delta_j\|_2^2}{\sqrt{\zeta_{2,n}}}$ and
$\frac{2\sqrt{n}\|\Delta_j\|_2}{\sqrt{\zeta_{2,n}}}$ are both 
$o_{P_{f_0}}(1)$ uniformly for $j$, so (\ref{cp:cr:strong:eqn1}) holds.

Next we prove the theorem. Consider expansion (\ref{an:important:eqn:cr:strong}).
Only focus on $T_3$.
Define $W(N)=2\sum_{1\le i<k\le N}\epsilon_i\epsilon_k V(R_{X_i},R_{X_k})$.
Let $W_{ik}=2\epsilon_i\epsilon_k V(R_{X_i},R_{X_k})$,
then $W(N)=\sum_{1\le i<k\le N}W_{ik}$. Note that $W(N)$
is clean in the sense of \cite{deJong87}.
Let $\sigma^2(N)=E_{f_0}\{W(N)^2\}$ and
$G_I$, $G_{II}$, $G_{IV}$ be defined as
$G_I=\sum_{i<j} E_{f_0}\{W_{ij}^4\}$,
$G_{II}=\sum_{i<j<k} (E_{f_0}\{W_{ij}^2W_{ik}^2\}+E_{f_0}\{W_{ji}^2W_{jk}^2\}+
E_{f_0}\{W_{ki}^2W_{kj}^2\})$, and
\begin{eqnarray*}
G_{IV}&=&\sum_{i<j<k<l} (E_{f_0}\{W_{ij} W_{ik} W_{lj} W_{lk}\}+
E_{f_0}\{W_{ij} W_{il} W_{kj} W_{kl}\}
+ E_{f_0}\{W_{ik}W_{il}W_{jk}W_{jl}\}).
\end{eqnarray*}
Since $\varphi_\nu$ are uniformly bounded,
we get that
$\|R_x\|_2^2=\sum_{\nu=1}^\infty\frac{|\varphi_\nu(x)|^2}{(1+N^{-1}\tau_\nu^2+\lambda\rho_\nu)^2}
\lesssim h^{-1}$,
where ``$\lesssim$" is free of $x$. This implies that
$G_I=O(N^2h^{-4})$ and $G_{II}=O(N^3 h^{-4})$.

It can also be shown that for pairwise distinct $i,k,t,l$,
\begin{eqnarray*}
E_{f_0}\{W_{ik} W_{il} W_{tk} W_{tl}\}
&=&2^4E_{f_0}\{\epsilon_i^2\epsilon_k^2\epsilon_t^2\epsilon_l^2V(R_{X_i},R_{X_k})V(R_{X_i},R_{X_l})V(R_{X_t},R_{X_{k}})
V(R_{X_t},R_{X_l})\}\\
&=&2^4\sum_{\nu=1}^\infty\frac{a_{N,\nu}^8}{(1+\lambda\rho_\nu)^8}=O(h^{-1}),
\end{eqnarray*}
which implies that $G_{IV}=O(N^4h^{-1})$.
In the mean time, a straight algebra leads to that
\begin{eqnarray*}
\sigma^2(N)&=&4{N \choose 2}\sum_{\nu=1}^\infty\frac{a_{N,\nu}^4}{(1+\lambda\rho_\nu)^4}
=4{N\choose 2}\sum_{\nu=1}^\infty\left(\frac{N}{\tau_\nu^2+N(1+\lambda\rho_\nu)}\right)^4=2N(N-1)\zeta_{4,N}
\asymp N^2h^{-1}.
\end{eqnarray*}
Since $Nh^2\asymp N^{1-2/(2m+\beta)}\rightarrow\infty$,
we get that $G_I,G_{II}$ and $G_{IV}$ are all of order $o(\sigma^4(N))$.
Then it follows by \cite{deJong87} that as $N\rightarrow\infty$,
$\frac{W(N)}{N\sqrt{2\zeta_{4,N}}}\overset{d}{\longrightarrow}N(0,1)$.
Since $\zeta_{4,N}\asymp h^{-1}$, the above equation leads to that $W(N)/N=O_{P_{f_0}}(h^{-1/2})$.

It follows by direct examination that
$Var_{f_0}\{\sum_{i=1}^N \epsilon_i^2 V(R_{X_i},R_{X_i})\}\le NE_{f_0}\{\epsilon_i^4\|R_{X_i}\|_2^4\}
=O(Nh^{-2})$,
leading to that
$\sum_{i=1}^N \epsilon_i^2 V(R_{X_i},R_{X_i})=
E_{f_0}\{\sum_{i=1}^N \epsilon_i^2 V(R_{X_i},R_{X_i})\}+O_{P_{f_0}}(N^{1/2}h^{-1})
=N\zeta_{2,N}+O_{P_{f_0}}(N^{1/2}h^{-1})$.
Therefore, it follows by Rate Condition (\textbf{R}),
i.e., $Nha_n^2=o(1)$, and the analysis on $T_1,T_2,T_3,T_4$ in (\ref{an:important:eqn:cr:strong}) 
that
\begin{eqnarray}\label{cp:cr:strong:eqn6}
Nh\|\widetilde{f}_{N,\lambda}-f_0\|_2^2&=&Nh\|T_3\|_2^2+O_{P_{f_0}}(Nha_n^2)+o_{P_{f_0}}(1)
=h\zeta_{2,N}+o_{P_{f_0}}(1).
\end{eqnarray}

In the end, note from (\ref{cp:cr:strong:eqn1}) and $\zeta_{k,n}\asymp n^{\alpha_1}$ for $\alpha_1=1/(2m+\beta)$
(see proof of Lemma \ref{a:limit:distribution}) that
$\frac{n}{s}\sum_{j=1}^s r_{j,n}(\alpha)^2=
\zeta_{1,n}+\sqrt{2\zeta_{2,n}}z_\alpha+o_{P_{f_0}}(\sqrt{\zeta_{2,n}})$,
which leads to that
\begin{equation}\label{cp:cr:strong:eqn4}
Nr_N(\alpha)^2=\zeta_{1,N}+\sqrt{2\zeta_{2,N}}z_\alpha+o_{P_{f_0}}(h^{-1/2}).
\end{equation}
Therefore,
$Nhr_N(\alpha)^2=h\zeta_{1,N}(1+o_{P_{f_0}}(1))$.
Since $\lim\inf_{N\rightarrow\infty}(h\zeta_{1,N}-h\zeta_{2,N})>0$,
we get by (\ref{breve:vs:tilde}) that, 
with $P_{f_0}$-probability approaching one, $\|\breve{f}_{N,\lambda}-f_0\|_2\le r_N(\alpha)$.
Meanwhile, it follows by \cite{SCBigRate} that
$\|\widehat{f}^{or}_{N,\lambda}-f_0-S_{N,\lambda}(f_0)\|_2=O_{P_{f_0}}(a_N)$
and $\|\frac{1}{s}\sum_{j=1}^s\widehat{f}_{j,n}-f_0-\frac{1}{s}\sum_{j=1}^sS_{j,n}(f_0)\|_2=O_{P_{f_0}}(a_n)$,
where $S_{N,\lambda}(f_0)=\frac{1}{N}\sum_{i=1}^N\epsilon_i K_{X_i}-\mathcal{P}_\lambda f_0$.
Note that $S_{N,\lambda}(f_0)=\frac{1}{s}\sum_{j=1}^sS_{j,n}(f_0)$,
which leads to $\|\widehat{f}^{or}_{N,\lambda}-\frac{1}{s}\sum_{j=1}^s\widehat{f}_{j,n}\|_2=O_{P_{f_0}}(a_n+a_N)$.
Since $a_{N,\nu}\le 1$, we get that
\begin{eqnarray}\label{cp:cr:strong:eqn5}
N\|\widetilde{f}^{or}_{N,\lambda}-\widetilde{f}_{N,\lambda}\|^2&=&
N\sum_{\nu=1}^\infty a_{N,\nu}^2V\left(\widehat{f}^{or}_{N,\lambda}-\frac{1}{s}\sum_{j=1}^s\widehat{f}_{j,n},
\varphi_\nu\right)^2(1+\lambda\rho_\nu)\nonumber\\
&\le&N\sum_{\nu=1}^\infty V\left(\widehat{f}^{or}_{N,\lambda}-\frac{1}{s}\sum_{j=1}^s\widehat{f}_{j,n},\varphi_\nu
\right)^2(1+\lambda\rho_\nu)
\nonumber\\
&=&N\|\widehat{f}^{or}_{N,\lambda}-\frac{1}{s}\sum_{j=1}^s\widehat{f}_{j,n}\|^2
=O_{P_{f_0}}(Na_n^2+Na_N^2)\\
&=&o_{P_{f_0}}(h^{-1/2}),\,\,\,\,(\textrm{by condition $Nh^{1/2}a_n^2+Nh^{1/2}a_N^2=o(1)$})\nonumber
\end{eqnarray}
Using (\ref{breve:vs:tilde}) we get that 
$N\|\widetilde{f}^{or}_{N,\lambda}-\breve{f}_{N,\lambda}\|_2^2=o_{P_{f_0}}(h^{-1/2})$.
Since 
$E\{|\langle W_N,\widetilde{f}^{or}_{N,\lambda}-\breve{f}_{N,\lambda}\rangle_2|^2|\textbf{D}\}
=\sum_{\nu\ge1}b_{N,\nu}^2V(\widetilde{f}^{or}_{N,\lambda}-\breve{f}_{N,\lambda},\varphi_\nu)^2
\le \|\widetilde{f}^{or}_{N,\lambda}-\breve{f}_{N,\lambda}\|_2^2/N=o_{P_{f_0}}(N^{-2}h^{-1/2})$,
we have that
$N\|W^{or}-\breve{f}_{N,\lambda}\|_2^2=N\|W_N\|_2^2+o_{P_{f_0}}(h^{-1/2})$.
It follows by
$P\left(\frac{N\|W_N\|_2^2-\zeta_{1,N}}{\sqrt{2\zeta_{2,N}}}\le z_\alpha\right)\rightarrow 1-\alpha$,
(\ref{full:bvm:expression}) and (\ref{cp:cr:strong:eqn4}) that
$P(R_N(\alpha)|\textbf{D})=1-\alpha+o_{P_{f_0}}(1)$.
This completes the proof.
\end{proof}

\subsection{Proofs in Section \ref{sec:cr:weak}}\label{app:proofs:in:cr:weak}
Before proving Theorem \ref{cp:cr:weak}, let us present a preliminary lemma.

\begin{lemma}\label{a:limit:distribution:weak}
As $N\rightarrow\infty$,
$N\|W_N\|_\omega^2\overset{d}{\rightarrow}\sum_{\nu=1}^\infty d_\nu\eta_\nu^2,
\,\,\,\,\textrm{and}\,\,\,\,
n\|W_n\|_\omega^2\overset{d}{\rightarrow}\sum_{\nu=1}^\infty d_\nu\eta_\nu^2$,
where $\eta_\nu$ are independent standard normal random variables.
\end{lemma}

\begin{proof}[Proof of Theorem \ref{cp:cr:weak}]
By direct examinations,
one can show that Rate Conditions ($\textbf{R}'$):
$n\widetilde{r}_n^2b_n=o(1)$,
$N\widetilde{r}_{N}^2b_N=o(1)$, $Na_N^2=o(1)$ and $Na_n^2=o(1)$
are all satisfied.

We first have the following fact:
\begin{equation}\label{cp:cr:weak:eqn1}
\max_{1\le j\le s}|\sqrt{n}r_{\omega,j,n}(\alpha)-\sqrt{c_\alpha}|=o_{P_{f_0}}(1),
\end{equation}
where $c_\alpha>0$ satisfies $P(\sum_{\nu=1}^\infty d_\nu\eta_\nu^2\le c_\alpha)=1-\alpha$
with $\eta_\nu$ being independent standard normal random variables.
It follows from (\ref{cp:cr:weak:eqn1}) that
\begin{equation}\label{cp:cr:weak:eqn2}
Nr_{\omega,N}(\alpha)^2=c_\alpha+o_{P_{f_0}}(1).
\end{equation}
By Theorem \ref{uniform:post:mean} and the condition 
$s=o(N^{\frac{4m^2+2m\beta-12m+1}{8m(2m+\beta)}}(\log{N})^{-\frac{3}{2}})$
we have the following 
$\max_{1\le j\le s}n\|\Delta_j\|_\omega^2=\max_{1\le j\le s}n\|\Delta_j\|_2^2=O_{P_{f_0}}(nL_N^2)=o_{P_{f_0}}(1)$.
Also, for arbitrarily small $\varepsilon\in(0,1)$,
$P(|\langle W_n,\Delta_j\rangle_\omega|^2\ge\|\Delta_j\|_\omega^2/(n\varepsilon)|\textbf{D}_j)\le\varepsilon$.
The proof of (\ref{cp:cr:weak:eqn1}) is then similar to the proof of (\ref{cp:cr:strong:eqn1})
and details are omitted. 

Let $T_1,T_2,T_3,T_4$ be defined in (\ref{an:important:eqn:cr:strong}).
It follows from the proof of Theorem \ref{cp:cr:strong} that
$\|T_1\|_\omega^2\le \|T_1\|_2^2=O_{P_{f_0}}(a_n^2)$,
so $N\|T_1\|_\omega^2=O_{P_{f_0}}(Na_n^2)=o_{P_{f_0}}(1)$ due to the condition $Na_n^2=o(1)$.
It follows by condition $h\asymp N^{-1/(2m+\beta)}$, dominated convergence theorem
and direct examinations,
\begin{eqnarray*}
\|T_2\|_\omega^2&=&\sum_{\nu=1}^\infty d_\nu (a_{N,\nu}-1)^2 |f_\nu^0|^2
\asymp N^{-2}\sum_{\nu=1}^\infty d_\nu\frac{\nu^{2m+\beta+1}}{(1+(h\nu)^{2m}+(h\nu)^{2m+\beta})^2}\times
\nu^{2m+\beta-1}|f_\nu^0|^2\\
&\lesssim& N^{-1}\sum_{\nu=1}^\infty\frac{(h\nu)^{2m+\beta+1}}{(1+(h\nu)^{2m}+(h\nu)^{2m+\beta})^2}\times
\nu^{2m+\beta-1}|f_\nu^0|^2=o(N^{-1}),
\end{eqnarray*}
and
\begin{eqnarray*}
\|T_4\|_\omega^2&=&\sum_{\nu=1}^\infty d_\nu a_{N,\nu}^2\left(\frac{\lambda\rho_\nu}{1+\lambda\rho_\nu}\right)^2
|f_\nu^0|^2
\lesssim\sum_{\nu=1}^\infty d_\nu \frac{(h\nu)^{2m-\beta+1}}{(1+(h\nu)^{2m}+(h\nu)^{2m+\beta})^2}
\times |f_\nu^0|^2 (h\nu)^{2m+\beta-1}\\
&\lesssim&h^{2m+\beta}\sum_{\nu=1}^\infty\frac{(h\nu)^{2m-\beta}}{(1+(h\nu)^{2m}+(h\nu)^{2m+\beta})^2}
\times|f_\nu^0|^2\nu^{2m+\beta-1}=o(N^{-1}).
\end{eqnarray*}

By direct examination it can be shown that
$T_3=\frac{1}{N}\sum_{i=1}^N\epsilon_i
\sum_{\nu=1}^\infty\frac{\varphi_\nu(X_i)\varphi_\nu}{1+\lambda\rho_\nu+N^{-1}\tau_\nu^2}$.
It follows by \cite{SC14} that, as $N\rightarrow\infty$,
$N\|T_3\|_\omega^2\overset{d}{\rightarrow}\sum_{\nu=1}^\infty d_\nu\eta_\nu^2$.
By the above analysis on $T_1$ through $T_4$, 
and $N\|\breve{f}_{N,\lambda}-\widetilde{f}_{N,\lambda}\|_\omega^2=O_{P_{f_0}}(Ns^2L_N^2)=o_{P_{f_0}}(1)$,
we get that
$N\|\breve{f}_{N,\lambda}-f_0\|_\omega^2\overset{d}{\rightarrow}\sum_{\nu=1}^\infty d_\nu\eta_\nu^2$.
It follows by (\ref{cp:cr:weak:eqn2}) that $\lim_{N\rightarrow\infty}P_{f_0}(f_0\in R_N^\omega(\alpha))=1-\alpha$.

It follows by $N\|\widetilde{f}^{or}_{N,\lambda}-\widetilde{f}_{N,\lambda}\|_2^2
=O_{P_{f_0}}(Na_N^2+Na_n^2)=o_{P_{f_0}}(1)$
(see (\ref{cp:cr:strong:eqn5})), $P(N\|W_N\|_2^2\le c_\alpha)\rightarrow 1-\alpha$, (\ref{cp:cr:weak:eqn2})
and (\ref{full:bvm:expression}) that $P(R_N^\omega(\alpha)|\textbf{D})=1-\alpha+o_{P_{f_0}}(1)$.
Proof is completed.
\end{proof}

\subsection{Computational Details}\label{subsection:computational_details}
In this subsection, we provide some computational details relating to Section \ref{sec:toy:agg:proc}. For convenience, we rewrite model \eqref{sec:toy:np:bayes:model} as following:
\begin{equation*}
Y_{ji} = f(X_{ji}) + \epsilon_{ji}, \ j=1,\ldots,s, \ i=1,\ldots, n.
\end{equation*}

\emph{Calculation of posterior means.} In order to calculate the posterior mean $\breve{f}_{j,n}$, we have to generate samples of $f$ from its posterior distribution $P(f|\{Y_{ji}, X_{ji}\}_{i=1}^n)$. In practice, directly sampling the function $f$ from $P(f|\{Y_{ji}, X_{ji}\}_{i=1}^n)$ is impossible. Instead, we generate some samples from $(f(X_{j1}), \ldots, f(X_{jn}))^{\top}$. As $n$ is large, $(f(X_{j1}), \ldots, f(X_{jn}))^{\top}$ can represent the whole curve of $f$. Firstly, let us derive the posterior distribution for $(f(X_{j1}), \ldots, f(X_{jn}))^{\top}$. For the $j$-th subsample, the likelihood function is written by
\begin{equation*}
Y_{j1}, \ldots, Y_{jn} | X_{j1}, \ldots, X_{jn}\sim N((f(X_{j1}), \ldots, f(X_{jn}))^{\top}, I_n).
\end{equation*}
Since $f$ follows a GP prior with mean zero and covariance function $K_0$, where $K_0$ is given in \eqref{eqn:kernal}, the prior of  $(f(X_{j1}), \ldots, f(X_{jn}))^{\top}$ is multivariate Gaussian:
\begin{equation*}
(f(X_{j1}), \ldots, f(X_{jn}))^{\top} \sim N(0, K_j),
\end{equation*}
where $K_j$ is the covariance matrix satisfying
\[
K_j =
\begin{bmatrix}
K_0(X_{j1}, X_{j1}),     & \cdots & K_0(X_{j1}, X_{jn}) \\
 \vdots   & \ddots & \vdots\\
K_0(X_{jn}, X_{j1})    & \cdots & K_0(X_{jn}, X_{jn})
\end{bmatrix}.
\]
$K_0(x, x')$ involves an infinite summation which is practically infeasible. Instead, the infinite sum is approximated by a finite one, i.e.,
\begin{equation*}
K_0(x, x') \approx 2\sum_{k=1}^M\frac{\cos(2\pi k(x-x'))}{(2\pi k)^{2m+\beta}+n\lambda(2\pi k)^{2m}}.
\end{equation*} 
In our numerical study, we found that  $M=100$ can already provide a good approximation. Due to the conjugacy, the posterior distribution of $(f(X_{j1}), \ldots, f(X_{jn}))^{\top}$ also follows a multivariate Gaussian distribution 
\begin{equation*}\label{eqn:posterior_Gau}
(f(X_{j1}), \ldots, f(X_{jn}))^{\top} \Big| \{Y_{ji}, X_{ji}\}_{i=1}^n\sim N\Big(K_j(K_j+\frac{1}{n}I_n)^{-1}(Y_{j1}, \ldots, Y_{jn})^{\top}, K_j(K_j+\frac{1}{n}I_n)^{-1}\frac{1}{n}\Big).
\end{equation*}
Next we generate $M$ independent samples, denoted $(f^{(l)}(X_{j1}), \ldots, f^{(l)}(X_{jn}))^{\top}, l =1,\ldots,M,$  from above multivariate Gaussian
distribution. Therefore, the posterior mean can be approximated by 
$$
\Big(\breve{f}_{jn}(X_{j1}), \ldots, \breve{f}_{jn}(X_{jn})\Big)^{\top} = \Big(\frac{1}{M}\sum_{l=1}^Mf^{(l)}(X_{j1}), \ldots, \frac{1}{M}\sum_{l=1}^Mf^{(l)}(X_{jn})\Big)^{\top}.
$$

\emph{Calculation of posterior radius.} Once we have $M$ independent samples $\{(f^{(l)}(X_{j1}), \ldots,f^{(l)}(X_{jn}) )\}_{l=1}^M$, we are able to approximate $\|f^{(l)}-\breve{f}_{j,n}\|_{L^2}$ by
\begin{equation*}
L_l = \Big(\frac{1}{n}\sum_{i=1}^n(f^{(l)}(X_{ji}) -\breve{f}_{jn}(X_{ji}))^2\Big)^{\tfrac{1}{2}}, \ \text{for} \ l=1,\ldots, M.
\end{equation*}
Finally, the radius $r_{j,n}(\alpha)$ is approximated by the upper $\alpha$-th percentile of $\{L_1, \ldots, L_M\}$. 

\emph{Calculation of the integral.} We approximate \eqref{eqn:coeff_intergral} by
\begin{equation*}
\breve{f}_{j,n,k}\approx\frac{\sqrt{2}}{n}\sum_{i=1}^n \breve{f}_{j,n}(X_{ji})\cos(2\pi kX_{ji})dx, \ \breve{g}_{j,n,k}\approx\frac{\sqrt{2}}{n}\sum_{i=1}^n\breve{f}_{j,n}(X_{ji})\sin(2\pi kX_{ji})dx.
\end{equation*}
In \eqref{eqn:radius_weight}, $C_k$ and $D_k$ also involve two integrals. Since they are independent of samples, any numerical method for integral calculation is applicable. We also approximate $\breve{f}_{N,\lambda}(x)$ in \eqref{special:aggregation:mean} by 
\begin{equation*}
\breve{f}_{N,\lambda}(x)\approx\sum_{k=1}^M w_{s,N,\lambda,k}
\left\{\breve{f}_{N,\lambda,k}\sqrt{2}\cos(2\pi kx)
+\breve{g}_{N,\lambda,k}\sqrt{2}\sin(2\pi kx)\right\}.
\end{equation*}

\newpage

\setcounter{page}{1}
\begin{frontmatter}
\begin{center}
\textit{Supplementary document to}
\end{center}
\title{Nonparametric Bayesian Aggregation for Massive
Data}
\runtitle{Supplement to Bayesian Aggregation}
\end{frontmatter}

\setcounter{subsection}{0}
\renewcommand{\thesubsection}{S.\arabic{section}.\arabic{subsection}}
\setcounter{subsubsection}{0}
\renewcommand{\thesubsubsection}{S.\arabic{section}.\arabic{subsubsection}}
\setcounter{equation}{0}
\renewcommand{\theequation}{S.\arabic{equation}}
\setcounter{lemma}{0}
\renewcommand{\thelemma}{S.\arabic{lemma}}
\setcounter{proposition}{0}
\renewcommand{\theproposition}{S.\arabic{proposition}}

This supplementary document is structured as follows.
\begin{itemize}
\item Section \ref{suppl:sec:proofs:two:lemmas} contains the proofs of Lemmas
\ref{a:limit:distribution} and \ref{a:limit:distribution:weak}.
\item Section \ref{suppl:sec:two:sections:app} contains the proofs of the main results in Section \ref{sec:ci:functional} and \ref{sec:asymp:post:infer} that were not included in the main paper.
\item Section \ref{sec:proof:thm:refine} proves Proposition \ref{refined:contraction:rate}, i.e., a uniform contraction rate result. Preliminary results relevant to the proof of Proposition \ref{refined:contraction:rate} are provided in Section \ref{app:sec:prelim}.
\item Section \ref{suppl:sec:intial:rate} includes a result that characterizes the posterior tail moments of $\|f-f_0\|^a$ for any $a\ge0$.
\item Section \ref{sec:suppl:additional:simulation} includes additional simulation results supplementary to Section \ref{sec:simulations}.
\end{itemize}

\subsection{Proofs of Lemmas \ref{a:limit:distribution} and \ref{a:limit:distribution:weak}}
\label{suppl:sec:proofs:two:lemmas} 

\begin{proof}[Proof of Lemma \ref{a:limit:distribution}]
We only show the first limit distribution
since the proof of the second one is similar.

Let $\eta_\nu=\tau_\nu v_\nu$. Then $\eta_\nu$ is a sequence of \textit{iid}
standard normals.
Note that
\[
\|W_n\|_2^2=\sum_{\nu=1}^\infty\frac{\eta_\nu^2}{\tau_\nu^2+n(1+\lambda\rho_\nu)}.
\]
Let $U_n=(n\|W_n\|_2^2-\zeta_{1,n})/\sqrt{2\zeta_{2,n}}$,
then we have
\[
U_n=\frac{1}{\sqrt{2\zeta_{2,n}}}\sum_{\nu=1}^\infty\frac{n(\eta_\nu^2-1)}{\tau_\nu^2+n(1+\lambda\rho_\nu)}.
\]
By straightforward calculations and Taylor's expansion of
$\log(1-x)$, it can be shown that
the logarithm of the moment generating function of $U_n$ equals
\begin{equation}\label{a:limit:distribution:eqn1}
\log{E\{\exp(t U_n)\}}=t^2/2+O\left(t^3\zeta_{2,n}^{-3/2}\zeta_{3,n}\right).
\end{equation}

Without loss of generality, assume that $N=n^a$ for some $a\ge1$.
Then $\alpha_1:=\min\{1/(2m+\beta),a/(2m+\beta)\}=1/(2m+\beta)$.
It follows by \cite[Lemma S.1]{SC14} that
$\zeta_{2,n}\asymp n^{\alpha_1}$ and $\zeta_{3,n}\asymp n^{\alpha_1}$,
so the remainder term in (\ref{a:limit:distribution:eqn1})
is $O(n^{-\alpha_1/2})=o(1)$.
So $\lim_{n\rightarrow\infty}E\{\exp(t U_n)\}=\exp(t^2/2)$. Proof is completed.
\end{proof}

\begin{proof}[Proof of Lemma \ref{a:limit:distribution:weak}]
The proof follows by moment generating function approach and direct calculations.
\end{proof}

\subsection{Proofs in Sections \ref{sec:ci:functional} and \ref{sec:asymp:post:infer}}\label{suppl:sec:two:sections:app}
This section contains the proofs in Sections \ref{sec:ci:functional} and \ref{sec:asymp:post:infer}.

\subsubsection*{Proofs in Section \ref{sec:ci:functional}}

\begin{proof}[Proof of Theorem \ref{cp:ci:functional}]
Recall in the proof of Theorem \ref{cp:cr:weak}
we showed that Rate Conditions ($\textbf{R}'$)
are satisfied.

It is easy to see that
\begin{equation}\label{ci:functional:eqn1}
F(W_n)\overset{d}{=}N(0,\theta_{1,n}^2),\,\,\,\,
\textrm{and}\,\,\,\,F(W_N)\overset{d}{=}N(0,\theta_{1,N}^2).
\end{equation}

For $1\le j\le s$, define
$R_{j,n}^F(\alpha)=\{f\in S^m(\mathbb{I}): |F(f)-F(\breve{f}_{j,n})|\le r_{F,j,n}(\alpha)\}$.
It follows by Theorem \ref{uniform:bvm:thm} that
$\max_{1\le j\le s}|1-\alpha-P_{0j}(R_{j,n}^F(\alpha))|=o_{P_{f_0}}(1)$.
Since $s=o(N^{\frac{4m^2+2m\beta-12m+1}{8m(2m+\beta)}}(\log{N})^{-\frac{3}{2}})$, it can be examined that
$NL_N^2=o(1)$. Together with the condition
$h^{-r}\lesssim N\theta_{1,N}^2$ and the fact
$\theta_{k,N}\le\theta_{k,n}$, one can verify that
$h^{-r}\lesssim N\theta_{1,N}^2\le N\theta_{1,n}^2=o(L_N^{-2}\theta_{1,n}^2)$.
So we have by (\ref{boundedness:F}) and Theorem \ref{uniform:post:mean} that
\[
\max_{1\le j\le s}|F(\Delta_j)|=O_{P_{f_0}}(h^{-r/2}L_N)=o_{P_{f_0}}(\theta_{1,n}).
\]
Combined with (\ref{ci:functional:eqn1}) we get that
\begin{eqnarray*}
P_{0j}(R_{j,n}^F(\alpha))&=&P(|F(W_n)-F(\Delta_j)|\le r_{F,j,n}(\alpha)|\textbf{D}_j)\\
&=&\Phi\left(\frac{r_{F,j,n}(\alpha)+F(\Delta_j)}{\theta_{1,n}}\right)
+\Phi\left(\frac{r_{F,j,n}(\alpha)-F(\Delta_j)}{\theta_{1,n}}\right)-1\\
&=&2\Phi\left(\frac{r_{F,j,n}(\alpha)}{\theta_{1,n}}\right)-1+o_{P_{f_0}}(1),\,\,
\textrm{uniformly for $1\le j\le s$.}
\end{eqnarray*}
The above argument leads to 
$\Phi(r_{F,j,n}(\alpha)/\theta_{1,n})=1-\alpha/2+o_{P_{f_0}}(1)$ uniformly for $1\le j\le s$,
which further leads to the following
\begin{equation}\label{ci:functional:eqn2}
\max_{1\le j\le s}|r_{F,j,n}(\alpha)/\theta_{1,n}-z_{\alpha/2}|=o_{P_{f_0}}(1).
\end{equation}

Consider the decomposition (\ref{an:important:eqn:cr:strong})
with $T_1,T_2,T_3,T_4$ being defined therein.
It follows by (\ref{cp:cr:strong:eqn3.5}) and rate condition $Na_n^2=o(1)$ that
$N\|T_1\|^2=O_{P_{f_0}}(Na_n^2)=o_{P_{f_0}}(1)$.
Meanwhile, it follows by Condition (\textbf{S}$'$),
$N^{-1}\asymp h^{2m+\beta}$ and $\lambda=h^{2m}$ and direct examinations that
\begin{eqnarray*}
N\|T_2\|^2&=&N\sum_{\nu=1}^\infty(a_{N,\nu}-1)^2|f_\nu^0|^2(1+\lambda\rho_\nu)\\
&\asymp&N\sum_{\nu=1}^\infty\left(\frac{\nu^{2m+\beta}}{\nu^{2m+\beta}+N(1+\lambda\nu^{2m})}\right)^2
|f_\nu^0|^2(1+\lambda\nu^{2m})\\
&\asymp&\sum_{\nu=1}^\infty\frac{(h\nu)^{2m+\beta}+(h\nu)^{4m+\beta}}{(1+(h\nu)^{2m}+(h\nu)^{2m+\beta})^2}
\times|f_\nu^0|^2\nu^{2m+\beta}=o(1),
\end{eqnarray*}
and
\begin{eqnarray*}
N\|T_4\|^2&=&N\sum_{\nu=1}^\infty a_{N,\nu}^2\left(\frac{\lambda\rho_\nu}{1+\lambda\rho_\nu}\right)^2
|f_\nu^0|^2(1+\lambda\rho_\nu)\\
&\asymp&\sum_{\nu=1}^\infty\frac{(h\nu)^{2m-\beta}}{1+(h\nu)^{2m}}\times |f_\nu^0|^2\nu^{2m+\beta}=o(1).
\end{eqnarray*}
By (\ref{breve:vs:tilde}) and $Ns^2L_N^2=o(1)$ we get
$\|\breve{f}_{N,\lambda}-\widetilde{f}_{N,\lambda}\|=o_{P_{f_0}}(N^{-1/2})$.
Therefore, $\|\breve{f}_{N,\lambda}-f_0-T_3\|\le\|\breve{f}_{N,\lambda}-\widetilde{f}_{N,\lambda}\|
+\|T_1+T_2+T_4\|=o_{P_{f_0}}(N^{-1/2})$.
If follows from (\ref{boundedness:F}) that
$|F(\breve{f}_{N,\lambda}-f_0)-F(T_3)|=o_{P_{f_0}}(h^{-r/2}N^{-1/2})$.

Note that $F(T_3)=\frac{1}{N}\sum_{i=1}^N\epsilon_i F(R_{X_i})$,
where the kernel $R$ is defined in the proof of Theorem \ref{cp:cr:strong}.
We will derive asymptotic distribution for $F(T_3)$.
Let $s_N^2=Var_{f_0}(\sum_{i=1}^N\epsilon_i F(R_{X_i}))$.
It is easy to show that
\[
s_N^2=N^3\sum_{\nu=1}^\infty\frac{F(\varphi_\nu)^2}{(\tau_\nu^2+N(1+\lambda\rho_\nu))^2}=N^3\theta_{2,N}^2.
\]
Clearly, by uniform boundedness of $\varphi_\nu$ and $F(\varphi_\nu)$, we get
\[
|F(R_x)|=|\sum_{\nu=1}^\infty a_{N,\nu}\frac{\varphi_\nu(x)F(\varphi_\nu)}{1+\lambda\rho_\nu}|\lesssim h^{-1},
\]
where the ``$\lesssim$" is free of $x\in\mathbb{I}$, and
\begin{equation}\label{ci:functional:eqn3}
E_{f_0}\{\epsilon^2 F(R_X)^2\}=N^2\sum_{\nu=1}^\infty\frac{F(\varphi_\nu)^2}{(\tau_\nu^2+N(1+\lambda\rho_\nu))^2}
=N^2\theta_{2,N}^2.
\end{equation}
Then for any $\delta>0$, by condition $E_{f_0}\{\epsilon^4|X\}\le M_4$ a.s.,
\begin{eqnarray*}
&&\frac{1}{s_N^2}\sum_{i=1}^N E_{f_0}\{\epsilon_i^2 F(R_{X_i})^2 I(|\epsilon_i F(R_{X_i})|\ge \delta s_N)\}\\
&\le&\frac{N}{s_N^2}(\delta s_N)^{-2} E_{f_0}\{\epsilon^4 F(R_X)^4\}\\
&\lesssim&\frac{N}{s_N^2}(\delta s_N)^{-2}h^{-2}E_{f_0}\{\epsilon^2 F(R_X)^2\}
\lesssim\delta^{-2}N^{-1}h^{-2+r}=o(1),
\end{eqnarray*}
where the last $o(1)$-term follows by $h\asymp h^\ast$ and $2-r<2m+\beta$.
By Lindeberg's central limit theorem, as $N\rightarrow\infty$,
\begin{equation}\label{ci:functional:eqn4}
\frac{F(T_3)}{\sqrt{N}\theta_{2,N}}=\frac{1}{s_N}\sum_{i=1}^N\epsilon_i F(R_{X_i})\overset{d}{\rightarrow}N(0,1).
\end{equation}

By condition $N^2\theta_{2,N}^2\gtrsim h^{-r}$,
we have
\[
\bigg|\frac{F(\breve{f}_{N,\lambda}-f_0-T_3)}{\sqrt{N}\theta_{2,N}}\bigg|=
o_{P_{f_0}}\left(\frac{h^{-r/2}N^{-1/2}}{\sqrt{N}\theta_{2,N}}\right)
=o_{P_{f_0}}(1).
\]
It follows by (\ref{ci:functional:eqn2}) that
\begin{equation}\label{ci:functional:eqn5}
r_{F,N}(\alpha)=\theta_{1,N}\sqrt{\frac{1}{s}\sum_{j=1}^s r_{F,j,n}(\alpha)^2/\theta_{1,n}^2}
=\theta_{1,N}z_{\alpha/2}(1+o_{P_{f_0}}(1)),
\end{equation}
leading to that
\[
\frac{r_{F,N}(\alpha)}{\sqrt{N}\theta_{2,N}}=\frac{\theta_{1,N}}{\sqrt{N}\theta_{2,N}}
\times z_{\alpha/2}(1+o_{P_{f_0}}(1)).
\]
It can be shown that
\[
\frac{\theta_{1,N}^2}{N\theta_{2,N}^2}=\frac{\sum_{\nu=1}^\infty
\frac{F(\varphi_\nu)^2}{1+\lambda\rho_\nu+N^{-1}\tau_\nu^2}}
{\sum_{\nu=1}^\infty\frac{F(\varphi_\nu)^2}{(1+\lambda\rho_\nu+N^{-1}\tau_\nu^2)^2}}\ge1,
\]
together with (\ref{ci:functional:eqn4}) we get that
\begin{eqnarray}\label{ci:functional:eqn6}
&&P_{f_0}(|F(f_0)-F(\breve{f}_{N,\lambda})|\le r_{F,N}(\alpha))\nonumber\\
&=&P_{f_0}\left(\bigg|\frac{F(\breve{f}_{N,\lambda}-f_0-T_3)}{\sqrt{N}\theta_{2,N}}
+\frac{F(T_3)}{\sqrt{N}\theta_{2,N}}\bigg|\le\frac{r_{F,N}(\alpha)}{\sqrt{N}\theta_{2,N}}\right)\nonumber\\
&\ge&P_{f_0}\left(\bigg|\frac{F(\breve{f}_{N,\lambda}-f_0-T_3)}{\sqrt{N}\theta_{2,N}}
+\frac{F(T_3)}{\sqrt{N}\theta_{2,N}}\bigg|\le z_{\alpha/2}(1+o_{P_{f_0}}(1))\right)\nonumber\\
&\rightarrow& 1-\alpha.
\end{eqnarray}
Notice that when $0<\sum_{\nu=1}^\infty F(\varphi_\nu)^2<\infty$,
$\frac{\theta_{1,N}^2}{N\theta_{2,N}^2}\rightarrow1$, leading to that
the probability in (\ref{ci:functional:eqn6}) approaches exactly $1-\alpha$.

In the end, we show that $P(R_N^F(\alpha)|\textbf{D})=1-\alpha+o_{P_{f_0}}(1)$,
where $R_N^F(\alpha)=\{f\in S^m(\mathbb{I}): |F(f)-F(\breve{f}_{N,\lambda})|\le r_{F,N}(\alpha)\}$.
By rate condition $N(a_N^2+a_n^2)=o(1)$,
proof of (\ref{cp:cr:strong:eqn5}) leading to
$\|\widetilde{f}^{or}_{N,\lambda}-\widetilde{f}_{N,\lambda}\|=O_{P_{f_0}}(a_N+a_n)$,
and (\ref{boundedness:F}) we have
\begin{eqnarray*}
\frac{F(\widetilde{f}^{or}_{N,\lambda}-\widetilde{f}_{N,\lambda})}{\theta_{1,N}}=O_{P_{f_0}}
\left(\frac{h^{-r/2}(a_N+a_n)}{\theta_{1,N}}\right)=o_{P_{f_0}}(1),
\end{eqnarray*}
where the last $o(1)$-term follows by condition $N\theta_{1,N}^2\gtrsim h^{-r}$ and Rate Condition ($\textbf{R}'$).
From (\ref{ci:functional:eqn5}) we get that
\begin{eqnarray}\label{ci:functional:eqn7}
P_0(R_N^F(\alpha))&=&P(W^{or}\in R_N^F(\alpha)|\textbf{D})\nonumber\\
&=&P(|F(W^{or})-F(\breve{f}_{N,\lambda})|\le r_{F,N}(\alpha)|\textbf{D})\nonumber\\
&=&P\left(\bigg|\frac{F(\widetilde{f}^{or}_{N,\lambda}-\widetilde{f}_{N,\lambda})}{\theta_{1,N}}
+\frac{F(W_N)}{\theta_{1,N}}\bigg|\le\frac{r_{F,N}(\alpha)}{\theta_{1,N}}\bigg|\textbf{D}\right)
\nonumber\\
&=&1-\alpha+o_{P_{f_0}}(1).
\end{eqnarray}
So it follows from (\ref{full:bvm:expression}) that
$P(R_N^F(\alpha)|\textbf{D})=1-\alpha+o_{P_{f_0}}(1)$.
Proof is completed.
\end{proof}

\subsubsection*{Proofs in Section \ref{sec:asymp:post:infer}}

\begin{proof}[Proof of Theorem \ref{cp:cr:strong:asymp}]
It follows from (\ref{cp:cr:strong:eqn4}) that $r_N(\alpha)-r_N^\dag(\alpha)=o_{P_{f_0}}(N^{-1}h^{-1/2})$, which together with (\ref{cp:cr:strong:eqn6}) leads to that
$\lim_{n\rightarrow\infty}P_{f_0}(f_0\in R_N^\dag(\alpha))=1$.
It follow from Lemma \ref{a:limit:distribution}, (\ref{full:bvm:expression}) and the proof of Theorem \ref{cp:cr:strong} that $P(R_N^\dag(\alpha)|\textbf{D})=1-\alpha+o_{P_{f_0}}(1)$.

It follows from (\ref{cp:cr:weak:eqn2}) that
$r_{\omega,N}(\alpha)^2-r_{\omega,N}^\dag(\alpha)^2=o_{P_{f_0}}(N^{-1})$.
Then the desired results on $R^{\dag\omega}_N(\alpha)$ directly follow from the proof of Theorem \ref{cp:cr:weak}.

It follows by (\ref{ci:functional:eqn5}) that $r_{F,N}^\dag(\alpha)=r_{F,N}(\alpha)(1+o_{P_{f_0}}(1))$.
Then the desired results on $CI_N^{\dag F}(\alpha)$ follow from (\ref{ci:functional:eqn6})
and (\ref{ci:functional:eqn7}).
\end{proof}

\subsection{Proofs of Proposition \ref{refined:contraction:rate} and relevant results}\label{sec:proof:thm:refine}

The goal of this section is to prove Proposition \ref{refined:contraction:rate}
and relevant results.
Before proofs, we exactly describe the Fr\'{e}chet derivatives of the likelihood function
that will be technically useful.
Suppose that $(Y,X)$ follows model (\ref{basic:model}) based on $f$.
Let $g, g_k\in S^m(\mathbb{I})$ for $k=1,2$.
For $j=1,2,\ldots,s$,
the Fr\'{e}chet derivative of $\ell_{jn}$ can be identified as
\[
D\ell_{jn}(g)g_1=\frac{1}{n}\sum_{i\in I_j}(Y_i-g(X_i))
\langle K_{X_i},g_1\rangle-\langle
\mathcal{P}_\lambda g,g_1\rangle:=\langle S_{j,n}(g),g_1 \rangle.
\]
Define $S_\lambda(g)=E\{S_{j,n}(g)\}$.
We also use $DS_\lambda$ and $D^2S_\lambda$
to represent the second- and third-order
Fr\'{e}chet derivatives of $S_\lambda$.
Note that $S_{j,n}(\widehat{f}_{j,n})=0$, and
$S_{j,n}(f)$ can be expressed as
\begin{eqnarray}\label{score}
S_{j,n}(f)=\frac{1}{n} \sum_{i\in I_j}(Y_i-f(X_i))
K_{X_i}-\mathcal{P}_\lambda f.
\end{eqnarray}
The Fr\'{e}chet derivative of $S_{j,n}$ is
denoted $DS_{j,n}(g)g_1g_2$.
These derivatives can be explicitly written as
\[
D^2 \ell_{jn}(g)
g_1g_2:= DS_{j,n}(g)g_1g_2
=-\frac{1}{n}\sum_{i\in I_j}g_1(X_i)g_2(X_i)-\langle \mathcal{P}_\lambda g_1,g_2\rangle,
\]

The proof of Proposition \ref{refined:contraction:rate} requires a series of preliminary lemmas.
Define
$H^m(b)=\{f\in S^m(\mathbb{I}): J(f)\le b^2\}$.
We first state a basic lemma about a concentration phenomenon of 
smoothing spline estimates in the distributed setup. 

\begin{lemma}\label{basic:thm}
	If $b,r,h,M$ are positives satisfying the following Rate Condition (\textbf{H}):
	\begin{enumerate}[(i)]
		\item\label{rate:basic:thm:ii} $h^{1/2}r\le 1$,
		\item\label{rate:basic:thm:iii} $c_K^2M^{1/2}rh^{-1/2}B(h)\le 1/2$, where $B(h)=A(h,2)$
		with $A(h,\varepsilon)$ given in (\ref{function:A}),
	\end{enumerate}
	then, for any $1\le j\le s$, the following two results hold:
	\begin{enumerate}[(a)]
		\item\label{basic:thm:i}
		$\sup_{f\in H^m(b)}P_f\left(\|\widehat{f}_{j,n}-f\|\ge \delta_n\right)\le
		2\exp(-Mnhr^2)$,
		where $\delta_n=bh^m+2c_K(C_\epsilon+M)r$ with $C_\epsilon=E\{(|\epsilon|+1)^2\exp(|\epsilon|+1)\}$
		an absolute constant;
		\item\label{basic:thm:ii}
		$\sup_{f\in H^m(b)}P_f\left(\|\widehat{f}_{j,n}-f-S_{j,n}(f)\|>a_n\right)
		\le 2\exp(-Mnhr^2),
		$
		where
		$a_n=c_K^2 M^{1/2}h^{-1/2}rB(h)\delta_n$.
		Here, $S_{j,n}(f)$ is the Fr\'{e}chet derivative of the likelihood function $\ell_{jn}(f)$; 
		see (\ref{score}) for its exact expression.
	\end{enumerate}
\end{lemma}

\begin{lemma}\label{pre:lemma:1} For any fixed constants $M>1$ and $b>0$, let
\begin{equation}\label{b:r:deltan}
r=(nh/\log{2s})^{-1/2},
\delta_n=bh^m+2c_K(C_\epsilon+M)r,
\end{equation}
\begin{equation}\label{an:bn}
a_n=c_K^2 M^{1/2}h^{-1/2}rB(h)\delta_n.
\end{equation}
Then as $n\rightarrow\infty$,
\[
P_{f_0}\left(\max_{1\le j\le s}\|\widehat{f}_{j,n}-f_0\|\ge\delta_n\right)\le
6sN^{-M}\rightarrow0,
\]
and
\[
P_{f_0}\left(\max_{1\le j\le s}\|\widehat{f}_{j,n}-f_0-S_{j,n}(f_0)\|>a_n\right)
\le 8sN^{-M}\rightarrow0.
\]
\end{lemma} 

\begin{proof}[Proof of Lemma \ref{pre:lemma:1}]
The result is a straightforward consequence of Lemma \ref{basic:thm}.
\end{proof}

\begin{lemma}\label{pre:lemma:2} It holds that
\begin{equation}\label{refined:rate:eqn1}
\max_{1\le j\le s}\|\widehat{f}_{j,n}-f_0-S_{j,n}(f_0)\|=O_{P_{f_0}}(a_n).
\end{equation}
\end{lemma}

\begin{proof}[Proof of Lemma \ref{pre:lemma:2}]
The proof follows by Lemma \ref{pre:lemma:1}, and simple fact
that $B(h)\lesssim h^{-\frac{2m-1}{4m}}$.
\end{proof}

\begin{lemma}\label{pre:lemma:3}
Under Condition (\textbf{S}), we get
$\max_{1\le j\le s}\|\widehat{f}_{j,n}-f_0\|=O_{P_{f_0}}(\widetilde{r}_n)$.
\end{lemma}

\begin{proof}[Proof of Lemma \ref{pre:lemma:3}]
Recall that
\[
S_{j,n}(f_0)=-\frac{1}{n} \sum_{i\in I_j}(Y_i-f_0(X_i))
K_{X_i}-\mathcal{P}_\lambda f_0.
\]
It was shown by \cite{SC13} that
$\mathcal{P}_\lambda\varphi_\nu=\frac{\lambda\varphi_\nu}{1+\lambda\varphi_\nu}\varphi_\nu$.
Since $f_0$ satisfies Condition (\textbf{S}),
\begin{eqnarray*}
\|\mathcal{P}_\lambda f_0\|^2&=&\langle\sum_{\nu=1}^\infty
f_\nu^0\frac{\lambda\rho_\nu}{1+\lambda\rho_\nu}\varphi_\nu,
\sum_{\nu=1}^\infty f_\nu^0\frac{\lambda\rho_\nu}{1+\lambda\rho_\nu}\varphi_\nu\rangle\\
&=&\sum_{\nu=1}^\infty|f_\nu^0|^2\frac{\lambda^2\rho_\nu^2}{1+\lambda\rho_\nu}\\
&=&\lambda^{1+\frac{\beta-1}{2m}}\sum_{\nu=1}^\infty|f_\nu^0|^2\rho_\nu^{1+\frac{\beta-1}{2m}}
\frac{(\lambda\rho_\nu)^{1-\frac{\beta-1}{2m}}}{1+\lambda\rho_\nu}=O(h^{2m+\beta-1}),
\end{eqnarray*}
where the last equation follows by $\lambda=h^{2m}$,
$\sup_{x\ge0}\frac{x^{1-\frac{\beta-1}{2m}}}{1+x}<\infty$, and
Condition (\textbf{S}).
On the other side, it follows by
the proof of (\ref{contraction:mapping:eqn0}) that
\begin{eqnarray*}
&&P_{f_0}\left(\max_{1\le j\le s}\|\sum_{i\in I_j}(Y_i-f_0(X_i))
K_{X_i}\|\ge L(M)n(nh/\log{2s})^{-1/2}\right)\\
&\le& 2s\exp\left(-Mnh(nh/\log{2s})^{-1}\right)=(2s)^{1-M}\rightarrow0,\,\,\textrm{as $M\rightarrow\infty$,}
\end{eqnarray*}
where $L(M):= c_K(C_\epsilon+M)$.
This implies that
\[
\max_{1\le j\le s}\|\sum_{i\in I_j}(Y_i-f_0(X_i))
K_{X_i}\|=O_{P_{f_0}}(n(nh/\log{2s})^{-1/2}),
\]
and hence,
\[
\max_{1\le j\le s}\|S_{j,n}(f_0)\|=
O_{P_{f_0}}((nh/\log{2s})^{-1/2}+h^{m+\frac{\beta-1}{2}})=O_{P_{f_0}}(\widetilde{r}_n).
\]
Together with (\ref{refined:rate:eqn1}) of Lemma \ref{pre:lemma:2} and the rate condition $a_n\lesssim \widetilde{r}_n$,
we get that
$\max_{1\le j\le s}\|\widehat{f}_{j,n}-f_0\|=O_{P_{f_0}}(\widetilde{r}_n)$.
\end{proof}
Consider a function class
\begin{equation}\label{space:G:C}
\mathcal{G}=\{g\in S^m(\mathbb{I}): \|g\|_\infty\le1,
J(g,g)\le c_K^{-2} h^{-2m+1}\}.
\end{equation}

\begin{lemma}\label{pre:lemma:4}
For any fixed constant $M>1$, as $n\rightarrow\infty$,
\[
P_{f_0}\left(\max_{1\le j\le s}\sup_{g\in\mathcal{G}}\|Z_{j,n}(g)\|\le
B(h)\sqrt{M\log{N}}\right)\rightarrow1,
\]
where
$Z_{j,n}(g)=\frac{1}{\sqrt{n}}\sum_{i\in I_j}[\psi_{j,n}(Z_i;g)K_{X_i}
-E\{\psi_{j,n}(Z_i;g)K_{X_i}\}]$,
$\psi_{j,n}(Z_i;g)=c_K^{-1}h^{1/2}g(X_i)$.
\end{lemma}
\begin{proof}[Proof of Lemma \ref{pre:lemma:4}]
It is easy to see that $\psi_{j,n}(Z_i;g)$ satisfies the Lipschitz continuity condition (\ref{Lip:cont:psi}).
Then the result directly follows by Lemma \ref{basic:lemma:UCT} (see appendix).
\end{proof}

\begin{lemma}\label{pre:lemma:5}
For $j=1,\ldots,s$,
\begin{enumerate}[(1).]
\item $\ell_{jn}(f)-\ell_{jn}(\widehat{f}_{j,n})=I_{j,n}(f)$, 
where $I_{j,n}(f)=\int_0^1\int_0^1sDS_{j,n}(\widehat{f}_{j,n}+ss'(f-\widehat{f}_{j,n}))
(f-\widehat{f}_{j,n})(f-\widehat{f}_{j,n})dsds'$ for any $f\in S^m(\mathbb{I})$;
\item $I_{j,n}(f)=T_{j}(f)-\frac{1}{2}\|f-\widehat{f}_{j,n}\|^2$,
where recall that (see \ref{Tj1:Tj2})
\begin{eqnarray}
T_{j}(f)&=&-\frac{1}{2n}\sum_{i\in I_j}
[(f-\widehat{f}_{j,n})(X_i)^2-E_X\{(f-\widehat{f}_{j,n})(X)^2\}].
\end{eqnarray}
\end{enumerate}
\end{lemma}
\begin{proof}[Proof of Lemma \ref{pre:lemma:5}]
Let $\Delta f=f-\widehat{f}_{j,n}$.
Therefore,
\begin{eqnarray*}
I_{j,n}(f)&=&-\frac{1}{n}\int_0^1\int_0^1s\sum_{i\in I_j}(\Delta f)(X_i)^2 dsds'
-\lambda J(\Delta f,\Delta f)/2\nonumber\\
&=&-\frac{1}{2n}\sum_{i\in I_j}[(\Delta f)(X_i)^2-E_X\{(\Delta f)(X)^2\}]
-\frac{1}{2}\|\Delta f\|^2\\
&=& T_{j}(f)-\frac{1}{2}\|\Delta f\|^2.
\end{eqnarray*}

By Taylor's expansion in terms of Fr\'{e}chet derivatives,
$\ell_{jn}(f)-\ell_{jn}(\widehat{f}_{j,n})=S_{j,n}(\widehat{f}_{j,n})(f-\widehat{f}_{j,n})
+I_{j,n}(f)=I_{j,n}(f)$.
\end{proof}

\begin{lemma}\label{pre:lemma:6}
There exists a universal constant $c_3>0$ s.t.
\[
\Pi(\|f-f_0\|\le\widetilde{r}_n)
\ge\exp(-c_3\widetilde{r}_n^{-\frac{2}{2m+\beta-1}}),
\]
where recall that $\Pi$ is the probability measure induced by $G$.
\end{lemma}
\begin{proof}[Proof of Lemma \ref{pre:lemma:6}]
Note that $\lambda\le\widetilde{r}_n^{\frac{4m}{2m+\beta-1}}$.
Then it follows by Lemma \ref{concentration:lemma} (with $d_n$ therein replaced by $\widetilde{r}_n$)
and the proof of Theorem \ref{initial:contraction:rate} that
\begin{eqnarray*}
\Pi(\|f-f_0\|\le\widetilde{r}_n)&=&P(\|G-f_0\|\le\widetilde{r}_n)\\
&\ge&P(V(G-f_0)\le \widetilde{r}_n^2/2,\lambda J(G-f_0)\le\widetilde{r}_n^2/2)\\
&\ge&P(V(G-f_0)\le \widetilde{r}_n^2/2,J(G-f_0)\le\widetilde{r}_n^{\frac{2(\beta-1)}{2m+\beta-1}}/2)\\
&=&P(\widetilde{V}(\widetilde{G}-\widetilde{f}_0)
\le\widetilde{r}_n^2/2,\widetilde{J}(\widetilde{G}-\widetilde{f}_0)\le\widetilde{r}_n^{\frac{2(\beta-1)}{2m+\beta-1}}/2)\\
&\ge&P(\widetilde{V}(\widetilde{G}-\omega)\le (1/\sqrt{2}-1/2)^2\widetilde{r}_n^2,
\widetilde{J}(\widetilde{G}-\omega)\le (1/\sqrt{2}-1/2)^2\widetilde{r}_n^{\frac{2(\beta-1)}{2m+\beta-1}})\\
&\ge&\exp(-\|\omega\|_\beta^2/2)\\
&&\times P(\widetilde{V}(\widetilde{G})\le (1/\sqrt{2}-1/2)^2\widetilde{r}_n^2,
\widetilde{J}(\widetilde{G})\le (1/\sqrt{2}-1/2)^2\widetilde{r}_n^{\frac{2(\beta-1)}{2m+\beta-1}})\\
&\ge&\exp(-\|\omega\|_\beta^2/2)
P(\widetilde{V}(\widetilde{G})\le (1/\sqrt{2}-1/2)^2\widetilde{r}_n^2/2)\\
&&\times P(\widetilde{J}(\widetilde{G})\le (1/\sqrt{2}-1/2)^2\widetilde{r}_n^{\frac{2(\beta-1)}{2m+\beta-1}}/2)\\
&\ge&\exp(-c_3\widetilde{r}_n^{-\frac{2}{2m+\beta-1}}),
\end{eqnarray*}
where $c_3>0$ is a universal constant.
\end{proof}

\begin{proof}[Proof of Proposition \ref{refined:contraction:rate}]
Fix any $\varepsilon\in(0,1)$.
Let $M_1$ be a large constant so that (thanks to Lemma \ref{pre:lemma:3}) the event
\begin{equation}\label{event:En'}
\mathcal{E}_n'=\{\max_{1\le j\le s}\|\widehat{f}_{j,n}-f_0\|\le M_1\widetilde{r}_n\}
\end{equation}
has probability approaching one. Meanwhile, for a fixed constant $M>1$, define
\begin{equation}\label{event:En''}
\mathcal{E}''_n=\left\{\max_{1\le j\le s}\sup_{g\in\mathcal{G}}\|Z_{j,n}(g)\|\le
B(h)\sqrt{M\log{N}}\right\}.
\end{equation}
By Lemma \ref{pre:lemma:4} we have that $\mathcal{E}_n''$
has $P_{f_0}$-probability approaching one.
Thus, it holds that, when $n$ becomes large,
$P_{f_0}(\mathcal{E}_n)\ge 1-\varepsilon/2$, where $\mathcal{E}_n:=\mathcal{E}_n'\cap\mathcal{E}_n''$.
In the rest of the proof we simply assume that $\mathcal{E}_n$ holds.

For some positive constant $M_0$,
it follows by Theorem \ref{initial:contraction:rate} that
$$\max_{1\le j\le s}E\{\|f-f_0\|^aI(\|f-f_0\|\ge M_0 r_n)|\textbf{D}_j\}=O_{P_{f_0}}(s^2\exp(-nr_n^2)).$$ 
Let $C'>M_1$ be a constant to be further determined later,
then we have that
\begin{eqnarray*}
&&\max_{1\le j\le s}E\{\|f-f_0\|^aI(\|f-f_0\|\ge 2C'\widetilde{r}_n)|\textbf{D}_j\}\\
&\le&\max_{1\le j\le s}E\{\|f-f_0\|^aI(\|f-f_0\|\ge M_0r_n)|\textbf{D}_j\}\\
&&+\max_{1\le j\le s}E\{\|f-f_0\|^aI(2C'\widetilde{r}_n\le\|f-f_0\|\le M_0 r_n)|\textbf{D}_j\}.
\end{eqnarray*}
The first term is $O_{P_{f_0}}(s^2\exp(-nr_n^2))$.
Thus, when $n$ is sufficiently large,
\[
P_{f_0}\left(\max_{1\le j\le s}E\{\|f-f_0\|^aI(\|f-f_0\|\ge M_0r_n)|\textbf{D}_j\}
\ge M's^2\exp(-nr_n^2)/2\right)\le\varepsilon/2
\]
for a large constant $M'>0$.

Next we only need to handle the second term.
Let $\Delta f=f-\widehat{f}_{j,n}$. It follows by Lemma \ref{pre:lemma:5} that
$I_{j,n}(f)=T_{j}(f)-\frac{1}{2}\|\Delta f\|^2$,
and
$\ell_{jn}(f)-\ell_{jn}(\widehat{f}_{j,n})=I_{j,n}(f)$.
Therefore,
\begin{eqnarray*}
&&E\{\|f-f_0\|^aI(f\in A_n)|\textbf{D}_j\}\\
&=&\frac{\int_{A_n}\|f-f_0\|^a\exp(n(\ell_{jn}(f)-\ell_{jn}(\widehat{f}_{j,n})))d\Pi(f)}
{\int_{S^m(\mathbb{I})}\exp(n(\ell_{jn}(f)-\ell_{jn}(\widehat{f}_{j,n})))d\Pi(f)}
=\frac{\int_{A_n}\|f-f_0\|^a\exp(nI_{j,n}(f))d\Pi(f)}{\int_{S^m(\mathbb{I})}\exp(nI_{j,n}(f))d\Pi(f)},
\end{eqnarray*}
where $A_n=\{f\in S^m(\mathbb{I}): 2C'\widetilde{r}_n\le\|f-f_0\|\le M_0 r_n\}$.

Let
\[
J_{j1}=\int_{S^m(\mathbb{I})}\exp(nI_{j,n}(f))d\Pi(f),\,\,
J_{j2}=\int_{A_n}\|f-f_0\|^a\exp(nI_{j,n}(f))d\Pi(f).
\]
Then on $\mathcal{E}_n$ and for $\|f-f_0\|\le\widetilde{r}_n$,
we have $\|f-\widehat{f}_{j,n}\|\le \|f-f_0\|+\|\widehat{f}_{j,n}-f_0\|\le (M_1+1)\widetilde{r}_n$.

Let $d_n=c_K(M_1+1)h^{-1/2}\widetilde{r}_n$.
It follows by similar arguments as above (\ref{contraction:mapping:eqn1}) that
$d_n^{-1}\Delta f\in\mathcal{G}$.
Note that on $\mathcal{E}_n$ and for $\|f-f_0\|\le\widetilde{r}_n$,
for all $1\le j\le s$,
\begin{eqnarray}\label{eqn:Tj2}
|T_{j}(f)|&=&\frac{1}{2n}\bigg|
\sum_{i\in I_j}[(\Delta f)(X_i)^2-E_X\{(\Delta f)(X)^2\}]\bigg|\nonumber\\
&=&
\frac{1}{2n}\bigg|\langle\sum_{i\in I_j}[(\Delta f)(X_i)K_{X_i}
-E_X\{(\Delta f)(X)K_X\}],\Delta f\rangle\bigg|\nonumber\\
&\le&\frac{1}{2n}\|\Delta f\|\times\|\sum_{i\in I_j}[(\Delta f)(X_i)K_{X_i}
-E_X\{(\Delta f)(X)K_X\}]\|\nonumber\\
&=&\frac{c_K h^{-1/2} d_n\|\Delta f\|}{2\sqrt{n}}\times\|Z_{j,n}(d_n^{-1}\Delta f)\|\nonumber\\
&\le&\frac{c_K h^{-1/2} d_n\|\Delta f\|}{2\sqrt{n}} B(h)\sqrt{M\log{N}}\nonumber\\
&\le&D(c_K,M,M_1)\times n^{-1/2}h^{-\frac{6m-1}{4m}}\widetilde{r}_n^2\sqrt{\log{N}}
\le D(c_K,M,M_1)\times \widetilde{r}_n^2b_n,
\end{eqnarray}
where $D(c_K,M,M_1)$ is constant depending only on $c_K,M_1,M$.

It follows that on $\mathcal{E}_n$ and for all $1\le j\le s$,
\begin{eqnarray*}
J_{j1}&\ge&\int_{\|f-f_0\|\le\widetilde{r}_n}\exp(nI_{j,n}(f))d\Pi(f)\\
&=&\int_{\|f-f_0\|\le\widetilde{r}_n}\exp\left(nT_{j}(f)-\frac{n}{2}\|f-\widehat{f}_{j,n}\|^2\right)d\Pi(f)\\
&\ge&\exp\left(-[D(c_K,M,M_1)b_n+(M_1+1)^2/2]n\widetilde{r}_n^2\right)
\Pi(\|f-f_0\|\le\widetilde{r}_n).
\end{eqnarray*}

Since $\Pi(\|f-f_0\|\le\widetilde{r}_n)
\ge\exp(-c_3\widetilde{r}_n^{-\frac{2}{2m+\beta-1}})$ (Lemma \ref{pre:lemma:6}),
together with
$\widetilde{r}_n\ge (nh)^{-1/2}+h^{m+\frac{\beta-1}{2}}\ge 2n^{-\frac{2m+\beta-1}{2(2m+\beta)}}$,
we get that
$n\widetilde{r}_n^{2+\frac{2}{2m+\beta-1}}\ge
n(4n^{-\frac{2m+\beta-1}{2m+\beta}})^{1+\frac{1}{2m+\beta-1}}=4$.
Therefore, $\widetilde{r}_n^{-\frac{2}{2m+\beta-1}}\le n\widetilde{r}_n^2/4$,
leading to
\begin{equation}\label{refined:rate:small:ball}
\Pi(\|f-f_0\|\le\widetilde{r}_n)\ge\exp\left(-\frac{c_3}{4}n\widetilde{r}_n^2\right).
\end{equation}
This implies by rate conditions $b_n\le1$ that, on $\mathcal{E}_n$ and for any $1\le j\le s$,
\begin{eqnarray*}
J_{j1}&\ge&\exp\left(-[D(c_K,M,M_1)b_n+(M_1+1)^2/2+c_3/4]
n\widetilde{r}_n^2\right)\\
&\ge&\exp\left(-[D(c_K,M,M_1)+(M_1+1)^2/2+c_3/4]
n\widetilde{r}_n^2\right).
\end{eqnarray*}

Next we handle $J_{j2}$. The idea is similar to how we handle $J_{j1}$
but with technical difference. Let $\Delta f=f-\widehat{f}_{j,n}$.
Note that $\widetilde{r}_n^2\le r_n^2\log(2s)$, and hence, on $\mathcal{E}_n$,
for any $f\in A_n$, i.e.,
$\|f-f_0\|\le M_0r_n$,
we get that
$\|\Delta f\|=\|\widehat{f}_{j,n}-f\|\le\|\widehat{f}_{j,n}-f_0\|+\|f-f_0\|\le
M_1\widetilde{r}_n+M_0r_n\le(M_0+M_1)r_n\sqrt{\log(2s)}$.
Let $d_{*n}=c_K(M_0+M_1)h^{-1/2}r_n\sqrt{\log(2s)}$.
Then $d_{*n}^{-1}\Delta f\in\mathcal{G}$.
Using previous similar arguments handling $T_{j}(f)$, we have that on $\mathcal{E}_n$,
for any $f\in A_n$ and $1\le j\le s$,
\begin{eqnarray*}
|T_{j}(f)|
&\le&\frac{\|\Delta f\|}{2\sqrt{n}}c_Kh^{-1/2}d_{*n}\cdot B(h)\sqrt{M\log{N}}\\
&\le&\frac{1}{2}c_K^2(M_0+M_1)^2M^{1/2}n^{-1/2}h^{-1}r_n^2B(h)(\log{N})^{3/2}\\
&\le&D(c_K,M,M_0,M_1)\times n^{-1/2}r_n^2h^{-\frac{6m-1}{4m}}(\log{N})^{3/2}\\
&=&D(c_K,M,M_0,M_1)\times r_n^2 b_n\le D(c_K,M,M_0,M_1)\times\widetilde{r}_n^2,
\end{eqnarray*}
where $D(c_K,M,M_0,M_1)$ is constant only depending on $c_K,M,M_0,M_1$
and the last inequality follows by rate condition
$r_n^2 b_n\le\widetilde{r}_n^2$.
It is easy to see that on $\mathcal{E}_n$ and for any $f\in A_n$ and $1\le j\le s$,
$\|\widehat{f}_{j,n}-f\|\ge \|f-f_0\|-\|\widehat{f}_{j,n}-f_0\|\ge (2C'-M_1)\widetilde{r}_n$,
leading to that
\begin{eqnarray*}
J_{j2}&\le&\exp\left(-\left(\frac{(2C'-M_1)^2}{2}-D(c_K,M,M_0,M_1)\right)
n\widetilde{r}_n^2\right)C(a,\Pi),
\end{eqnarray*}
where $C(a,\Pi)=\int_{S^m(\mathbb{I})}\|f-f_0\|^ad\Pi(f)$ is the $a$th prior moment of $\|f-f_0\|$
which is finite. Choose $C'>M_1$ to be large such that
\begin{eqnarray*}
\frac{(2C'-M_1)^2}{2}&\ge&1+D(c_K,M,M_1)
+D(c_K,M,M_0,M_1)+(M_1+1)^2/2+c_3/4.
\end{eqnarray*}
Therefore, on $\mathcal{E}_n$,
\[
\max_{1\le j\le s}E\{\|f-f_0\|^aI(f\in A_n)|\textbf{D}_j\}
\le\frac{\max_{1\le j\le s}J_{j2}}{\min_{1\le j\le s}J_{j1}}
\le\exp(-n\widetilde{r}_n^2)C(a,\Pi).
\]
So we get that
\begin{eqnarray*}
P_{f_0}\left(\max_{1\le j\le s}
E\{\|f-f_0\|^aI(f\in A_n)|\textbf{D}_j\}\ge \exp(-n\widetilde{r}_n^2)
C(a,\Pi)\right)
\le P_{f_0}(\mathcal{E}_n^c)\le\varepsilon/2.
\end{eqnarray*}
By $\widetilde{r}_n^2\le r_n^2\log(2s)$, the above leads to that
\begin{eqnarray*}
&&P_{f_0}\left(\max_{1\le j\le s}
E\{\|f-f_0\|^aI(\|f-f_0\|\ge 2C'\widetilde{r}_n)|\textbf{D}_j\}\right.\\
&&\left.\ge (M'+C(a,\Pi))s^2\exp(-n\widetilde{r}_n^2/\log(2s))\right)\le\varepsilon.
\end{eqnarray*}
Proof is completed.
\end{proof}

\subsection{Proofs of other results in Section \ref{sec:proof:thm:refine}}\label{app:sec:prelim}

Let $N(\varepsilon,\mathcal{G},\|\cdot\|_\infty)$
be the $\varepsilon$-packing number in terms of supremum norm,
where recall that the space $\mathcal{G}$
is defined in (\ref{space:G:C}).
The following result can be found in \cite{VG00}.
\begin{lemma}\label{entropy}
There exists a universal constant $c_0>0$ s.t.
for any $\varepsilon>0$,
\[
\log{N(\varepsilon,\mathcal{G},\|\cdot\|_\infty)}\le c_0
(\sqrt{2}c_K^{-1})^{1/m}h^{-\frac{2m-1}{2m}}\varepsilon^{-1/m}.
\]
\end{lemma}
For $r\ge0$, define $\Psi(r)=\int_0^r\sqrt{\log(1+\exp(x^{-1/m}))}dx$.
For arbitrary $\varepsilon>0$, define
\begin{eqnarray}\label{function:A}
A(h,\varepsilon)
&=&\frac{32\sqrt{6}}{\tau}
\sqrt{2}c_K^{-1}c_0^m h^{-(2m-1)/2}\Psi\left(\frac{1}{2\sqrt{2}}c_Kc_0^{-m}h^{(2m-1)/2}\varepsilon\right)\nonumber\\
&&+\frac{10\sqrt{24}\varepsilon}{\tau}
\sqrt{\log\left(1+\exp\left(2c_0((\sqrt{2})^{-1}c_Kh^{(2m-1)/2}\varepsilon)^{-1/m}\right)\right)},
\end{eqnarray}
where $\tau=\sqrt{\log{1.5}}\approx 0.6368$. 

We have the following useful lemma.
\begin{lemma}\label{basic:lemma:UCT}
For any $1\le j\le s$ and $f\in S^m(\mathbb{I})$,
suppose that $\psi_{j,n,f}(z;g)$ is a measurable function defined upon
$z=(y,x)\in\mathcal{Y}\times\mathbb{I}$ and $g\in\mathcal{G}$
satisfying $\psi_{j,n,f}(z;0)=0$ and the following Lipschitz continuity condition:
for any $i\in I_j$ and $g_1,g_2\in\mathcal{G}$,
\begin{equation}\label{Lip:cont:psi}
|\psi_{j,n,f}(Z_i;g_1)-\psi_{j,n,f}(Z_i;g_2)|\le
c_K^{-1} h^{1/2} \|g_1-g_2\|_\infty.
\end{equation}
Then for any constant $t\ge0$ and $n\ge1$,
\[
\sup_{f\in S^m(\mathbb{I})}
P_f\left(\sup_{g\in\mathcal{G}}\|Z_{j,n,f}(g)\|_f>t\right)
\le 2\exp\left(-\frac{t^2}{B(h)^2}\right),
\]
where $B(h)=A(h,2)$ and
\[
Z_{j,n,f}(g)=\frac{1}{\sqrt{n}}\sum_{i\in I_j}[\psi_{j,n,f}(Z_i;g)K_{X_i}
-E_f\{\psi_{j,n,f}(Z_i;g)K_{X_i}\}].
\]
\end{lemma}
\begin{proof}[Proof of Lemma \ref{basic:lemma:UCT}]
For any $f\in S^m(\mathbb{I})$
and $n\ge1$, and any $g_1,g_2\in\mathcal{G}$,
we get that
\begin{eqnarray*}
\|(\psi_{j,n,f}(Z_i;g_1)-\psi_{j,n,f}(Z_i;g_2)) K_{X_i}\|\le c_K^{-1}h^{1/2}\|g_1-g_2\|_\infty c_Kh^{-1/2}=\|g_1-g_2\|_\infty.
\end{eqnarray*}
By Theorem 3.5 of \cite{P94}, for any $t>0$,
$P_f\left(\|Z_{j,n,f}(g_1)-Z_{j,n,f}(g_2)\|\ge t\right)\le
2\exp\left(-\frac{t^2}{8\|g_1-g_2\|_\infty^2}\right)$.
Then by Lemma 8.1 in \cite{K08}, we have
\[
\left\|\|Z_{j,n,f}(g_1)-Z_{j,n,f}(g_2)\|\right\|_{\psi_2}\le
\sqrt{24}\|g_1-g_2\|_\infty,
\]
where $\|\cdot\|_{\psi_2}$ denotes the
Orlicz norm associated with $\psi_2(s):=\exp(s^2)-1$.
Recall $\tau=\sqrt{\log{1.5}}\approx 0.6368$.
Define $\phi(x)=\psi_2(\tau x)$.
Then it can be shown by elementary calculus that
$\phi(1)\le 1/2$, and for any $x,y\ge1$,
$\phi(x)\phi(y)\le\phi(xy)$.
By a careful examination of the proof of Lemma 8.2,
it can be shown that for any random variables $\xi_1,\ldots,\xi_l$,
\begin{equation}\label{lemma8.2:MK}
\|\max_{1\le i\le l}\xi_i\|_{\psi_2}\le \frac{2}{\tau}\psi_2^{-1}(l)\max_{1\le i\le l}
\|\xi_i\|_{\psi_2}.
\end{equation}

Next we use a ``chaining" argument. Let
$T_0\subset T_1\subset T_2\subset\cdots\subset T_\infty:= \mathcal{G}$
be a sequence of finite nested sets satisfying the following properties:
\begin{itemize}
\item for any $T_q$ and any $s,t\in T_q$, $\|s-t\|_\infty\ge \varepsilon 2^{-q}$;
each $T_q$ is ``maximal" in the sense that if one adds any point in $T_q$,
then the inequality will fail;
\item the cardinality of $T_q$ is upper bounded by
\begin{eqnarray*}
\log{|T_q|}\le \log{N(\varepsilon 2^{-q},\mathcal{G},\|\cdot\|_\infty)}\le c_0(\sqrt{2}c_K^{-1})^{1/m}h^{-(2m-1)/(2m)}(\varepsilon 2^{-q})^{-1/m},
\end{eqnarray*}
where $c_0>0$ is absolute constant;
\item each element $t_{q+1}\in T_{q+1}$ is uniquely linked to an element $t_q\in T_q$
which satisfies $\|t_q-t_{q+1}\|_\infty\le\varepsilon 2^{-q}$.
\end{itemize}

For arbitrary $s_{k+1},t_{k+1}\in T_{k+1}$
with $\|s_{k+1}-t_{k+1}\|_\infty\le\varepsilon$,
choose two chains (both being of length $k+2$) $t_q$ and $s_q$
with $t_q,s_q\in T_q$ for $0\le q\le k+1$.
The ending points $s_0$ and $t_0$
satisfy
\begin{eqnarray*}
\|s_0-t_0\|_\infty&\le&
\sum_{q=0}^k[\|s_q-s_{q+1}\|_\infty+\|t_q-t_{q+1}\|_\infty]
+\|s_{k+1}-t_{k+1}\|_\infty\\
&\le& 2\sum_{q=0}^k\varepsilon 2^{-q}+\varepsilon\le 5\varepsilon,
\end{eqnarray*}
and hence, $\left\|\|Z_{j,n,f}(s_0)-Z_{j,n,f}(t_0)\|_f\right\|_{\psi_2}
\le 5\sqrt{24}\varepsilon$.
It follows by the proof of Theorem
8.4 of \cite{K08} and (\ref{lemma8.2:MK}) that
\begin{eqnarray*}
&&\left\|\max_{s_{k+1},t_{k+1}\in T_{k+1}}
\|Z_{j,n,f}(s_{k+1})-Z_{j,n,f}(t_{k+1})-(Z_{j,n,f}(s_0)-Z_{j,n,f}(t_0))\|\right\|_{\psi_2}\\
&\le&2\sum_{q=0}^k\left\|\max_{\substack{u\in T_{q+1},
v\in T_q\\ \textrm{$u,v$ link each other}}}
\|Z_{j,n,f}(u)-Z_{j,n,f}(v)\|\right\|_{\psi_2}\\
&\le&\frac{4}{\tau}\sum_{q=0}^k\psi_2^{-1}(N(2^{-q-1}\varepsilon,\mathcal{G},\|\cdot\|_\infty))\\
&&\times\max_{\substack{u\in T_{q+1},
v\in T_q\\ \textrm{$u,v$ link each other}}}
\left\|\|Z_{j,n,f}(u)-Z_{j,n,f}(v)\|\right\|_{\psi_2}\\
&\le&\frac{4\sqrt{24}}{\tau}\sum_{q=0}^k
\sqrt{\log\left(1+N(\varepsilon 2^{-q-1},\mathcal{G},\|\cdot\|_\infty)\right)}\varepsilon 2^{-q}\\
&\le&
\frac{8\sqrt{24}}{\tau}\sum_{q=1}^{k+1}
\sqrt{\log\left(1+\exp\left(c_0c_K^{-1/m}h^{-(2m-1)/(2m)}(\varepsilon 2^{-q})^{-1/m}\right)\right)}
\varepsilon 2^{-q}\\
&\le&
\frac{32\sqrt{6}}{\tau}\int_0^{\varepsilon/2}
\sqrt{\log\left(1+\exp\left(c_0c_K^{-1/m}h^{-(2m-1)/(2m)}x^{-1/m}\right)\right)}dx\\
&=&\frac{32\sqrt{6}}{\tau}
c_K^{-1}c_0^m h^{-(2m-1)/2}\Psi\left(\frac{1}{2}c_Kc_0^{-m}h^{(2m-1)/2}\varepsilon\right).
\end{eqnarray*}
On the other hand,
\begin{eqnarray*}
\left\|\max_{\substack{u,v\in T_0
\\ \|u-v\|_\infty\le 5\varepsilon}}\|Z_{j,n,f}(u)-Z_{j,n,f}(v)\|_f\right\|_{\psi_2}&\le&
\frac{2}{\tau}\psi_2(|T_0|^2)\max_{\substack{u,v\in T_0
\\ \|u-v\|_\infty\le 5\varepsilon}}
\left\|\|Z_{j,n,f}(u)-Z_{j,n,f}(v)\|_f\right\|_{\psi_2}\\
&\le&\frac{2}{\tau}\psi_2^{-1}(N(\varepsilon,\mathcal{G},\|\cdot\|_\infty)^2)
(5\sqrt{24}\varepsilon).
\end{eqnarray*}
Therefore,
\begin{eqnarray*}
\left\|\max_{\substack{
s,t\in T_{k+1}\\
\|s-t\|_\infty\le\varepsilon}}\|Z_{j,n,f}(s)-Z_{j,n,f}(t)\|\right\|_{\psi_2}&\le&\frac{32\sqrt{6}}{\tau}
c_K^{-1}c_0^m h^{-(2m-1)/2}\Psi\left(\frac{1}{2}c_Kc_0^{-m}h^{(2m-1)/2}\varepsilon\right)\\
&&+\frac{2}{\tau}\psi_2^{-1}(N(\varepsilon,\mathcal{G},\|\cdot\|_\infty)^2)
(5\sqrt{24}\varepsilon)\\
&\le&
\frac{32\sqrt{6}}{\tau}
c_K^{-1}c_0^m h^{-(2m-1)/2}\Psi\left(\frac{1}{2}c_Kc_0^{-m}h^{(2m-1)/2}\varepsilon\right)\\
&&+\frac{10\sqrt{24}\varepsilon}{\tau}
\sqrt{\log\left(1+\exp\left(2c_0(c_Kh^{(2m-1)/2}\varepsilon)^{-1/m}\right)\right)}\\
&=&A(h,\varepsilon).
\end{eqnarray*}

Now for any $g_1,g_2\in\mathcal{G}$ with $\|g_1-g_2\|_\infty\le\varepsilon/2$.
Let $k\ge2$, hence, $2^{1-k}\le 1-\|g_1-g_2\|_\infty/\varepsilon$.
Since $T_k$ is ``maximal", there exist $s_k,t_k\in T_k$ s.t.
$\max\{\|g_1-s_k\|_\infty,\|g_2-t_k\|_\infty\}\le \varepsilon 2^{-k}$.
It is easy to see that $\|s_k-t_k\|_\infty\le\varepsilon$.
So
\begin{eqnarray*}
\|Z_{j,n,f}(g_1)-Z_{j,n,f}(g_2)\|&\le& \|Z_{j,n,f}(g_1)-Z_{j,n,f}(s_k)\|+\|Z_{j,n,f}(g_2)-Z_{j,n,f}(t_k)\|\\
&&+\|Z_{j,n,f}(s_k)-Z_{j,n,f}(t_k)\|\\
&\le&4\sqrt{n}\varepsilon 2^{-k}+\max_{\substack{u,v\in T_k
\\ \|u-v\|_\infty\le\varepsilon}}\|Z_{j,n,f}(u)-Z_{j,n,f}(v)\|.
\end{eqnarray*}
Therefore, letting $k\rightarrow\infty$ we get that
\begin{eqnarray*}
&&\left\|\sup_{\substack{g_1,g_2\in\mathcal{G}
\\ \|g_1-g_2\|_\infty\le\varepsilon/2}}
\|Z_{j,n,f}(g_1)-Z_{j,n,f}(g_2)\|\right\|_{\psi_2}\\
&\le&4\sqrt{n}\varepsilon 2^{-k}/\sqrt{\log{2}}
+\left\|\max_{\substack{u,v\in T_k
\\ \|u-v\|_\infty\le\varepsilon}}\|Z_{j,n,f}(u)-Z_{j,n,f}(v)\|\right\|_{\psi_2}\\
&\le&4\sqrt{n}\varepsilon 2^{-k}/\sqrt{\log{2}}+A(h,\varepsilon)\rightarrow A(h,\varepsilon).
\end{eqnarray*}
Taking $\varepsilon=2$ in the above inequality, we get that
\begin{eqnarray*}
\left\|\sup_{\substack{g_1,g_2\in\mathcal{G}\\
\|g_1-g_2\|_\infty\le1}}\|Z_{j,n,f}(g_1)-Z_{j,n,f}(g_2)\|\right\|_{\psi_2}
\le A(h,2)= B(h).
\end{eqnarray*}
By Lemma 8.1 in \cite{K08}, we have
\[
P_f\left(\sup_{\substack{g\in\mathcal{G}}}\|Z_{j,n,f}(g)\|\ge t\right)\le
2\exp\left(-\frac{t^2}{B(h)^2}\right).
\]
Note that the right hand side in the above
does not depend on $f$. This completes the proof.
\end{proof}

\begin{proof}[Proof of Lemma \ref{basic:thm}]
Let $f\in H^m(b)$ be the parameter based on which the data are drawn.
It is easy to see that
$DS_\lambda(f)g=-E\{g(X)K_X\}-\mathcal{P}_\lambda g,\,\,\,\,\forall g\in S^m(\mathbb{I})$.
Therefore,
for any $g,\widetilde{g}\in S^m(\mathbb{I})$,
$\langle DS_\lambda(f)g,\widetilde{g}\rangle=-\langle g,\widetilde{g}\rangle$,
implying $DS_\lambda(f)=-id$.

The proof of (\ref{basic:thm:i}) is finished in two parts.

\textbf{Part I}:
For any $f\in S^m(\mathbb{I})$,
define an operator mapping $S^m(\mathbb{I})$ to $S^m(\mathbb{I})$:
\[
T_{1f}(g)=g+S_\lambda(f+g),\,\,g\in S^m(\mathbb{I}).
\]
First observe that, under $P_f$ with $f\in H^m(b)$,
\[
\|S_\lambda(f)\|=\|\mathcal{P}_\lambda f\|
=\sup_{\|g\|=1}|\langle \mathcal{P}_\lambda f,g\rangle|\le \sqrt{\lambda J(f)}
\le h^mb.
\]
Let $r_{1n}=bh^m$. Let $\mathbb{B}(r_{1n})=\{g\in S^m(\mathbb{I}): \|g\|\le r_{1n}\}$
be the $r_{1n}$-ball. For any $g\in \mathbb{B}(r_{1n})$,
using $DS_\lambda(f)=-id$, it is easy to see that
$\|T_{1f}(g)\|=\|S_\lambda(f)\|\le bh^m=r_{1n}$.
Therefore, $T_{1f}$ maps $\mathbb{B}(r_{1n})$ to itself.
For any $g_1,g_2\in\mathbb{B}(r_{1n})$, by Taylor's expansion we have
\begin{eqnarray*}
\|T_{1f}(g_1)-T_{1f}(g_2)\|
&=&\|g_1-g_2+S_\lambda(f+g_1)-S_\lambda(f+g_2)\|\\
&=&\|g_1-g_2+\int_0^1DS_\lambda(f+g_2+sg)gds\|=0.
\end{eqnarray*}
This shows that
$T_{1f}$ is a contraction mapping which maps $\mathbb{B}(r_{1n})$ into $\mathbb{B}(r_{1n})$.
By contraction mapping theorem (see \cite{R76}),
$T_{1f}$ has a unique fixed point $g'\in\mathbb{B}(r_{1n})$
satisfying $T_{1f}(g')=g'$. Let $f_\lambda=f+g'$.
Then $S_\lambda(f_\lambda)=0$ and
$\|f_\lambda-f\|\le r_{1n}$.

\textbf{Part II}:
For any $f\in H^m(b)$, under (\ref{basic:model}) with $f$ being the truth,
let $f_\lambda$ be the function obtained in \textbf{Part I}
s.t. $\|f_\lambda-f\|\le r_{1n}$.
Define an operator
\[
T_{2f}(g)=g+S_{j,n}(f_\lambda+g),\,\,g\in S^m(\mathbb{I}).
\]
Rewrite $T_{2f}$ as
\begin{eqnarray*}
T_{2f}(g)&=&[DS_{j,n}(f_\lambda)g-DS_\lambda(f_\lambda)g]+S_{j,n}(f_\lambda).
\end{eqnarray*}
Denote the above two terms by $I_{1f},I_{2f}$, respectively.

For any $i\in I_j$,
let $R_i=(Y_i-f_\lambda(X_i))K_{X_i}-E_f\{(Y-f_\lambda(X))K_X\}$.
Obviously,
\begin{eqnarray*}
\|E_f\{(Y-f_\lambda(X))K_X\}\|&=&\sup_{\|g\|=1}|\langle E_f\{(Y-f_\lambda(X))K_X\},g\rangle|\\
&=&\sup_{\|g\|=1}|E_f\{(Y-f_\lambda(X))g(X)\}|
\le\|f-f_\lambda\|\le r_{1n}.
\end{eqnarray*}
Therefore,
$\|R_i\|
\le c_Kh^{-1/2}|Y_i-f_\lambda(X_i)|+r_{1n}$
which leads to that
\begin{eqnarray*}
E\left\{\exp\left(\frac{\|R_i\|}{c_Kh^{-1/2}}\right)\right\}\le
E\left(\exp(|\epsilon_i|+1)\right)\le C_\epsilon,
\end{eqnarray*}
where $C_\epsilon=E\{(|\epsilon|+1)^2\exp(|\epsilon|+1)\}$.
Let $\delta=hr/c_K$. By condition $rh^{1/2}\le1$, we have 
\begin{eqnarray*}
E\{\exp(\delta\|R_i\|)-1-\delta\|R_i\|\}&\le&E\{(\delta\|R_i\|)^2\exp(\delta\|R_i\|)\}
\le c_K^2C_\epsilon\delta^2h^{-1}.
\end{eqnarray*}
It follows by Theorem 3.2 of \cite{P94} that, for $L(M):=c_K(C_\epsilon+M)$,
\begin{eqnarray}\label{contraction:mapping:eqn0}
P_f\left(\|\sum_{i\in I_j} R_i\|_f\ge L(M)nr\right)
&\le& 2\exp\left(-L(M)\delta n r+c_K^2C_\epsilon nh^{-1}\delta^2\right)\nonumber\\
&=&2\exp(-Mnhr^2).
\end{eqnarray}
We note that the right hand side in (\ref{contraction:mapping:eqn0}) does not depend on $f$.
Moreover, it is easy to see that $S_{j,n}(f_\lambda)=S_{j,n}(f_\lambda)-S_\lambda(f_\lambda)
=\frac{1}{n}\sum_{i\in I_j}R_i$.
Let
\[
\mathcal{E}_{n,1}=\{\|S_{j,n}(f_\lambda)\|\le L(M)r\},
\]
then $\sup_{f\in H^m(C)}P_f(\mathcal{E}_{n,1}^c)\le2\exp(-Mnhr^2)$.
Define
$\psi_{j,n}(X_i;g)=c_K^{-1}h^{1/2}g(X_i),
\,\,i\in I_j$,
and $Z_{j,n}(g)=\frac{1}{\sqrt{n}}\sum_{i\in I_j}[\psi_{j,n}(X_i;g)K_{X_i}-
E_f\{\psi_{j,n}(X_i;g)K_{X_i}\}]$.
By Lemma \ref{basic:lemma:UCT},
$\sup_{f\in H^m(b)}P_f(\mathcal{E}_{n,2}^c)\le 2\exp(-Mnhr^2)$,
where $\mathcal{E}_{n,2}=\{\sup_{g\in\mathcal{G}}\|Z_{j,n}(g)\|\le \sqrt{Mnhr^2}B(h)\}$.

For any $g\in S^m(\mathbb{I})\backslash\{0\}$,
let $\bar{g}=g/d_n'$, where $d_n'=c_Kh^{-1/2}\|g\|$.
It follows that
\[
\|\bar{g}\|_\infty\le c_K h^{-1/2}\|\bar{g}\|=c_K h^{-1/2}\|g\|/d_n'=1,\,\,\textrm{and}
\]
\begin{eqnarray*}
J(\bar{g},\bar{g})&=&d_n'^{-2}J(g,g)
=h^{-2m}\frac{\lambda J(g,g)}{c_K^2 h^{-1}\|g\|^2}
\le h^{-2m}\frac{\|g\|^2}{c_K^2 h^{-1}\|g\|^2}
\le c_K^{-2}h^{-2m+1}.
\end{eqnarray*}
Therefore, $\bar{g}\in\mathcal{G}$.
Consequently, on $\mathcal{E}_{n,2}$,
for any $g\in S^m(\mathbb{I})\backslash\{0\}$, we get
$\|Z_{j,n}(\bar{g})\|\le\sqrt{Mnhr^2}B(h)$, which leads to that
\begin{eqnarray}\label{contraction:mapping:eqn1}
\|DS_{j,n}(f_\lambda)g-DS_\lambda(f_\lambda)g\|
&=&\frac{1}{n}\|\sum_{i\in I_j}[g(X_i)K_{X_i}
-E\{g(X_i)K_{X_i}\}]\|_f\nonumber\\
&\le& c_K^2M^{1/2}rh^{-1/2}B(h)\|g\|\le \|g\|/2,
\end{eqnarray}
where the last inequality follows by condition $c_K^2M^{1/2}rh^{-1/2}B(h)\le 1/2$.
Note that the above inequality also holds for $g=0$.

Let $r_{2n}=2L(M)r$.
Therefore, it follows by (\ref{contraction:mapping:eqn1}) that,
for any $f\in H^m(b)$, on $\mathcal{E}_n:=
\mathcal{E}_{n,1}\cap\mathcal{E}_{n,2}$
and for any $g\in\mathbb{B}(r_{2n})$,
$\|T_{2f}(g)\|\le \|g\|/2+r_{2n}/2\le r_{2n}$.
Meanwhile, for any $g_1,g_2\in\mathbb{B}(r_{2n})$, replacing
$g$ by $g_1-g_2$ in (\ref{contraction:mapping:eqn1}), we get that
$\|T_{2f}(g_1)-T_{2f}(g_2)\|\le \|g_1-g_2\|/2$.
Therefore, for any $f\in H^m(b)$, on $\mathcal{E}_n$,
$T_{2f}$ is a contraction mapping from $\mathbb{B}(r_{2n})$ to itself.
By contraction mapping theorem, there exists uniquely an element $g''\in\mathbb{B}(r_{2n})$
s.t. $T_{2f}(g'')=g''$.
Let $\widehat{f}_{j,n}=f_\lambda+g''$. Clearly, $S_{j,n}(\widehat{f}_{j,n})=0$,
and hence, $\widehat{f}_{j,n}$ is the maximizer of $\ell_{jn}$;
see (\ref{smoothing:spline:likhood}).
So we get that, on $\mathcal{E}_n$,
$\|\widehat{f}_{j,n}-f\|_f\le \|f_\lambda-f\|+\|\widehat{f}_{j,n}-f_\lambda\|
\le r_{1n}+r_{2n}=bh^m+2L(M)r$. The desired conclusion follows by the trivial fact:
$\sup_{f\in H^m(b)}P_f(\mathcal{E}_n^c)\le 4\exp(-Mnhr^2)$.
Proof of (\ref{basic:thm:i}) is completed.

Next we show (\ref{basic:thm:ii}).

For any $f\in H^m(b)$, let $\widehat{f}_{j,n}$
be the penalized MLE of $f$ obtained by (\ref{smoothing:spline:likhood}).
Let $g_n=\widehat{f}_{j,n}-f$,
$\delta_n=bh^m+2L(M)r$, $d_n'=c_Kh^{-1/2}\delta_n$. 

On $\mathcal{E}_n$,
we have $\|g_n\|_f\le \delta_n$. Let $\bar{g}=g_n/d_n'$. Clearly, $\bar{g}\in\mathcal{G}$.
Then we get that
\begin{eqnarray}\label{FBR:eqn1}
&&\|S_{j,n}(f+g_n)-S_{j,n}(f)-(S_\lambda(f+g_n)-S_\lambda(f))\|\nonumber\\
&=&\frac{1}{n}\|\sum_{i\in I_j}[
g_n(X_i)K_{X_i}
-E_X\{g_n(X)K_X\}]\|\nonumber\\
&=&\frac{c_K d_n'}{\sqrt{nh}}\|Z_{j,n}(\bar{g})\|
\le c_K^2 M^{1/2}h^{-1/2}rB(h)\delta_n=a_n.
\end{eqnarray}

Since $S_{j,n}(f+g_n)=0$ and $DS_\lambda(f)=-id$,
from (\ref{FBR:eqn1}) we have on $\mathcal{E}_n$,
\begin{eqnarray*}
a_n&\ge&\|S_{j,n}(f)+DS_\lambda(f)g_n+\int_0^1\int_0^1sD^2S_\lambda(f+ss'g_n)g_ng_ndsds'\|
=\|S_{j,n}(f)-g_n\|
\end{eqnarray*}
which implies that
$\|\widehat{f}_{j,n}-f-S_{n,\lambda}(f)\|\le a_n$.
Since $\sup_{f\in H^m(bC)}P_f(\mathcal{E}_n^c)\le 4\exp(-Mnhr^2)$,
proof of (\ref{basic:thm:ii}) is completed.
\end{proof}

\subsection{An initial contraction rate}\label{suppl:sec:intial:rate}
Theorem \ref{initial:contraction:rate} below states that the $s$ posterior measures uniformly
contract at rate $r_n=(nh)^{-1/2}+h^m$, where recall 
that $h=\lambda^{1/(2m)}$. This is an initial rate result that holds irrespective the diverging rate of $s$.

\begin{Theorem}\label{initial:contraction:rate} (An Initial Contraction Rate)
Suppose $f_0=\sum_{\nu=1}^\infty f_\nu^0\varphi_\nu$ satisfies
Condition (\textbf{S}).
Let $a\ge0$ be a fixed constant.
If $r_n=o(h^{3/2})$, $h^{1/2}\log{N}=o(1)$, $nh^{2m+1}\ge1$, then there exists a universal constant $M>0$ s.t.
$$\max_{1\le j\le s}E\{\|f-f_0\|^aI(\|f-f_0\|\ge Mr_n)|\textbf{D}_j\}=O_{P_{f_0}}(s^2\exp(-nr_n^2))$$ as $n\to\infty$, no matter $s$ is fixed or diverges at any rate.
\end{Theorem}

Before proving Theorem \ref{initial:contraction:rate},
we present a preliminary lemma.

Let $\{\widetilde{\varphi}_\nu:\nu\ge1\}$
be a bounded orthonormal basis of $L^2(\mathbb{I})$
under usual $L^2$ inner product.
For any $b\in[0,\beta]$, define
\[
\widetilde{H}_b=\{\sum_{\nu=1}^\infty f_\nu\widetilde{\varphi}_\nu:\sum_{\nu=1}^\infty
f_\nu^2\rho_\nu^{1+b/(2m)}<\infty\}.
\]
Then $\widetilde{H}_b$ can be viewed as a version of Sobolev space with regularity $m+b/2$.
Define $\widetilde{G}=\sum_{\nu=1}^\infty v_\nu\widetilde{\varphi}_\nu$, a centered GP,
and $\widetilde{f}_0=\sum_{\nu=1}^\infty f_\nu^0\widetilde{\varphi}_\nu$.
Define $\widetilde{V}(f,g)=\langle f,g\rangle_{L^2}=\int_0^1 f(x)g(x)dx$,
the usual $L^2$ inner product,
$\widetilde{J}(f)=\sum_{\nu=1}^\infty |\widetilde{V}(f,\widetilde{\varphi}_\nu)|^2\rho_\nu$,
a functional on $\widetilde{H}_0$.
For simplicity, denote $\widetilde{V}(f)=\widetilde{V}(f,f)$.
Clearly, $\widetilde{f}_0\in\widetilde{H}_\beta$.
Since $\widetilde{G}$ is a Gaussian process with covariance function
\[
\widetilde{r}(s,t)=E\{\widetilde{G}(s)\widetilde{G}(t)\}
=\sum_{\nu=1}^m\sigma_\nu^2\widetilde{\varphi}_\nu(s)\widetilde{\varphi}_\nu(t)
+\sum_{\nu>m}\rho_\nu^{-(1+\frac{\beta}{2m})}
\widetilde{\varphi}_\nu(s)\widetilde{\varphi}_\nu(t),
\]
it follows by \cite{VZ08} that $\widetilde{H}_\beta$ is the RKHS of $\widetilde{G}$.
For any $\widetilde{H}_b$ with $0\le b\le \beta$, define inner product
\[
\langle \sum_{\nu=1}^\infty f_\nu\widetilde{\varphi}_\nu,\sum_{\nu=1}^\infty g_\nu\widetilde{\varphi}_\nu\rangle_b
=\sum_{\nu=1}^m \sigma_\nu^{-2}f_\nu g_\nu+\sum_{\nu>m}f_\nu g_\nu\rho_\nu^{1+\frac{b}{2m}}.
\]

Let $\|\cdot\|_b$ be the norm corresponding to the above inner product.
The following lemma is used in the proof of Theorem \ref{initial:contraction:rate}.
Its proof can be found in \cite{SC14}.
\begin{lemma}\label{concentration:lemma}
Let $d_n$ be any positive sequence.
If Condition (\textbf{S}) holds,
then there exists $\omega\in \widetilde{H}_\beta$ such that
\begin{enumerate}[(i).]
\item\label{concentration:lemma:i} $\widetilde{V}(\omega-\widetilde{f}_0)\le \frac{1}{4}d_n^2$,
\item\label{concentration:lemma:ii} $\widetilde{J}(\omega-\widetilde{f}_0)\le\frac{1}{4}
d_n^{\frac{2(\beta-1)}{2m+\beta-1}}$,
\item\label{concentration:lemma:iii} $\|\omega\|_{\beta}^2=O(d_n^{-\frac{2}{2m+\beta-1}})$.
\end{enumerate}
\end{lemma}

To ease reading, we sketch the proof of Theorem \ref{initial:contraction:rate}. 
We first show the following result:
for any $\varepsilon>0$, as $n\rightarrow\infty$,
\begin{equation}\label{posterior:consistency}
\max_{1\le j\le s}\int_{\|f-f_0\|_\infty\ge \varepsilon}\|f-f_0\|^a dP(f|\textbf{D}_j)
=O_{P_{f_0}}(s^2\exp(-nr_n^2))
\end{equation}
To show (\ref{posterior:consistency}),
we can rewrite the posterior density of $f$ by
\[
p(f|\textbf{D}_j)=\frac{\prod_{i\in I_j} (p_f/p_{f_0})(Z_i)\exp(-n\lambda J(f)/2)d\Pi(f)}
{\int_{S^m(\mathbb{I})}\prod_{i\in I_j} (p_f/p_{f_0})(Z_i)\exp(-n\lambda J(f)/2)d\Pi(f)},
\,\,1\le j\le s,
\]
where recall that $p_f(z)$ is the probability density of $Z=(Y,X)$ under $f$.
For $1\le j\le s$, define
\begin{equation}\label{term:Ij1}
I_{j1}=\int_{S^m(\mathbb{I})}\prod_{i\in I_j}(p_f/p_{f_0})(Z_i)\exp(-n\lambda J(f)/2)d\Pi(f),
\end{equation}
\begin{equation}\label{term:Ij2}
I_{j2}=\int_{A_n}\|f-f_0\|^a\prod_{i\in I_j} (p_f/p_{f_0})(Z_i)\exp(-\frac{n\lambda}{2}J(f))d\Pi(f),
\end{equation}
\begin{equation}\label{term:Ij2'}
I_{j2}'=\int_{A_n'}\|f-f_0\|^a\prod_{i\in I_j}
(p_f/p_{f_0})(Z_i)\exp(-\frac{n\lambda}{2}J(f))d\Pi(f),
\end{equation}
where $A_n=\{f\in S^m(\mathbb{I}):\|f-f_0\|\ge 2\delta_n\}$ and
$A_n'=\{f\in S^m(\mathbb{I}):\|f-f_0\|\ge \sqrt{2}Mr_n\}$,
with the quantities $\delta_n,M$ specified later.
Using LeCam's uniformly consistent test \cite{GGV00},
we will show that $\max_{1\le j\le s}I_{j2}/I_{j1}$ is of an
exponential order (in the sense of $P_{f_0}$).
Then (\ref{posterior:consistency}) holds by taking $a=0$ in $I_{j2}$.
The proof of Theorem \ref{initial:contraction:rate}
will be completed by decomposing $I_{j2}'/I_{j1}$ into three terms
based on an auxiliary event $\{f\in S^m(\mathbb{I}): \|f-f_0\|_{\infty}\le\varepsilon\}$
with each term of an exponential order. 

\begin{proof}[Proof of Theorem \ref{initial:contraction:rate}]
Note that there exists a universal constant $c'>0$ such that
$\Psi(x)\le c' x^{1-1/(2m)}$ for any $0<x<1$.
Therefore, there exists a universal constant $c''>0$ s.t.
$B(h)\le c'' h^{-(2m-1)/(4m)}$.

Define $B_n=\{f\in S^m(\mathbb{I}): V(f-f_0)\le r_n^2, J(f-f_0)\le
r_n^{\frac{2(\beta-1)}
{2m+\beta-1}}\}$.
Then
\begin{eqnarray*}
I_{j1}&\ge& \int_{B_n}\prod_{i\in I_j}(p_f/p_{f_0})(Z_i)\exp(-n\lambda J(f)/2)d\Pi(f)\\
&=&\int_{B_n}\exp(\sum_{i\in I_j}R_i(f,f_0))\exp(-n\lambda J(f)/2)d\Pi(f),
\end{eqnarray*}
where $R_i(f,f_0)=\log\left(p_f(Z_i)/p_{f_0}(Z_i)\right)=Y_i(f(X_i)-f_0(X_i))-f(X_i)^2/2+f_0(X_i)^2/2$ for any $i\in I_j$.
Define $d\Pi^\ast(f)=d\Pi(f)/\Pi(B_n)$,
a reduced probability measure on $B_n$.
By Jensen's inequality,
\begin{eqnarray*}
&&\log\int_{B_n}\exp(\sum_{i\in I_j} R_i(f,f_0))\exp(-n\lambda J(f)/2)d\Pi^\ast(f)\\
&\ge& \int_{B_n}\left(\sum_{i\in I_j} R_i(f,f_0)-n\lambda J(f)/2\right)d\Pi^\ast(f)\\
&=&\int_{B_n}\sum_{i\in I_j} [R_i(f,f_0)-E_{f_0}\{R_i(f,f_0)\}]d\Pi^\ast(f)\\
&&+n\int_{B_n}E_{f_0}\{R_i(f,f_0)\}d\Pi^\ast(f)
-\int_{B_n}\frac{n\lambda J(f)}{2}d\Pi^\ast(f)\\
&:=& J_{j1}+J_{j2}+J_{j3}.
\end{eqnarray*}
For any $f\in B_n$,
$\|f-f_0\|^2=V(f-f_0)+\lambda J(f-f_0)
\le r_n^2+\lambda
r_n^{\frac{2(\beta-1)}
{2m+\beta-1}}$.
By \cite[Lemma A.9]{SC18} and the condition $h^{-3/2}r_n=o(1)$,
we can choose $n$
to be sufficiently large so that
$\|f-f_0\|_\infty\le ch^{-1/2}\|f-f_0\|\le
c\sqrt{h^{-1}r_n^2+h^{2m-1}}\le 1$.

It follows by Taylor's expansion and $E_{f_0}\{Y_i-f_0(X_i)|X_i\}=0$,
that for any $f\in B_n$,
\[
|E_{f_0}\{R_i(f,f_0)\}|=E_{f_0}\{(f(X)-f_0(X))^2\}/2
\le r_n^2/2.
\]
Therefore, $J_{j2}\ge -nr_n^2/2$ for any $1\le j\le s$.

Since $r_n^2=o(1)$, we can choose $n$ to be
large so that $|E_{f_0}\{R_i(f,f_0)\}|\le1$.
Meanwhile, for any $f\in B_n$, for some $s\in[0,1]$, we have
\begin{eqnarray*}
|R_i(f,f_0)|&=&|Y_i(f(X_i)-f_0(X_i))-f(X_i)^2/2+f_0(X_i)^2/2|\\
&=&|Y_i-f_0(X_i)
-\frac{1}{2}(f-f_0)(X_i)|\times |(f-f_0)(X_i)|\\
&\le&|Y_i-f_0(X_i)|+1/2=|\epsilon_i|+1/2.
\end{eqnarray*}
We have used $\|f-f_0\|_\infty\le1$ in the above inequalities.

For any $1\le i\le N$, define $A_i=\{|\epsilon_i|\le 2\log{N}\}$.
It is easy to check that
$P_{f_0}(\cup_{i=1}^N A_i^c)\rightarrow0$, as $N\rightarrow\infty$.
Define $\xi_i=\int_{B_n}R_i(f,f_0)d\Pi^\ast(f)\times I_{A_i}$,
we get that $|\xi_i|\le 2\log{N}+1/2$, a.s. It can also be shown by $r_n^2\gg 1/n\ge 1/N$ that,
as $n,N\to\infty$,
\begin{eqnarray*}
|E_{f_0}\{\int_{B_n}R_i(f,f_0)d\Pi^\ast(f)\times I_{A_i^c}\}|&\le&E_{f_0}\{(|\epsilon_i|+1/2)\times I_{A_i^c}\}\\
&\le& C_\epsilon (1/N+1/N^2)\le r_n^2,
\end{eqnarray*}
where $C_\epsilon$ is an absolute constant.

Let $\delta=1/(\sqrt{n}r_n)$. Note that by the condition $h^{1/2}\log{N}=o(1)$ we have
$\delta\log{N}=(\log{N})/(\sqrt{n}r_n)\le h^{1/2}\log{N}=o(1)$,
we can let $n$ be large so that $\delta (2\log{N}+1)\le1$.
Let $d_i=\xi_i-E_{f_0}\{\xi_i\}$ for $i\in I_j$,
then it is easy to see that
\[
|d_i|\le |\xi_i|+|E_{f_0}\{\xi_i\}|\le 2\log{N}+1,\,\, a.s.
\]
Let $e_i=E_{f_0}\{\exp(\delta|d_i|)-1-\delta|d_i|\}$.
It can be shown using inequality $\exp(x)-1-x\le x^2\exp(x)$ for $x\ge0$
and Cauchy-Schwartz inequality that
\begin{eqnarray*}
|e_i|&\le&E_{f_0}\{\delta^2 d_i^2\exp(\delta|d_i|)\}\\
&\le&e\delta^2E_{f_0}\{d_i^2\}\\
&\le&e\delta^2E_{f_0}\{\xi_i^2\}\\
&\le&e\delta^2\int_{B_n}E_{f_0}\{R_i(f,f_0)^2\}d\Pi^\ast(f)\\
&\le&e\delta^2\int_{B_n}E_{f_0}\{(|\epsilon_i|+1/2)^2(f-f_0)(X_i)^2\}d\Pi^\ast(f)\\
&\le& eC_\epsilon\delta^2r_n^2,
\end{eqnarray*}
where the last step follows from $V(f-f_0)\le r_n^2$ for any $f\in B_n$.
Therefore, it follows by \cite[Theorem 3.2]{P94} that
\begin{eqnarray}\label{initial:rate:eqn1}
&&P_{f_0}\left(\max_{1\le j\le s}|\sum_{i\in I_j}[\xi_i-E_{f_0}\{\xi_i\}]|\ge
4\sqrt{n}r_n\log{N}\right)\nonumber\\
&\le& sP_{f_0}\left(|\sum_{i\in I_j}[\xi_i-E_{f_0}\{\xi_i\}]|\ge
4\sqrt{n}r_n\log{N}\right)\nonumber\\
&\le&2s\exp(-4\sqrt{n}r_n(\log{N})\delta
+eC_\epsilon\delta^2nr_n^2)\nonumber\\
&\le& 2s/N^2\rightarrow0,\,\,\textrm{as $N\rightarrow\infty$.}
\end{eqnarray}
Since $\sqrt{n}r_n\gg \log{N}$, we can let $n$ be large so that
$4\sqrt{n}r_n\log{N}\le nr_n^2$.
Since on $\cap_{i=1}^N A_i$,
\[
J_{j1}=\sum_{i\in I_j}[\xi_i-E_{f_0}\{\xi_i\}]-nE_{f_0}\{\int_{B_n}R_i(f,f_0)d\Pi^\ast(f)
\times I_{A_i^c}\},
\]
we get from (\ref{initial:rate:eqn1}) that with $P_{f_0}$-probability approaching one,
for any $1\le j\le s$,
\[
J_{j1}\ge-4\sqrt{n}r_n\log{N}-nr_n^2\ge -2nr_n^2.
\]
Meanwhile, for any $f\in B_n$, $J(f)\le (1+J(f_0)^{1/2})^2$. Therefore,
$J_{j3}\ge -\frac{(1+J(f_0)^{1/2})^2}{2}n\lambda$.
So, with probability approaching one, for any $1\le j\le s$,
\[
I_{j1}\ge \exp\left(-5nr_n^2/2-\frac{(1+J(f_0)^{1/2})^2}{2}n\lambda\right)\Pi(B_n).
\]

To proceed, we need a lower bound for $\Pi(B_n)$.
It follows by Lemma \ref{concentration:lemma} by replacing $d_n$ therein by $r_n$,
by Gaussian correlation inequality (see Theorem 1.1 of \cite{L99}),
by Cameron-Martin theorem (see \cite{CM44} or \cite[eqn (4.18)]{KLL94})
and \cite[Example 4.5]{HJSD79} that
\begin{eqnarray}\label{initial:rate:eqn0}
\Pi(B_n)&=&P(V(G-f_0)\le r_n^2,J(G-f_0)\le
r_n^{\frac{2(\beta-1)}{2m+\beta-1}})\nonumber\\
&=&P(\widetilde{V}(\widetilde{G}-\widetilde{f}_0)\le r_n^2,
\widetilde{J}(\widetilde{G}-\widetilde{f}_0)\le
r_n^{\frac{2(\beta-1)}{2m+\beta-1}})\nonumber\\
&\ge& P(\widetilde{V}(\widetilde{G}-\omega)\le r_n^2/4,
\widetilde{J}(\widetilde{G}-\omega)\le
r_n^{\frac{2(\beta-1)}{2m+\beta-1}}/4)\nonumber\\
&\ge&\exp(-\frac{1}{2}\|\omega\|_\beta^2)P(\widetilde{V}(\widetilde{G})\le r_n^2/4,\widetilde{J}(\widetilde{G})\le
r_n^{\frac{2(\beta-1)}{2m+\beta-1}}/4)\nonumber\\
&\ge&\exp(-\frac{1}{2}\|\omega\|_\beta^2)P(\widetilde{V}(\widetilde{G})\le r_n^2/8)P(\widetilde{J}(\widetilde{G})\le
r_n^{\frac{2(\beta-1)}{2m+\beta-1}}/8)\nonumber\\
&\ge&\exp(-c_1 r_n^{-2/(2m+\beta-1)}),
\end{eqnarray}
where $c_1>0$ is a universal constant.

Since $\beta>1$ and $r_n^2=(nh)^{-1}+\lambda\ge n^{-2m/(2m+1)}$,
we get $r_n^2\ge\lambda$ and
$nr_n^{\frac{2(2m+\beta)}{2m+\beta-1}}\ge n^{1-\frac{2m(2m+\beta)}{(2m+1)(2m+\beta-1)}}>1$,
so $nr_n^2>r_n^{-\frac{2}{2m+\beta-1}}$.
Consequently, with $P_{f_0}$-probability approaching one
\begin{equation}\label{initial:rate:eqn2}
\min_{1\le j\le s}I_{j1}\ge \exp(-c_2nr_n^2),
\end{equation}
where $c_2=5/2+(1+J(f_0)^{1/2})^2/2+c_1$.

Let $b=2\sqrt{c_2+1}$ and $C\ge b^2/4$.
Next we examine
$I_{j2}$ defined in (\ref{term:Ij2})
with $A_n=\{f\in S^m(\mathbb{I}):\|f-f_0\|\ge 2\delta_n\}$, for
$\delta_n=bh^m+2c_K(C_\epsilon+C)r$,
$r=r_nh^{-1/2}$.
By the condition $h^{-3/2}r_n=o(1)$ and $B(h)\lesssim h^{-(2m-1)/(4m)}$
it can be easily checked that the Rate Condition (\textbf{H}):
is satisfied (when $n$ becomes large)
with $M$ therein replaced by $C$.
For $1\le j\le s$, define test $\phi_{j,n}=I(\|\widehat{f}_{j,n}-f_0\|\ge \delta_n)$.
It follows by part (\ref{basic:thm:i}) of Theorem \ref{basic:thm} that for any $1\le j\le s$,
\[
E_{f_0}\{\phi_{j,n}\}=P_{f_0}(\|\widehat{f}_{j,n}-f_0\|\ge \delta_n)
\le 2\exp(-Cnr_n^2),
\]
and 
\begin{eqnarray*}
\sup_{\substack{f\in H^m(b)\\
\|f-f_0\|\ge 2\delta_n}}E_f\{1-\phi_{j,n}\}&=&\sup_{\substack{f\in H^m(b)\\
\|f-f_0\|\ge 2\delta_n}}P_f(\|\widehat{f}_{j,n}-f_0\|<\delta_n)\\
&\le&\sup_{\substack{f\in H^m(b)\\
\|f-f_0\|\ge 2\delta_n}}P_f(\|\widehat{f}_{j,n}-f\|\ge \delta_n)
\le2\exp(-Cnr_n^2).
\end{eqnarray*}
An immediate consequence is 
$E_{f_0}\{\max_{1\le j\le s}\phi_{j,n}\}\le 2s\exp(-Cnr_n^2)$,
which implies $\max_{1\le j\le s}\phi_{j,n}=O_{P_{f_0}}(s\exp(-Cnr_n^2))$.

Note that for any $f\in A_n\backslash H^m(b)$,
$J(f)>b^2$.
Since $nh^{2m+1}\ge1$ leads to $r_n^2=(nh)^{-1}+\lambda\le 2\lambda$,
it then holds that, for any $1\le j\le s$,
\begin{eqnarray*}
&&E_{f_0}\{I_{j2}(1-\phi_{j,n})\}\\
&=&\int_{A_n}\|f-f_0\|^a E_f\{1-\phi_{j,n}\}\exp(-n\lambda J(f)/2)d\Pi(f)\\
&=&\int_{A_n\backslash H^m(b)}\|f-f_0\|^a E_f\{1-\phi_{j,n}\}\exp(-n\lambda J(f)/2)d\Pi(f)\\
&&+\int_{A_n\cap H^m(b)}\|f-f_0\|^a E_f\{1-\phi_{j,n}\}\exp(-n\lambda J(f)/2)d\Pi(f)\\
&\le&\left(\exp(-b^2n\lambda/2)+2\exp(-Cnr_n^2)\right)C(a,\Pi)\\
&\le &3\exp(-b^2nr_n^2/4)C(a,\Pi),
\end{eqnarray*}
where the last inequality follows by $C\ge b^2/4$ and $\lambda\ge r_n^2/2$.
So
\[
E_{f_0}\{\max_{1\le j\le s}I_{j2}(1-\phi_{j,n})\}
\le\sum_{j=1}^sE_{f_0}\{I_{j2}(1-\phi_{j,n})\}\le 3s\exp(-b^2nr_n^2/4)C(a,\Pi),
\]
which implies $\max_{1\le j\le s}I_{j2}(1-\phi_{j,n})=O_{P_{f_0}}(s\exp(-b^2nr_n^2/4))$.
On the other hand, as $n\rightarrow\infty$,
\begin{eqnarray*}
E_{f_0}\{\max_{1\le j\le s}I_{j2}\}\le s\int_{S^m(\mathbb{I})}\|f-f_0\|^2 d\Pi(f)
\end{eqnarray*}
which implies that $\max_{1\le j\le s}I_{j2}=o_{P_{f_0}}(s)$.
Therefore,
\begin{equation}\label{initial:rate:eqn6}
\max_{1\le j\le s}\frac{I_{j2}}{I_{j1}}\phi_{j,n}
\le\frac{\max_{1\le j\le s} I_{j2}\times\max_{1\le j\le s}\phi_{j,n}}{\min_{1\le j\le s}I_{j1}}
=O_{P_{f_0}}(s^2\exp(-nr_n^2)).
\end{equation}

By the above arguments and (\ref{initial:rate:eqn2}), we have
\begin{eqnarray*}
\max_{1\le j\le s}\int_{A_n}\|f-f_0\|^adP(f|\textbf{D}_j)
&=&\max_{1\le j\le s}\frac{I_{j2}}{I_{j1}}\\
&\le&\max_{1\le j\le s}\frac{I_{j2}}{I_{j1}}\phi_{j,n}
+\max_{1\le j\le s}\frac{I_{j2}(1-\phi_{j,n})}{I_{j1}}\\
&=&O_{P_{f_0}}(s^2\exp(-nr_n^2))+O_{P_{f_0}}(s\exp(-b^2nr_n^2/4)\exp(c_2nr_n^2))\\
&=&O_{P_{f_0}}(s^2\exp(-nr_n^2)).
\end{eqnarray*}
By condition $r_nh^{-3/2}=o(1)$ and the trivial fact $\delta_n\asymp r_nh^{-1/2}$,
we have that $h^{-1/2}\delta_n=o(1)$.
Therefore, eventually
$\int_{\|f-f_0\|_\infty\ge\varepsilon}\|f-f_0\|^adP(f|\textbf{D}_j)
\le\int_{A_n}\|f-f_0\|^adP(f|\textbf{D}_j)$ for all $1\le j\le s$,
which implies that (\ref{posterior:consistency}) holds.

Now we will prove the theorem. Let $I_{j2}'$ be defined as in (\ref{term:Ij2'})
with $A_n'=\{f\in S^m(\mathbb{I}):\|f-f_0\|\ge \sqrt{2}Mr_n\}$
for a fixed number satisfying
$M>\max\{2,J(f_0)^{1/2}+\sqrt{2(c_2+1)},1+\|f_0\|_\infty\}$
($M$ will be further described).
Let
$A_{n1}'=\{f\in S^m(\mathbb{I}):V(f-f_0)\ge M^2r_n^2,
\lambda J(f-f_0)\le M^2 r_n^2\}$
and
$A_{n2}'=\{f\in S^m(\mathbb{I}):\lambda J(f-f_0)\ge M^2 r_n^2\}$.
For any $f\in A_{n2}'$, it can be shown that
\[
Mr_n\le \sqrt{\lambda J(f-f_0)}\le \sqrt{\lambda} (J(f)^{1/2}+J(f_0)^{1/2})\le
(\lambda J(f))^{1/2}+J(f_0)^{1/2} r_n,
\]
which leads to $\lambda J(f)\ge (M-J(f_0)^{1/2})^2r_n^2$. So we have
\begin{eqnarray*}
&&E_{f_0}\{\max_{1\le j\le s}\int_{A_{n2}'}\|f-f_0\|^a\prod_{i\in I_j}
(p_f/p_{f_0})(Z_i)\exp(-\frac{n\lambda}{2}J(f))d\Pi(f)\}\\
&\le&\sum_{j=1}^sE_{f_0}\{\int_{A_{n2}'}\|f-f_0\|^a\prod_{i\in I_j} (p_f/p_{f_0})(Z_i)\exp(-\frac{n\lambda}{2}J(f))d\Pi(f)\}\\
&=&s\int_{A_{n2}'}\|f-f_0\|^a\exp(-\frac{n\lambda}{2}J(f))d\Pi(f)\}\\
&\le& s\exp(-(M-J(f_0)^{1/2})^2nr_n^2/2)C(a,\Pi),
\end{eqnarray*}
which leads to that
\begin{eqnarray}\label{initial:rate:eqn3}
&&\max_{1\le j\le s}\int_{A_{n2}'}\|f-f_0\|^a
\prod_{i\in I_j} (p_f/p_{f_0})(Z_i)\exp(-\frac{n\lambda}{2}J(f))d\Pi(f)\nonumber\\
&=&O_{P_{f_0}}(s\exp(-(M-J(f_0)^{1/2})^2nr_n^2/2)).
\end{eqnarray}
It follows from (\ref{initial:rate:eqn2}) and (\ref{initial:rate:eqn3}) that
\begin{eqnarray}\label{initial:rate:eqn5}
\max_{1\le j\le s}\frac{1}{I_{j1}}\int_{A_{n2}'}\|f-f_0\|^a
\prod_{i\in I_j} (p_f/p_{f_0})(Z_i)\exp(-\frac{n\lambda}{2}J(f))d\Pi(f)\nonumber\\
=O_{P_{f_0}}\left(s\exp(-(M-J(f_0)^{1/2})^2nr_n^2/2+c_2nr_n^2)\right)=O_{P_{f_0}}(s\exp(-nr_n^2)),
\end{eqnarray}
where the last inequality follows by
$(M-J(f_0)^{1/2})^2>2(c_2+1)$.

To continue, we need to build uniformly consistent test.
Let $d_H^2(P_f,P_g)=\frac{1}{2}\int (\sqrt{dP_f}-\sqrt{dP_g})^2$
be the squared Hellinger distance between the
two probability measures $P_f(z)$ and $P_g(z)$.
Recall that their corresponding probability density functions are
$p_f$ and $p_g$, respectively.
Nextwe  present a lemma showing the local equivalence of $V$ and $d_H^2$.

\begin{lemma}\label{local:equiv:V:d_H^2}
Let $\varepsilon\in(0,1)$ satisfy
$\varepsilon^2+32\varepsilon\exp(1/2)C_\epsilon\le2$,
where $C_\epsilon=E\{\exp(|\epsilon|)\}$.
Then for any $f,g\in S^m(\mathbb{I})$ satisfying $\|f-g\|_\infty\le\varepsilon$,
$V(f-g)/16\le d_H^2(P_f,P_g)\le 3V(f-g)/16$.
\end{lemma}

Let $\varepsilon$ satisfy the conditions in Lemma \ref{local:equiv:V:d_H^2}.
Define $\mathcal{F}_n=\{f\in S^m(\mathbb{I}): \|f-f_0\|_\infty\le\varepsilon/2,
J(f)\le (M+J(f_0)^{1/2})^2r_n^2\lambda^{-1}\}$.
Let $\mathcal{P}_n=\{P_f: f\in\mathcal{F}_n\}$
and $D(\delta,\mathcal{P}_n,d_H)$ be the $\delta$-packing number in terms of $d_H$.
Since $r_n^2\ge \lambda$ which leads to
$(M+J(f_0)^{1/2})r_nh^{-m}>M+J(f_0)^{1/2}>\varepsilon+\|f_0\|_\infty$,
it can be easily checked that $\mathcal{F}_n\subset (M+J(f_0)^{1/2})r_nh^{-m}\mathcal{T}$,
where $\mathcal{T}=\{f\in S^m(\mathbb{I}): \|f\|_\infty\le1, J(f)\le 1\}$.

For any $f,g\in\mathcal{F}_n$ with $\|f-g\|_\infty\le\varepsilon$,
it follows by Lemma \ref{local:equiv:V:d_H^2} that
$D(\delta,\mathcal{P}_n,d_H)\le D(4\delta/\sqrt{3},\mathcal{F}_n,d_V)$,
where $d_V$ is the distance induced by $V$, i.e., $d_V(f,g)=V^{1/2}(f-g)$.
And hence, it follows by \cite[Theorem 9.21]{K08} that
\begin{eqnarray*}
\log{D(\delta,\mathcal{P}_n,d_H)}&\le& \log{D(4\delta/\sqrt{3},\mathcal{F}_n,d_V)}\\
&\le&\log{D(4\delta/\sqrt{3},(M+J(f_0)^{1/2})r_nh^{-m}\mathcal{T},d_V)}\\
&\le& c_V \left(\frac{\delta}{(M+J(f_0)^{1/2})r_nh^{-m}}\right)^{-1/m},
\end{eqnarray*}
where $c_V$ is a universal constant only depending on the regularity level $m$.
This implies that for any $\delta>2r_n$,
\begin{eqnarray*}
\log{D(\delta/2,\mathcal{P}_n,d_H)}
&\le& \log{D(r_n,\mathcal{P}_n,d_H)}\\
&\le& c_V(M+J(f_0)^{1/2})^{1/m}h^{-1}\\
&\le& c_V(M+J(f_0)^{1/2})^{1/m}nr_n^2,
\end{eqnarray*}
where the last inequality follows by the fact $r_n^2\ge (nh)^{-1}$.
Thus, the right side of the above inequality is constant in $\delta$.
By \cite[Theorem 7.1]{GGV00}, with $\delta=Mr_n/4$, there exists
test $\widetilde{\phi}_{j,n}$ and a universal constant $k_0>0$ satisfying
\begin{eqnarray*}
E_{f_0}\{\widetilde{\phi}_{j,n}\}&=&P_{f_0}\widetilde{\phi}_{j,n}\\
&\le& \frac{\exp(c_V(M+J(f_0)^{1/2})^{1/m}nr_n^2)\exp(-k_0n\delta^2)}
{1-\exp(-k_0n\delta^2)}\\
&=& \frac{\exp(c_V(M+J(f_0)^{1/2})nr_n^2-k_0M^2nr_n^2/16)}{1-\exp(-k_0M^2nr_n^2/16)},
\end{eqnarray*}
and, combined with Lemma \ref{local:equiv:V:d_H^2},
\begin{eqnarray*}
\sup_{\substack{f\in\mathcal{F}_n\\ d_V(f,f_0)\ge 4\delta}}E_{f}\{1-\widetilde{\phi}_{j,n}\}
&=&\sup_{\substack{f\in\mathcal{F}_n\\ d_V(f,f_0)\ge 4\delta}}P_f\{1-\widetilde{\phi}_{j,n}\}\\
&\le&\sup_{\substack{f\in\mathcal{F}_n\\ d_H(P_f,P_{f_0})\ge \delta}}P_f\{1-\widetilde{\phi}_{j,n}\}\\
&\le&\exp(-k_0n\delta^2)=\exp(-k_0M^2nr_n^2/16).
\end{eqnarray*}
This implies that
\begin{eqnarray*}
&&E_{f_0}\{\max_{1\le j\le s}\int_{\substack{f\in\mathcal{F}_n\\
d_V(f,f_0)\ge 4\delta}}\|f-f_0\|^a\prod_{i\in I_j}(p_f/p_{f_0})(Z_i)\exp(-n\lambda J(f)/2)d\Pi(f)
(1-\widetilde{\phi}_{j,n})\}\\
&\le&\sum_{j=1}^s\int_{\substack{f\in\mathcal{F}_n\\
d_V(f,f_0)\ge 4\delta}}\|f-f_0\|^aE_{f_0}\{\prod_{i\in I_j}(p_f/p_{f_0})(Z_i)(1-\widetilde{\phi}_{j,n})\}d\Pi(f)\\
&=&\sum_{j=1}^s\int_{\substack{f\in\mathcal{F}_n\\
d_V(f,f_0)\ge 4\delta}}\|f-f_0\|^a E_f\{1-\widetilde{\phi}_{j,n}\}d\Pi(f)\\
&\le&s\exp(-k_0M^2nr_n^2/16)C(a,\Pi).
\end{eqnarray*}
Therefore,
\begin{eqnarray}\label{initial:rate:eqn4}
&&\max_{1\le j\le s}\int_{\substack{f\in\mathcal{F}_n\\
d_V(f,f_0)\ge 4\delta}}\|f-f_0\|^a\prod_{i\in I_j}(p_f/p_{f_0})(Z_i)\exp(-n\lambda J(f)/2)d\Pi(f)
(1-\widetilde{\phi}_{j,n})\nonumber\\
&=&O_{P_{f_0}}\left(s\exp(-k_0M^2nr_n^2/16)\right).\nonumber\\
\end{eqnarray}
Meanwhile, it follows by (\ref{initial:rate:eqn2}) and (\ref{initial:rate:eqn4}) that
\begin{eqnarray*}
&&\max_{1\le j\le s} 
\int_{A_{n1}',\|f-f_0\|_\infty\le\varepsilon/2}
\|f-f_0\|^adP(f|\textbf{D}_j)(1-\widetilde{\phi}_{j,n})\\
&\le&\max_{1\le j\le s} 
\int_{\mathcal{F}_n,d_V(f,f_0)\ge 4\delta}
\|f-f_0\|^adP(f|\textbf{D}_j)(1-\widetilde{\phi}_{j,n})\\
&\le&\frac{\max\limits_{1\le j\le s}\int_{\substack{f\in\mathcal{F}_n\\
d_V(f,f_0)\ge 4\delta}}\|f-f_0\|^a\prod_{i\in I_j}(p_f/p_{f_0})(Z_i)\exp(-n\lambda J(f)/2)d\Pi(f)
(1-\widetilde{\phi}_{j,n})}{\min\limits_{1\le j\le s}I_{j1}}\\
&=&O_{P_{f_0}}\left(s\exp(-k_0M^2nr_n^2/16+c_2nr_n^2)\right)
=O_{P_{f_0}}\left(s\exp(-nr_n^2)\right).
\end{eqnarray*}
Choose the constant $M$ to be even bigger so that
$c_V(M+J(f_0)^{1/2})+1+c_2<k_0M^2/16$.
Similar to (\ref{initial:rate:eqn6}) we get
\begin{eqnarray*}
&&\max_{1\le j\le s} \int_{A_{n1}',\|f-f_0\|_\infty\le\varepsilon/2}
\|f-f_0\|^adP(f|\textbf{D}_j)\widetilde{\phi}_{j,n}
=O_{P_{f_0}}(s^2\exp(-nr_n^2)).
\end{eqnarray*}
Therefore,
\begin{eqnarray}\label{initial:rate:eqn7}
&&\max_{1\le j\le s} \int_{A_{n1}',\|f-f_0\|_\infty\le\varepsilon/2}
\|f-f_0\|^adP(f|\textbf{D}_j)=O_{P_{f_0}}(s^2\exp(-nr_n^2)).
\end{eqnarray}
Together with (\ref{posterior:consistency}), (\ref{initial:rate:eqn6}) and (\ref{initial:rate:eqn7}),
we get
\begin{eqnarray*}
&&\max_{1\le j\le s}\int_{A_n'}\|f-f_0\|^adP(f|\textbf{D}_j)\\
&\le&\max_{1\le j\le s} \int_{A_{n1}'}\|f-f_0\|^adP(f|\textbf{D}_j)+
\max_{1\le j\le s} \int_{A_{n2}'}\|f-f_0\|^adP(f|\textbf{D}_j)\\
&\le&\max_{1\le j\le s} \int_{A_{n1}',\|f-f_0\|_\infty\le\varepsilon/2}
\|f-f_0\|^adP(f|\textbf{D}_j)
+\max_{1\le j\le s}\int_{\|f-f_0\|_\infty>\varepsilon/2}\|f-f_0\|^adP(f|\textbf{D}_j)\\
&&+\max_{1\le j\le s} \int_{A_{n2}'}\|f-f_0\|^adP(f|\textbf{D}_j)\\
&=&O_{P_{f_0}}(s^2\exp(-nr_n^2)).
\end{eqnarray*}
This completes the proof.
\end{proof}

\begin{proof}[Proof of Lemma \ref{local:equiv:V:d_H^2}]
For any $f,g\in S^m(\mathbb{I})$ with $\|f-g\|_\infty\le\varepsilon$,
define $\Delta_Z(f,g)=\frac{1}{2}[Y(f(X)-g(X))-f(X)^2/2+g(X)^2/2]$,
where recall and $Z=(Y,X)$.
It is easy to see by direct calculations that
$d_H^2(P_f,P_g)=1-E_g\{\exp(\Delta_Z(f,g))\}$.
By Taylor's expansion, for some random $t\in[0,1]$,
\begin{eqnarray*}
&&1-E_g\{\exp(\Delta_Z(f,g))\}\\
&=&-E_g\{\Delta_Z(f,g)\}-\frac{1}{2}E_g\{\Delta_Z(f,g)^2\}
-\frac{1}{6}E_g\{\exp(t\Delta_Z(f,g))\Delta_Z(f,g)^3\}.
\end{eqnarray*}
We will analyze the terms on the right side of the equation.

Define $\xi=Y-\dot{A}(g(X))$.
By \cite{M82} we get $E_g\{\xi|X\}=0$ and
$E_g\{\xi^2|X\}=1$.
By Taylor's expansion,
$\Delta_Z(f,g)=\frac{1}{2}[\xi(f(X)-g(X))-\frac{1}{2}(f(X)-g(X))^2$.
Then we get that
$-E_g\{\Delta_Z(f,g)\}=\frac{1}{4}V(f-g)$
and
\begin{eqnarray*}
E_g\{\Delta_Z(f,g)^2\}&=&E_g\{(\frac{1}{2}\xi(f(X)-g(X))-\frac{1}{4}(f(X)-g(X))^2)\}\\
&=&\frac{1}{4}E_g\{\xi^2(f(X)-g(X))^2\}-\frac{1}{4}E_g\{\xi(f(X)-g(X))^3\}
+\frac{1}{16}E_g\{(f(X)-g(X))^4\}\\
&=&\frac{1}{4}V(f-g)+\frac{1}{16}E_g\{(f(X)-g(X))^4\}.
\end{eqnarray*}
Since $\|f-g\|_\infty\le\varepsilon<1$ and
$|\Delta_Z(f,g)|\le\frac{1}{2}(|\xi|+1/2)|f(X)-g(X)|$, we get
\begin{eqnarray*}
&&|E_g\{\exp(t\Delta_Z(f,g))\Delta_Z(f,g)^3\}|\\
&\le&E_g\{\exp(|\Delta_Z(f,g)|)|\Delta_Z(f,g)|^3\}\\
&\le&E_g\{\exp(\varepsilon |\xi|/2+\varepsilon/4)(|\xi|/2+1/4)^3|f(X)-g(X)|^3\}\\
&=&6E_g\left\{\exp(\varepsilon |\xi|/2+\varepsilon/4)\times\frac{1}{3!}
\left(|\xi|/2+1/4\right)^3|f(X)-g(X)|^3\right\}\\
&\le&6E_g\{\exp(\varepsilon |\xi|/2+\varepsilon/4)\exp(|\xi|/2+1/4)|f(X)-g(X)|^3\}\\
&\le&6\exp(\varepsilon/4+/4)E_g\{\exp(|\xi|)|f(X)-g(X)|^3\}\\
&\le&6\varepsilon\exp(1/2)C_\epsilon V(f-g).
\end{eqnarray*}
It also holds that
$|E_g\{(f(X)-g(X))^4\}|\le \varepsilon^2 V(f-g)$.
Therefore, for any $f,g\in S^m(\mathbb{I})$
with $\|f-g\|_\infty\le\varepsilon$,
\begin{eqnarray*}
&&|d_H^2(P_f,P_g)-V(f-g)/8|\\
&=&|\frac{1}{32}E_g\{(f(X)-g(X))^4\}
+\frac{1}{6}E_g\{\exp(t\Delta_Z(f,g))\Delta_Z(f,g)^3\}|\\
&\le&\left(\varepsilon C_\epsilon\exp(1/2)+\varepsilon^2/32\right)V(f-g)
<V(f-g)/16,
\end{eqnarray*}
which implies $V(f-g)/16\le d_H^2(P_f,P_g)\le 3V(f-g)/16$. This proves Lemma \ref{local:equiv:V:d_H^2}.
\end{proof}

\newpage

\subsection{Additional Plots in Section \ref{sec:simulations}}\label{sec:suppl:additional:simulation}

\subsubsection*{Radius of the credible sets/intervals}\label{sec:suppl:additional:simulation:radius}

\begin{figure}[htp]
	\begin{center}
		\includegraphics[scale=0.65]{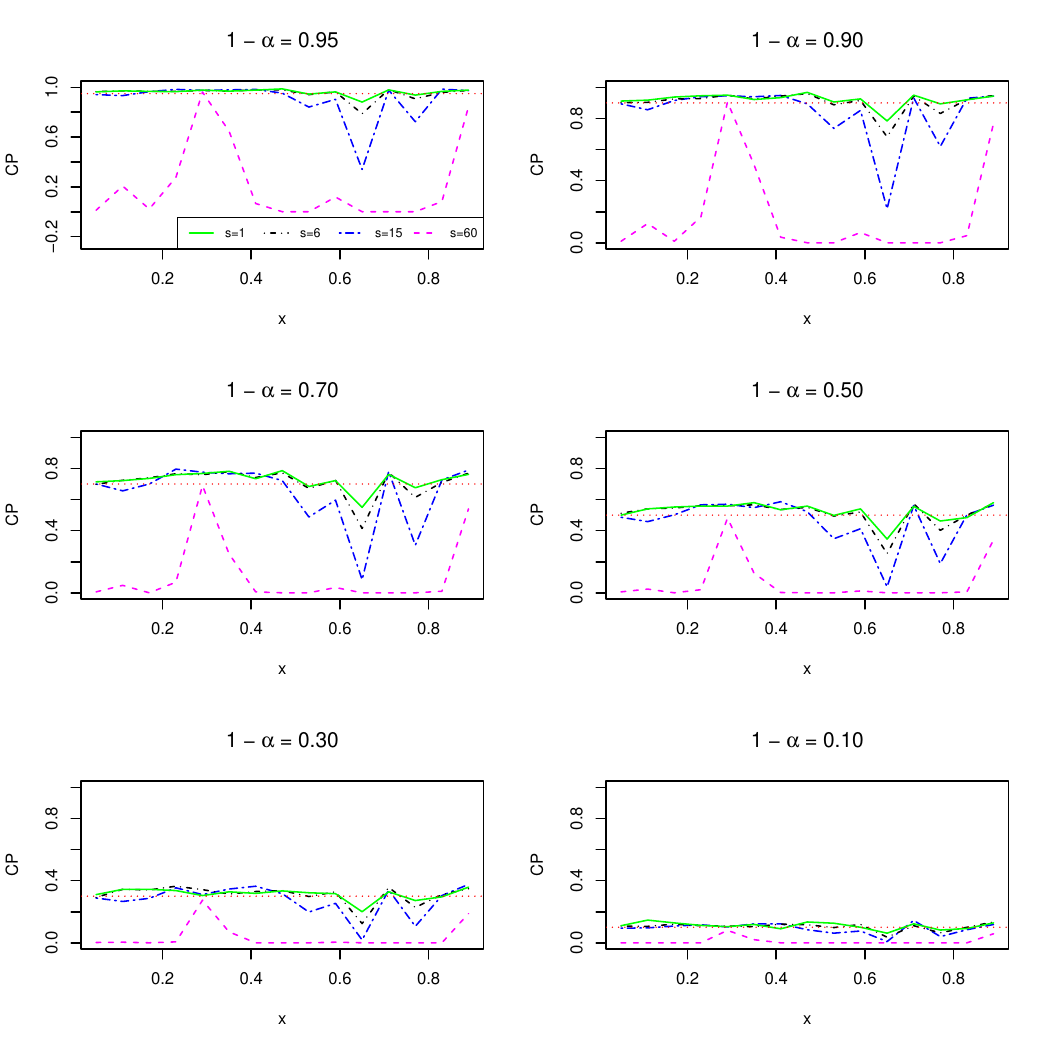}
		\caption{\textit{\footnotesize CP of $F_x(f)=f(x)$ against $x$ based on asymptotic theory.}}
		\label{fig:asymp:pointwise:CI}
	\end{center}
\end{figure}

\begin{figure}[htp]
	\begin{center}
		\includegraphics[scale=0.65]{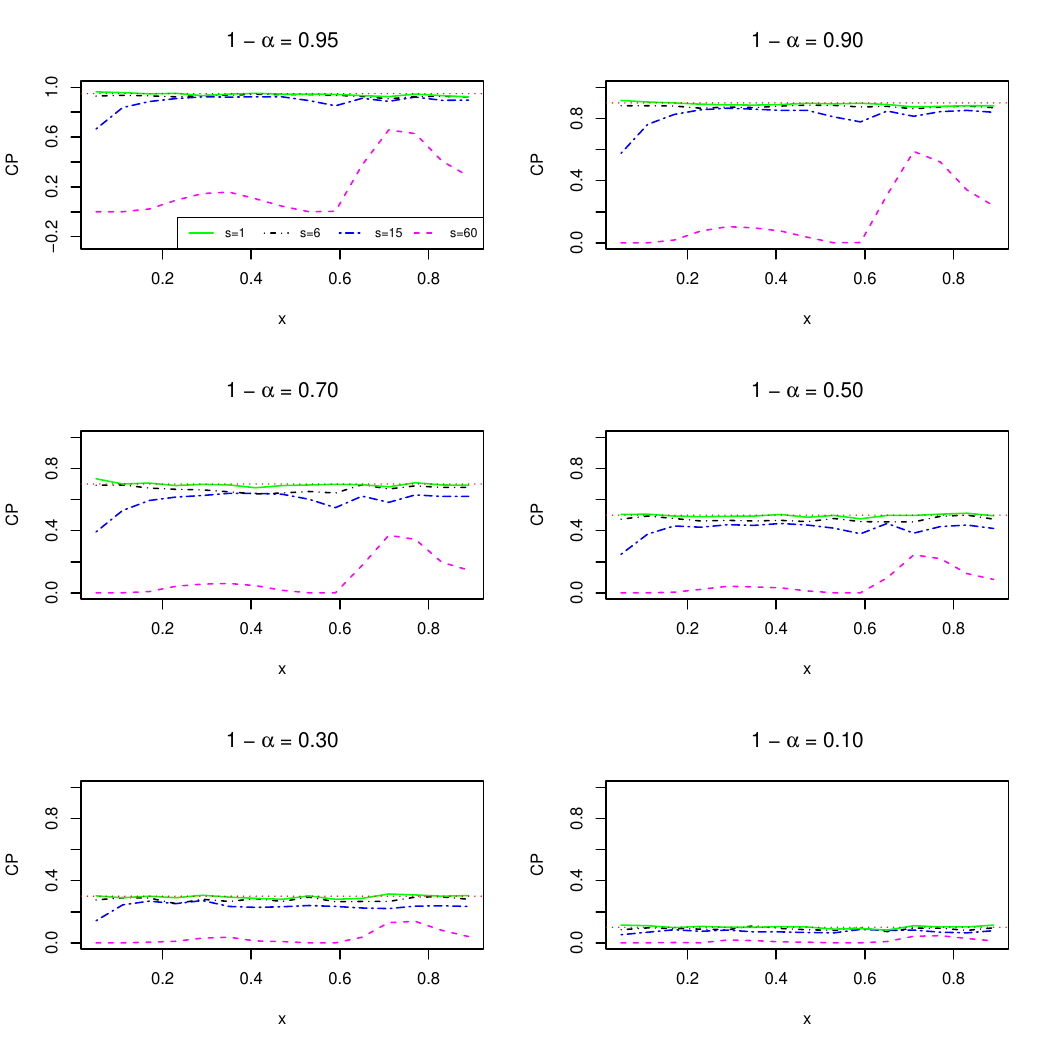}
		\caption{\textit{\footnotesize CP of $F_x(f)=\int_0^xf(z)dz$ against $x$ based on asymptotic theory.}}
		\label{fig:asymp:integral:CI}
	\end{center}
\end{figure}

\subsubsection*{Results on larger $N$}\label{sec:suppl:additional:simulation:larger:N}
Simulation results about credible regions/intervals in Section \ref{sec:simulations} are based on $N=1200$.
This section repeated the same study for $N=1800,2400$.
Results are summarized in following plots.

\begin{figure}[htp]
\begin{center}
\includegraphics[scale=0.65]{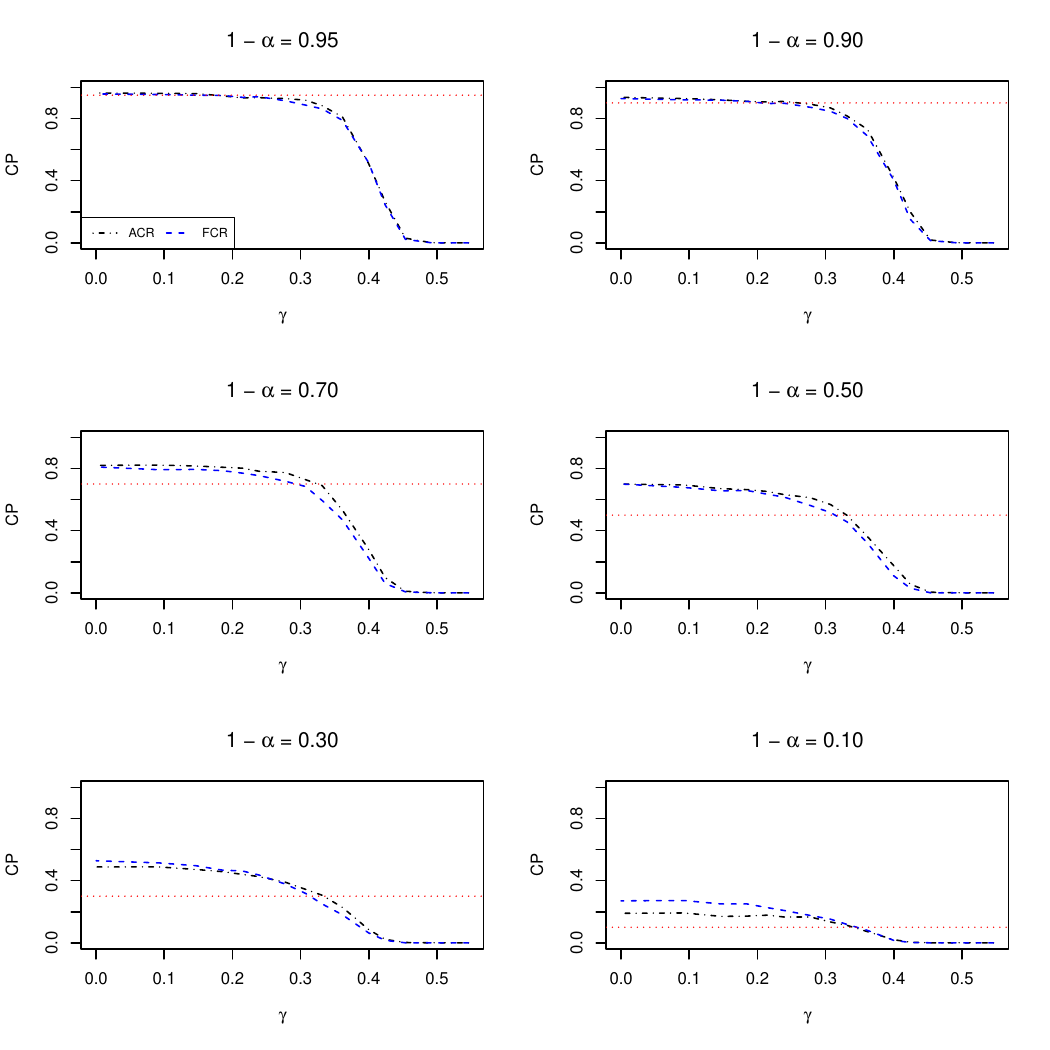}
\caption{\textit{\footnotesize $N=1800$: CP of ACR and FCR based on strong topology.}}
\label{fig:strong:cr:1800}
\end{center}
\end{figure}

\begin{figure}[htp]
\begin{center}
\includegraphics[scale=0.65]{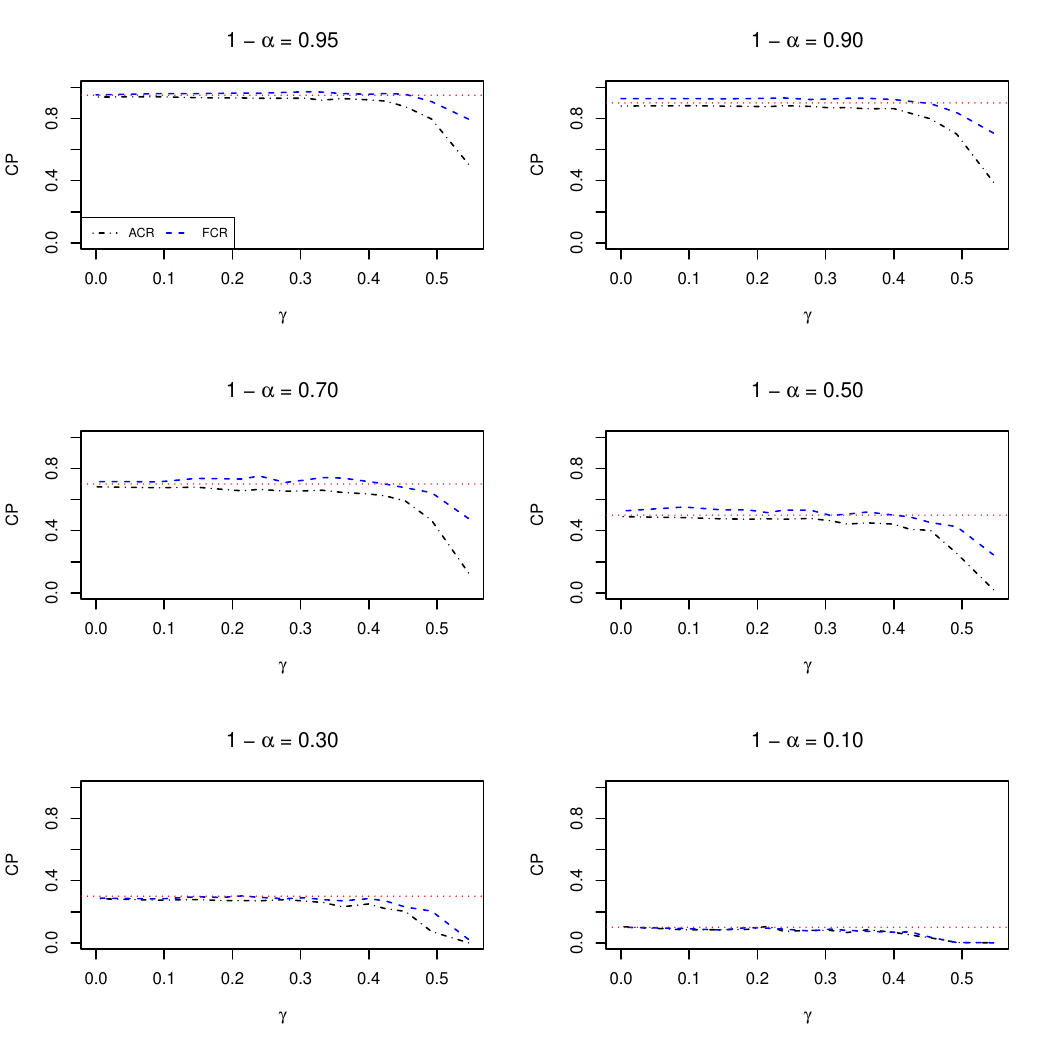}
\caption{\textit{\footnotesize $N=1800$: CP of ACR and FCR based on weak topology.}}
\label{fig:weak:cr:1800}
\end{center}
\end{figure}

\begin{figure}[htp]

\begin{center}

\includegraphics[scale=0.65]{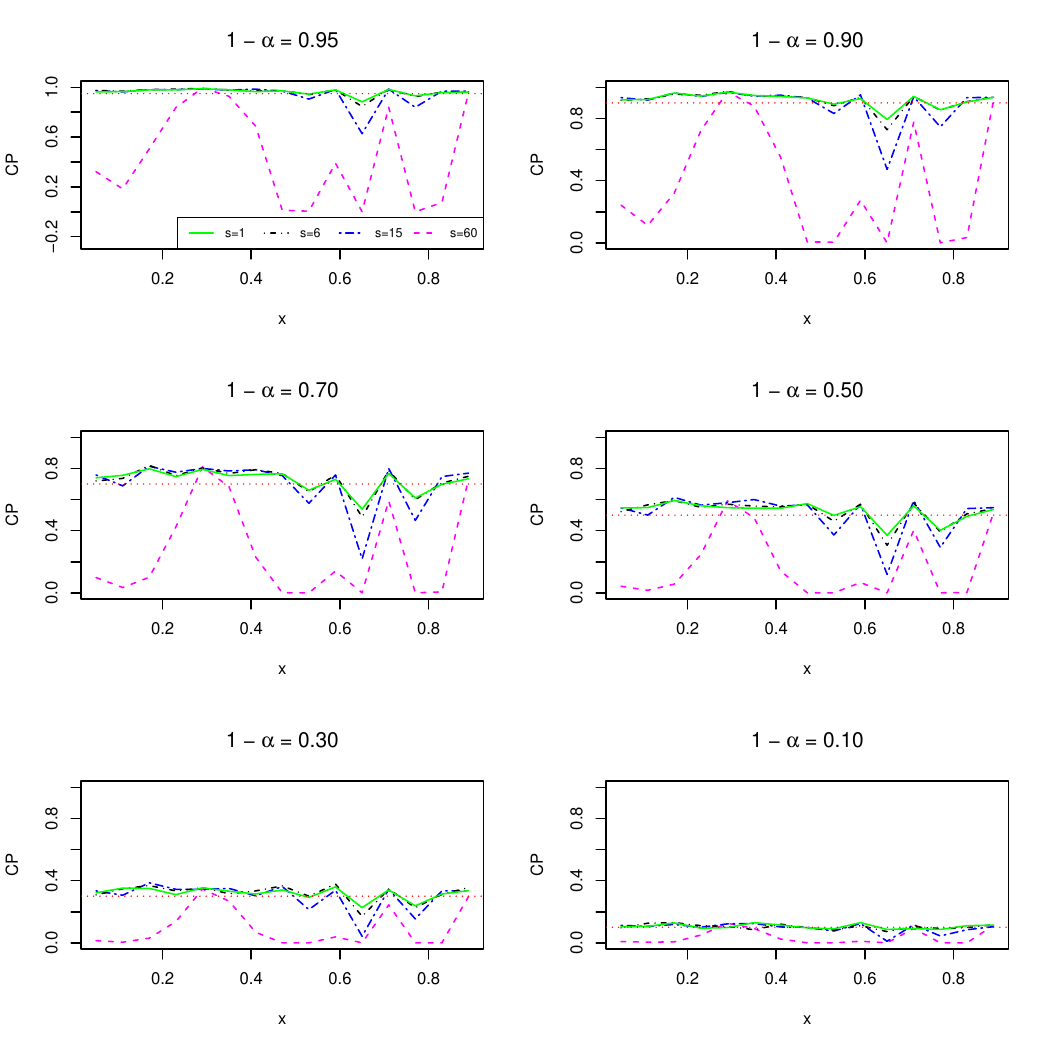}

\caption{\textit{\footnotesize $N=1800$: CP of $F_x(f)=f(x)$ against $x$ based on posterior samples of $f$.}}

\label{fig:finite:pointwise:CI:1800}

\end{center}

\end{figure}

\begin{figure}[htp]

\begin{center}

\includegraphics[scale=0.65]{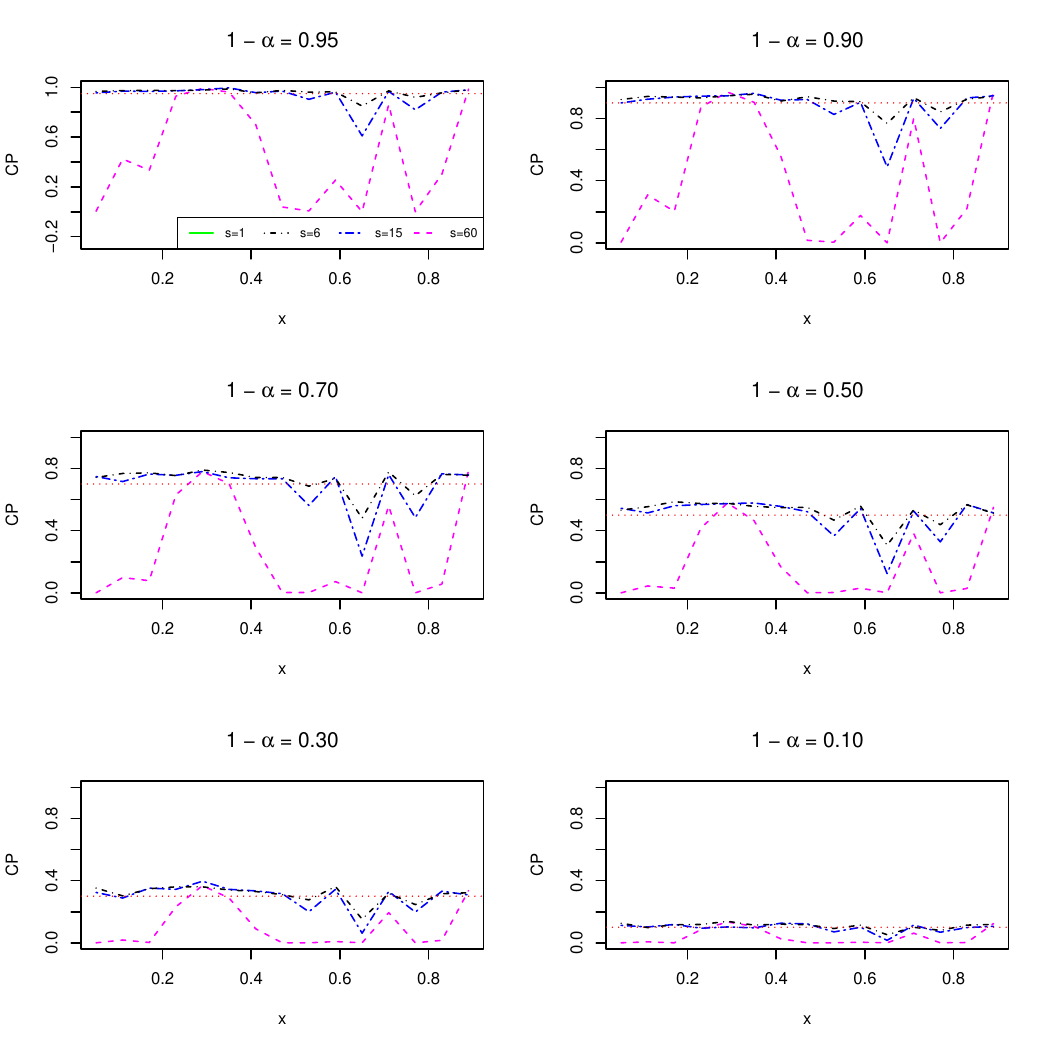}

\caption{\textit{\footnotesize $N=1800$: CP of $F_x(f)=f(x)$ against $x$ based on asymptotic theory.}}

\label{fig:asymp:pointwise:CI:1800}

\end{center}

\end{figure}

\begin{figure}[htp]

\begin{center}

\includegraphics[scale=0.65]{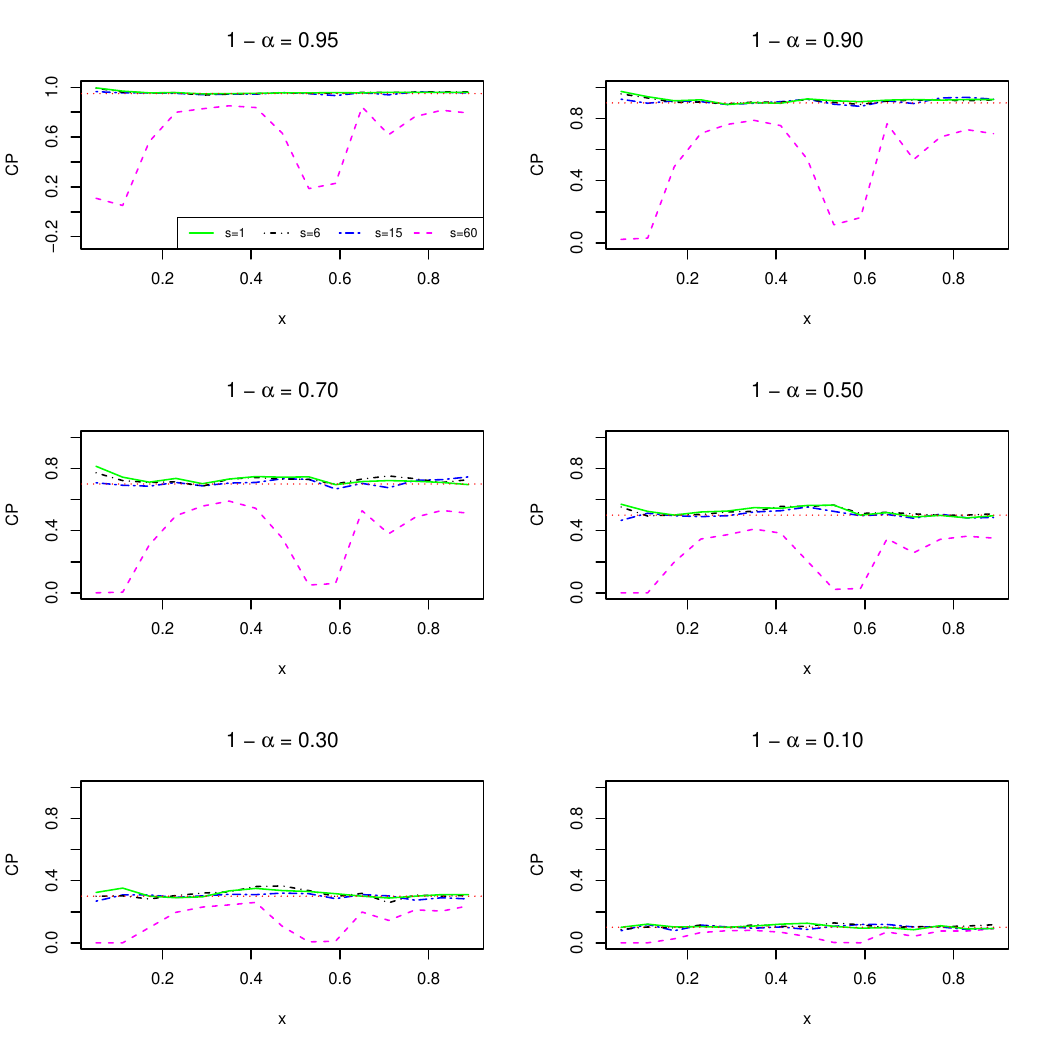}

\caption{\textit{\footnotesize $N=1800$: CP of $F_x(f)=\int_0^xf(z)dz$ against $x$ based on posterior samples of $f$.}}

\label{fig:finite:integral:CI:1800}

\end{center}

\end{figure}

\begin{figure}[htp]

\begin{center}

\includegraphics[scale=0.65]{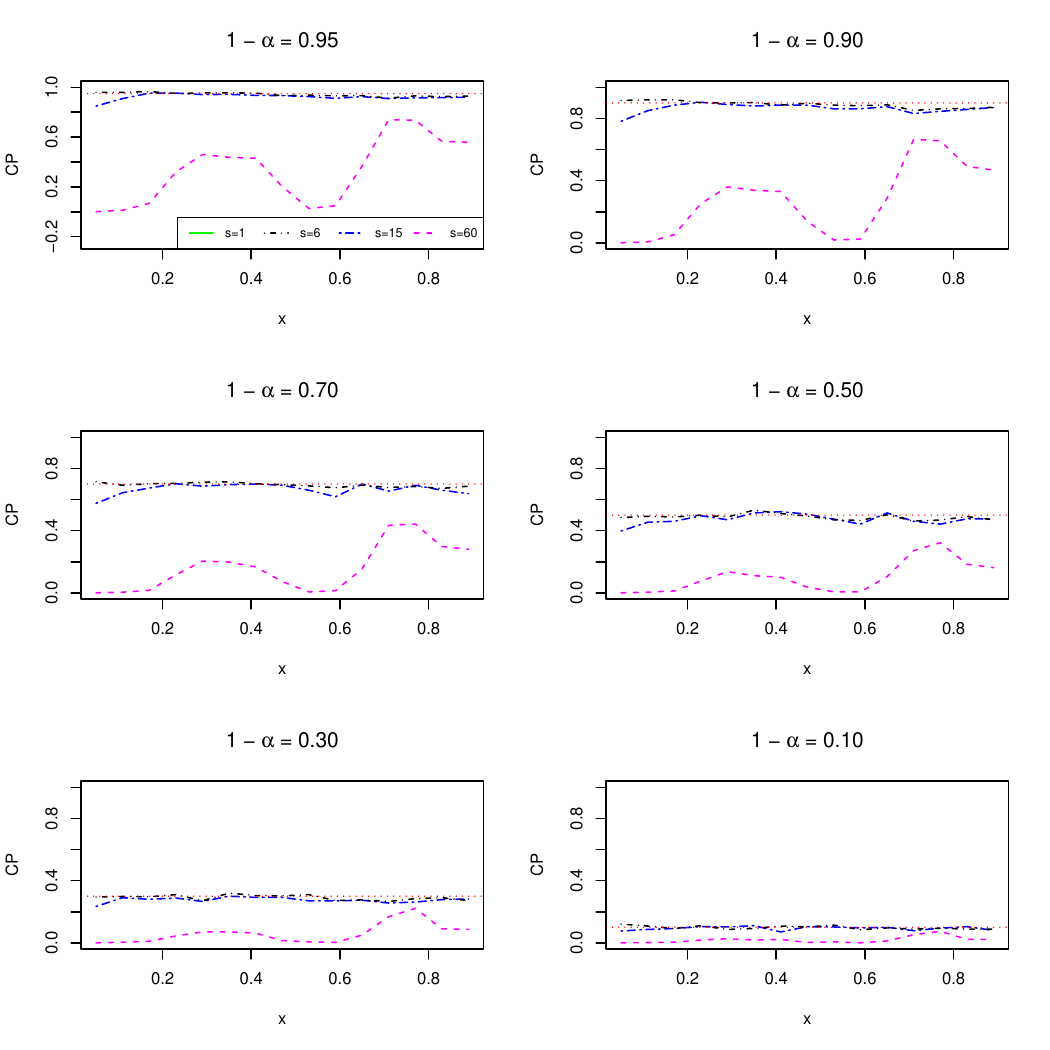}

\caption{\textit{\footnotesize $N=1800$: CP of $F_x(f)=\int_0^xf(z)dz$ against $x$ based on asymptotic theory.}}

\label{fig:asymp:integral:CI:1800}

\end{center}

\end{figure}

\begin{figure}[htp]
\begin{center}
\includegraphics[scale=0.65]{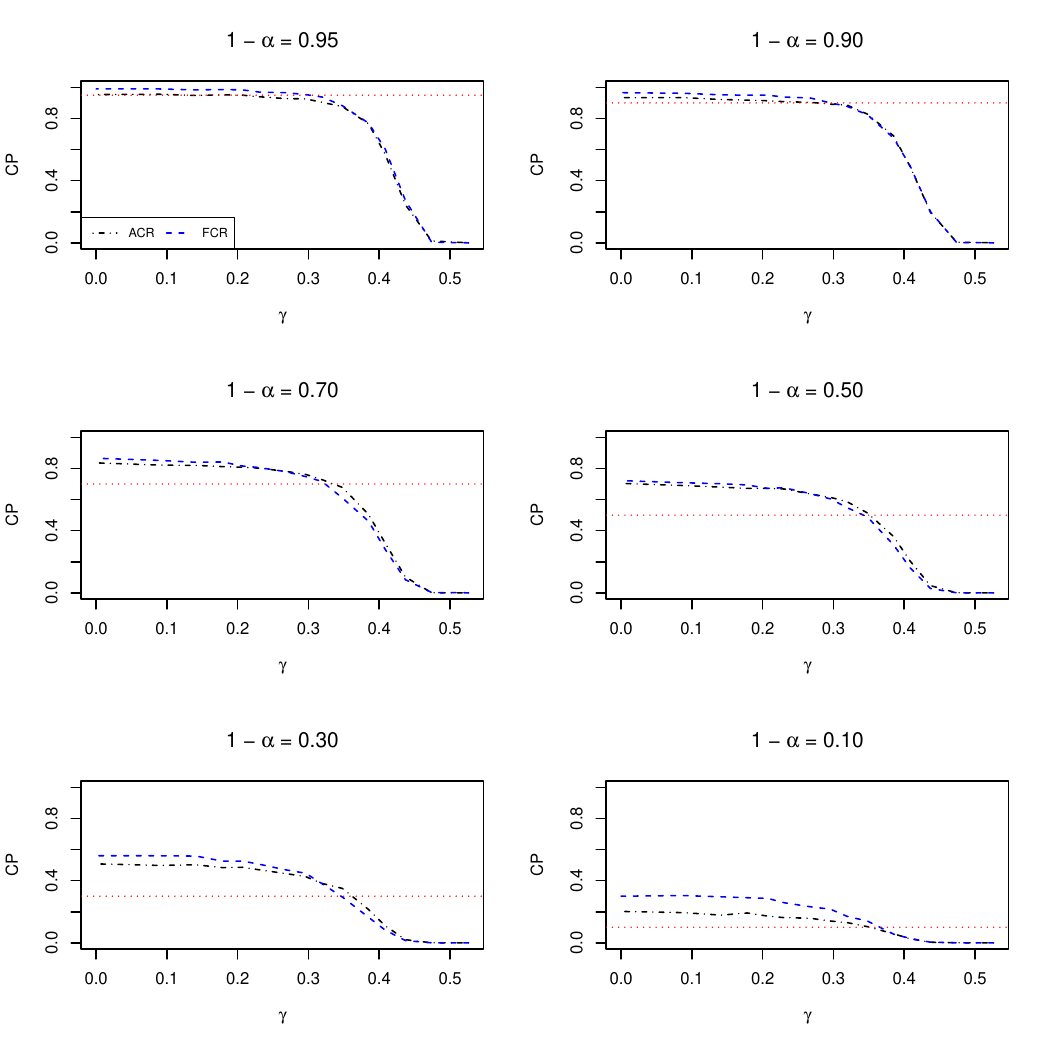}
\caption{\textit{\footnotesize $N=2400$: CP of ACR and FCR based on strong topology.}}
\label{fig:strong:cr:2400}
\end{center}
\end{figure}

\begin{figure}[htp]
\begin{center}
\includegraphics[scale=0.65]{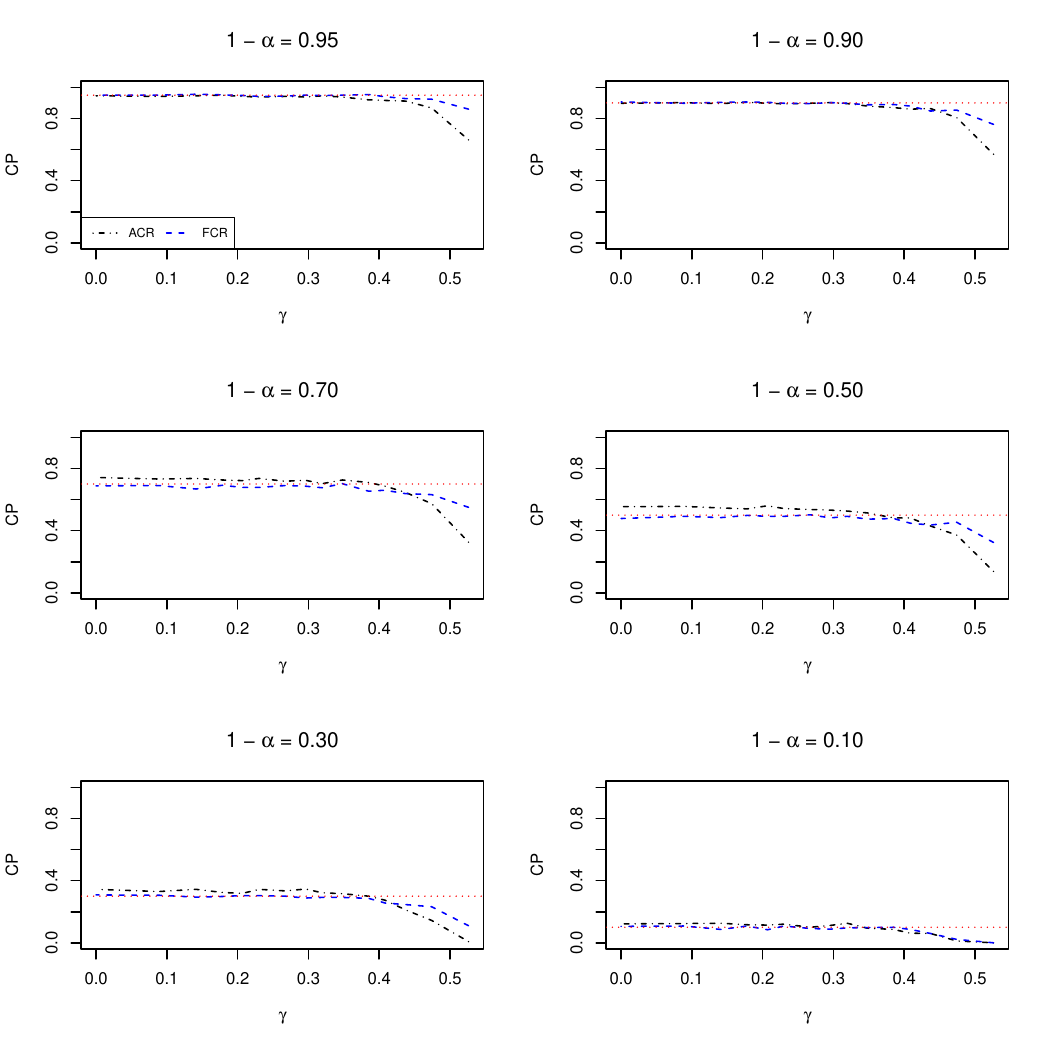}
\caption{\textit{\footnotesize $N=2400$: CP of ACR and FCR based on weak topology.}}
\label{fig:weak:cr:2400}
\end{center}
\end{figure}

\begin{figure}[htp]

\begin{center}

\includegraphics[scale=0.65]{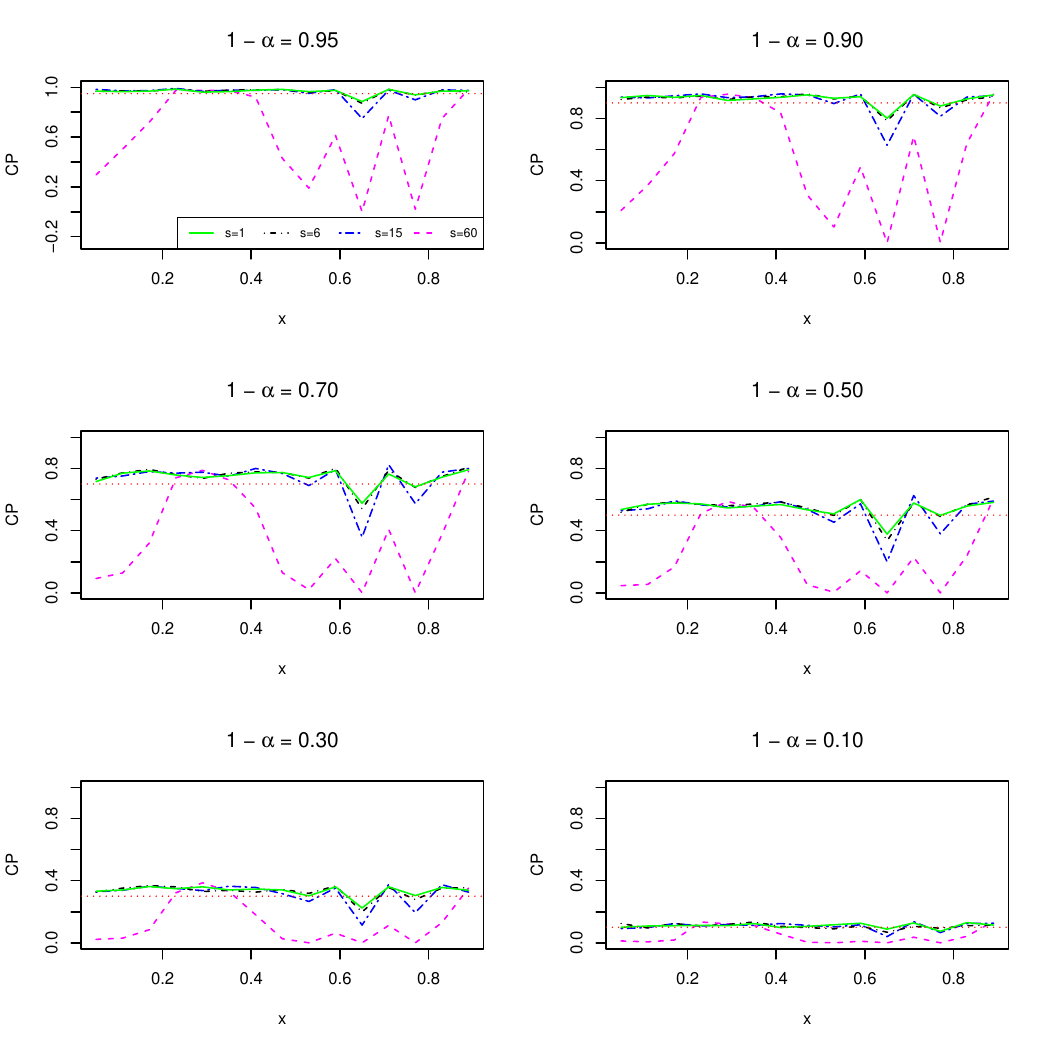}

\caption{\textit{\footnotesize $N=2400$: CP of $F_x(f)=f(x)$ against $x$ based on posterior samples of $f$.}}

\label{fig:finite:pointwise:CI:2400}

\end{center}

\end{figure}

\begin{figure}[htp]

\begin{center}

\includegraphics[scale=0.65]{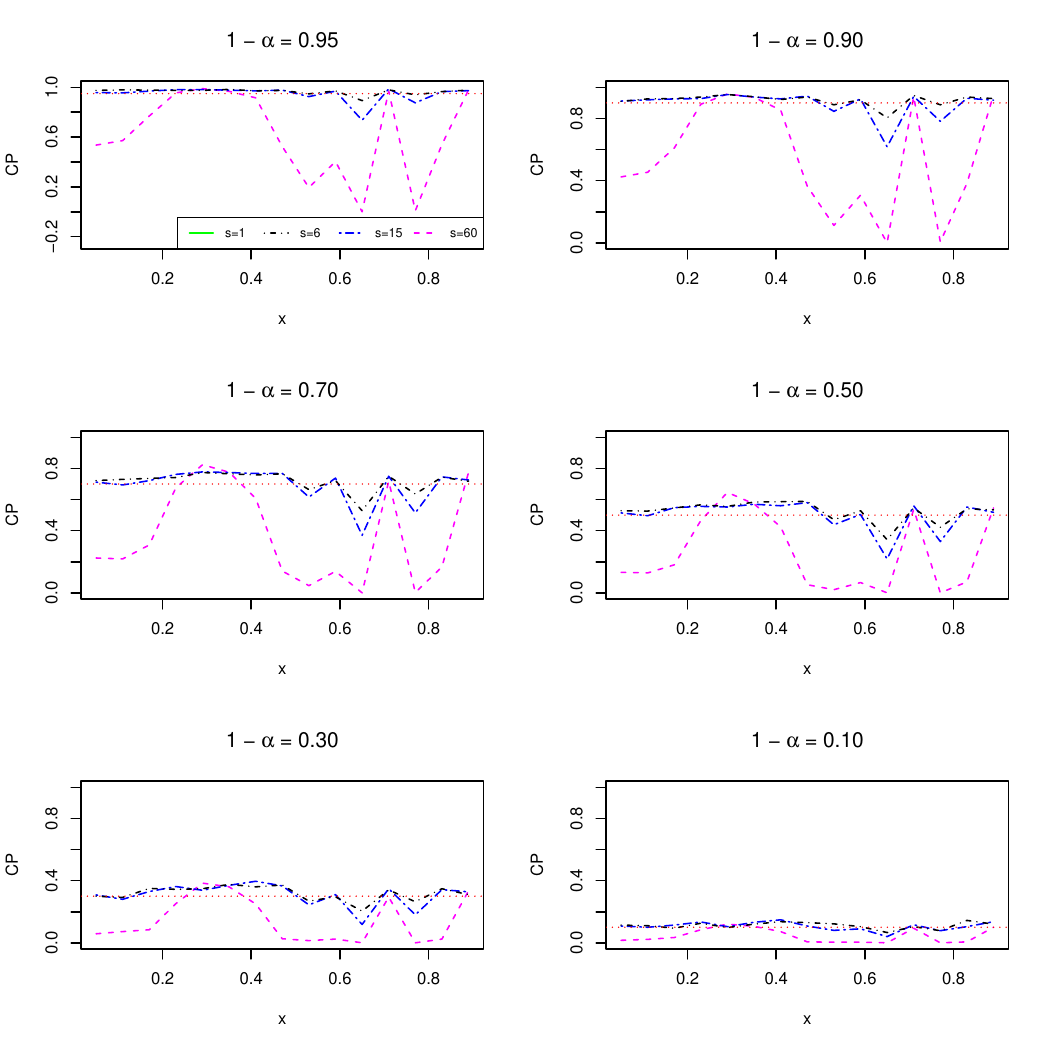}

\caption{\textit{\footnotesize $N=2400$: CP of $F_x(f)=f(x)$ against $x$ based on asymptotic theory.}}

\label{fig:asymp:pointwise:CI:2400}

\end{center}

\end{figure}

\begin{figure}[htp]

\begin{center}

\includegraphics[scale=0.65]{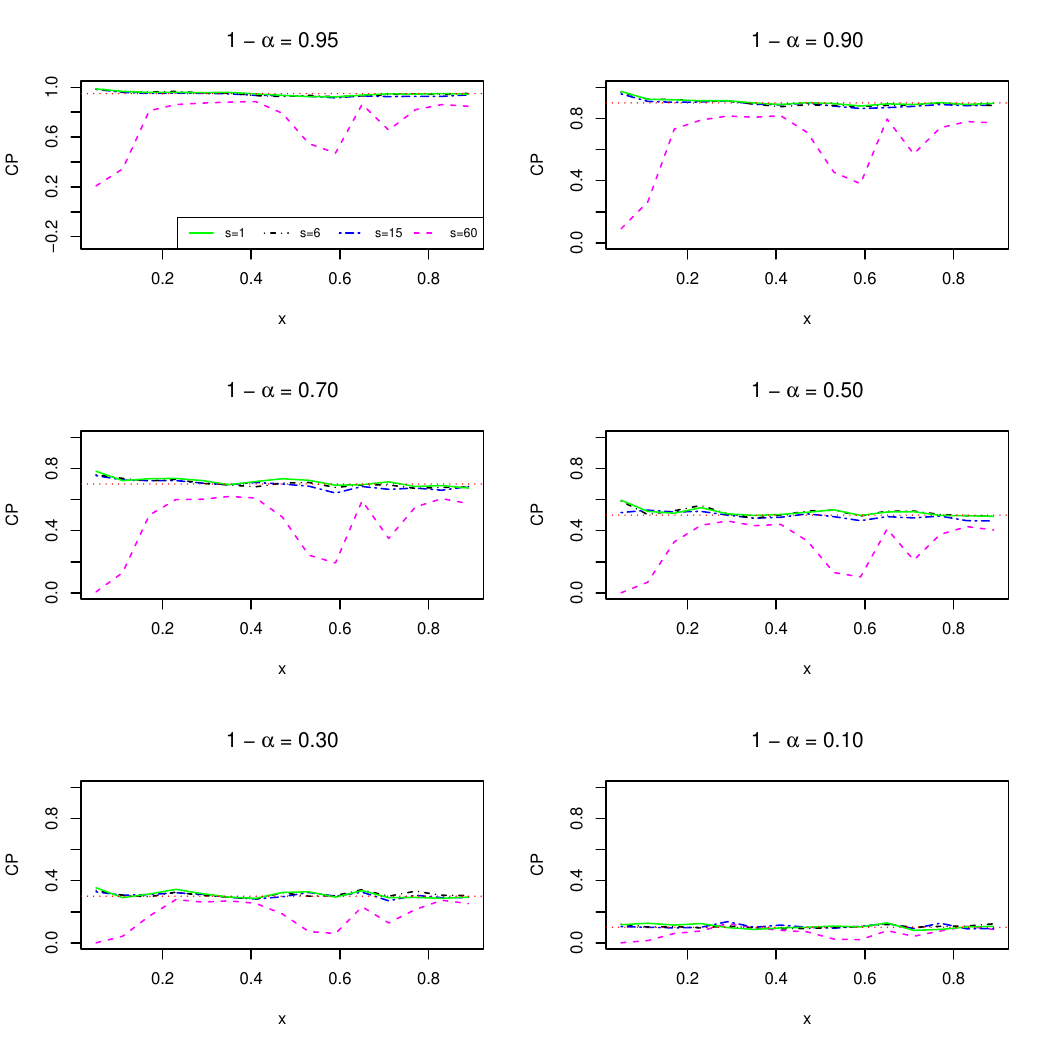}

\caption{\textit{\footnotesize $N=2400$: CP of $F_x(f)=\int_0^xf(z)dz$ against $x$ based on posterior samples of $f$.}}

\label{fig:finite:integral:CI:2400}

\end{center}

\end{figure}

\begin{figure}[htp]

\begin{center}

\includegraphics[scale=0.65]{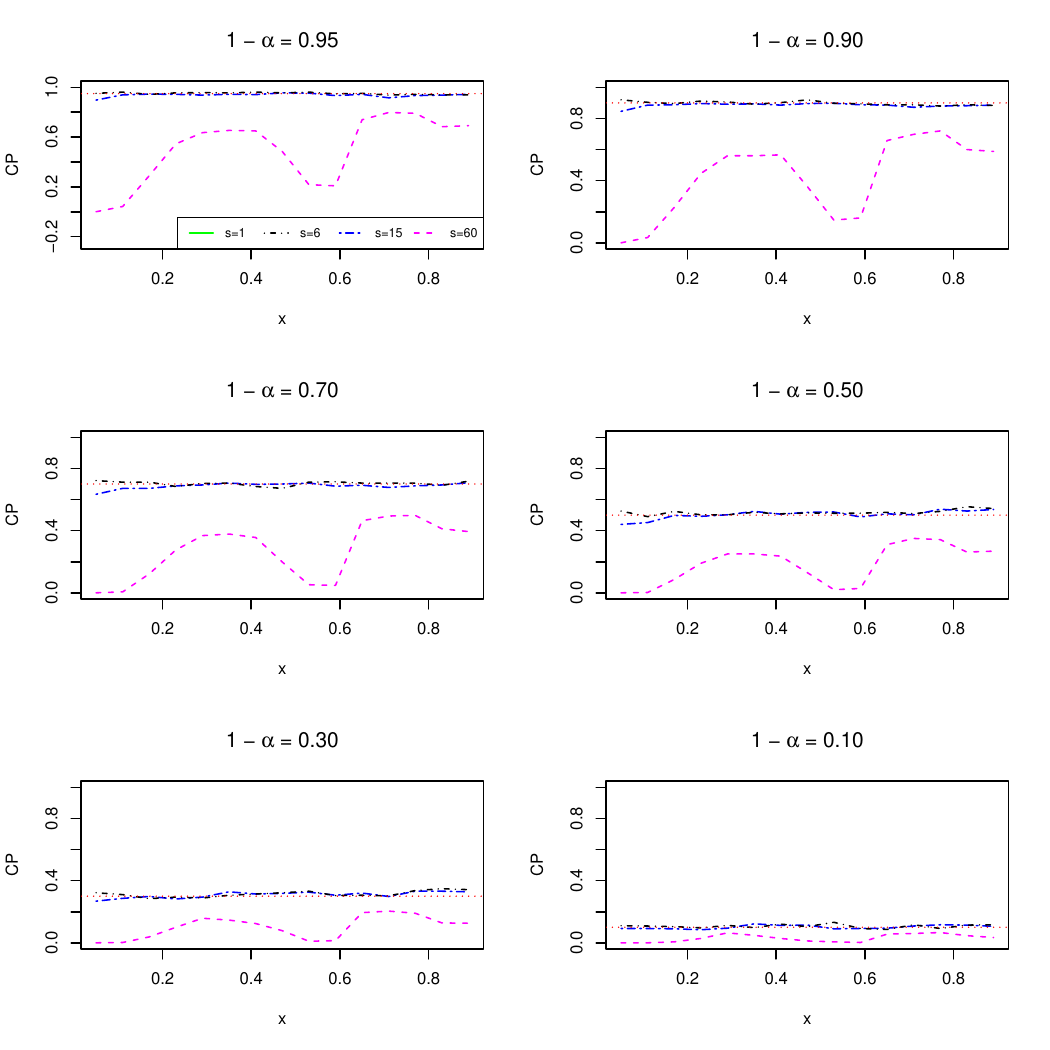}

\caption{\textit{\footnotesize $N=2400$: CP of $F_x(f)=\int_0^xf(z)dz$ against $x$ based on asymptotic theory.}}

\label{fig:asymp:integral:CI:2400}

\end{center}

\end{figure}

\end{document}